\documentclass{amsart}

\usepackage{amsmath,amssymb,amsthm,bbm} %
\usepackage{hyperref} %
\usepackage{cleveref} %
\usepackage{geometry} %
\usepackage{xcolor,enumerate} %
\usepackage{tikz} %
\usepackage{xfp} %
\usepackage{bm,todonotes}
\usepackage{csquotes}
\usepackage{caption}
\definecolor{theme}{RGB}{15, 87, 24} %
\definecolor{lighttheme}{RGB}{51, 153, 63} %
\definecolor{block}{RGB}{102, 60, 0} %
\definecolor{alert}{RGB}{212, 43, 43} %

\def\pagespace{3cm}
\hypersetup{
    colorlinks=true,
    linkcolor=theme,
    citecolor=lighttheme,
    urlcolor=lighttheme
}

\newtheorem{thm}{Theorem}[section]
\newtheorem{prop}[thm]{Proposition}
\newtheorem{lemma}[thm]{Lemma}
\newtheorem{cor}[thm]{Corollary}
\newtheorem{defn}[thm]{Definition}

\newtheorem{conj}[thm]{Conjecture}

\theoremstyle{remark}
\newtheorem{remark}[thm]{Remark}
\newtheorem{example}[thm]{Example}

\newcommand{\bR}{\mathbb{R}}
\newcommand{\Sym}{\mathfrak{S}}
\newcommand{\prob}{\mathbb{P}}
\newcommand{\expec}{\mathbb{E}}

\newcommand{\eps}{\varepsilon}
\DeclareMathOperator{\Leb}{Leb}
\DeclareMathOperator{\Poisson}{Poisson}

\DeclareMathOperator{\pat}{pat}
\DeclareMathOperator{\Perm}{Perm}

\DeclareMathOperator{\maj}{maj}

\newcommand\One{\mathbf{1}}

\newcommand{\LIS}[1]{\mathrm{LIS}\left(#1\right)}

\DeclareMathOperator{\Span}{Span}
\DeclareMathOperator{\Int}{Int}

\usepackage{listings}

\lstdefinelanguage{Wolfram}{
	morekeywords={ClearAll,Module,If,Return,Message,Length,Max,Round,Min,
		PermutationListQ,Subdivide,BinCounts,Transpose,Range,Clip,Sqrt,Pi,
		Ceiling,GaussianFilter,Padding,ConstantArray,Dimensions,Quantile,
		Select,Flatten,Framed,Image,Reverse,ColorSpace,ImageSize,FrameStyle,
		Directive,Black,AbsoluteThickness,FrameMargins,ContentPadding,
		Background,White,True,False},
	sensitive=true,
	morecomment=[s]{(*}{*)},
	morestring=[b]"
}

\title{{Permutons}}

\author[J. Borga]{Jacopo Borga}
\address[JB]{MIT, United States}

\author[V. Féray]{Valentin Féray}
\address[VF]{Université de Lorraine, CNRS, IECL, F-54000 Nancy, France}

\author[P. Winkler]{Peter Winkler}
\address[PW]{Dartmouth College, United States}

\date{}

\begin{document}

\vspace*{-1.6cm}

\maketitle

\vspace{-1cm}

\begin{abstract}
    Permutons are limit structures for large permutations.  Intrinsically, a permuton is a familiar mathematical object: a probability {measure} on the unit square with uniform marginals.  The compact space of permutons completes the set of permutations of all {sizes} under a natural topology; each permuton is not only a limit of permutations but also a generator of random permutations.
    
    Since their (relatively recent) discovery as a tool for understanding the behavior of large (possibly random) permutations, permutons have proved remarkably fruitful, with dozens of fascinating connections to combinatorics and probability.  We review the fundamentals of permutons, sometimes providing new, shorter proofs, and give an account of some of their myriad applications. The survey concludes with a long list of open problems. An extensive---but not comprehensive---bibliography is appended. 
\end{abstract}

\vspace{-0.5cm}

\tableofcontents

\section{Introduction}

\begin{figure}[t]
\[
\includegraphics[scale=0.7]{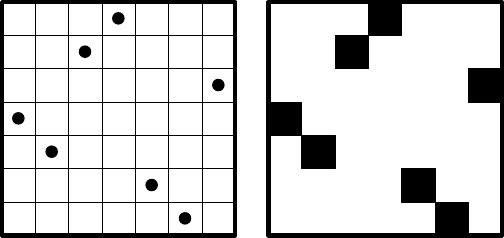}
\hspace{0.2cm}
\includegraphics[width=0.18\textwidth]{Mallows_100k}
\hspace{0.2cm}
\includegraphics[width=0.18\textwidth]{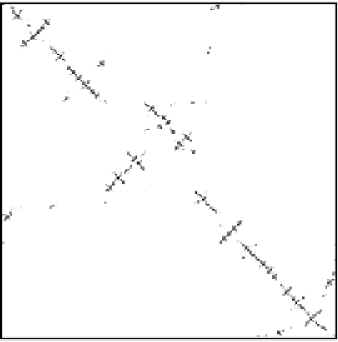}
\hspace{0.2cm}
\includegraphics[width=0.18\textwidth]{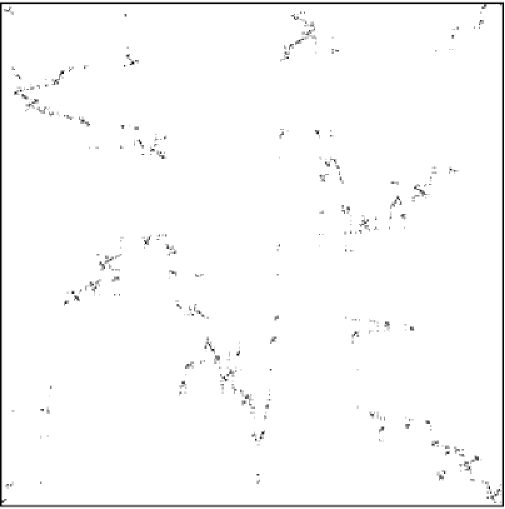}
\]
\caption{\textbf{From left to right}: (1) The diagram of the permutation $\sigma=4367215$; (2) the corresponding permuton where the black boxes represent Lebesgue measure;  (3) the permuton associated with a large Mallows permutation; (4) the permuton associated with a large uniform separable permutation; (5) the permuton associated with a large uniform Baxter permutation. See~\Cref{sect-pictures} for an explanation of how permutons associated with large permutations are represented.
\label{fig:fig-intro}}
\end{figure}

A permutation $\sigma$ of the numbers $\{1,2,\dots,n\}$ can be represented by its diagram: the set of points (see the first panel in~\Cref{fig:fig-intro})
\[
    \{(i,\sigma(i)):1\le i\le n\}.
\]
After rescaling the two coordinates by $n$, this diagram becomes a finite subset of the unit square.  The basic idea of \textbf{permutons} is to replace this finite picture with a probability measure on $[0,1]^2$. That is, one associates with a permutation $\sigma\in\Sym_n$ the probability measure obtained by putting Lebesgue measure of total  mass $1/n$ on each subsquare of $[0,1]^2$ corresponding to $(i/n,\sigma(i)/n)$; see the second panel in~\Cref{fig:fig-intro} and \eqref{eq:defn-perm-ass} for a precise definition.  In the limit $n\to\infty$, such measures may converge to probability measures on the square whose two marginals are uniform.  These limiting measures are called permutons.

This simple change of viewpoint has proved remarkably useful.  It turns permutation diagrams into probability measures and, in doing so, makes compactness available.  Questions about large permutations can then be rephrased as questions about a compact space of continuum objects.  For example, the convergence of permutations is equivalent to the convergence of all pattern densities, where the density of a pattern $\tau\in\Sym_k$ is the proportion of subsets of entries of $\sigma$ that appear in the same relative order as $\tau$.  Conversely, every permuton can be approximated by finite permutations, for instance, by sampling independent points from the permuton and reading their relative order. Thus, permutons play two complementary roles: they are limiting objects for large permutations, and they are also models from which random finite permutations can be sampled.

The analogy with graph limits is useful but should not be overplayed. {Permutons were introduced as permutation limits, in analogy with graphons, by Hoppen, Kohayakawa, Moreira, Ráth and Sampaio \cite{hoppen-permuton}; the term “permuton” was introduced later by Glebov, Grzesik, Klimošová and Král’~\cite{glebov15-forcible}.}  As in graphon theory, the limiting object records the asymptotic densities of all finite substructures; here, these substructures are permutation patterns.  At the same time, permutons have their own features.  In statistics, they are closely related to copulas: a permuton is precisely a probability measure on the unit square with uniform marginals, and its distribution function is a two-dimensional copula; see, for instance, \cite{gruebel2024copulas}.  In combinatorics, permutons provide a natural compactification of the set of permutations.  In probability, they allow one to study random permutations through weak convergence of measures, entropy, variational principles, and even continuum random geometry.

The purpose of this survey is to give an account of this developing theory. We have chosen a broad perspective, emphasizing the interaction between asymptotic combinatorics, probability, statistics, optimization, and random geometry. The survey is necessarily biased: many important topics are treated briefly, some are only mentioned in passing, and probably others have been entirely forgotten---we deeply apologize for that. When we believe it helps illuminate the main ideas, we also give new proofs of known results, often providing shorter or more direct arguments. Our goal is not only to review known results, but also to convey why permutons have become a useful language and why the area still contains many natural open problems.

\subsection{Permutons before permutons}\label{ssec:before}

Although the formal theory of permutation limits is quite recent, the idea of representing a large permutation by a measure appeared implicitly in several earlier contexts.  It is instructive to recall these precursors because they show that the same continuum object was reached from rather different directions.

One of the earliest appearances comes from \emph{non-parametric statistics}.
Let
\[
    (X_i,Y_i)_{1\le i\le n}
\]
be i.i.d.\ copies of a bivariate random vector $(X,Y)$ with a.s.\ no ties among the $X_i$'s and a.s.\ no ties among the $Y_i$'s.  If we sort the points by their $x$-coordinates and then record the relative order of their $y$-coordinates, we obtain a \emph{random} permutation $\sigma_n$.  When $X$ and $Y$ are independent, this permutation is uniform on $\Sym_n$.  Testing independence is therefore closely related to testing whether the induced random permutation looks asymptotically uniform, or, in modern language, whether it converges to the uniform permuton---and if this is the case, we say that the sequence $(\sigma_n)_n$ is quasirandom.  Hoeffding's work on non-parametric tests \cite{hoeffding1948test} already contains this philosophy: certain finite-dimensional rank statistics (of size 5) can force independence.  Later, Yanagimoto \cite{yanagimoto1970quasi_randomness} and then Král'--Pikhurko \cite{kral2013quasirandom} showed, in different language, that quasirandomness of permutations can be forced by densities of patterns of size $4$.  Only the latter paper uses the formal language of permutons, but the common idea is the same: finite pattern information can determine the limiting measure. See Section~\ref{sect:quasirandomness} for further discussion on quasirandomness.

A second source comes from \emph{probability theory}.  Deuschel and Zeitouni studied variational limiting curves for increasing records and longest increasing subsequences in i.i.d.\ planar samples~\cite{deuschel1995limiting_curves,deuschel1999increasing_subseq}.  If the common distribution of the sample has continuous marginals, applying the two marginal distribution functions to the two coordinates sends the sample to the unit square and produces a measure with uniform marginals.  This transformation does not change the relative order of the points. Hence, the induced permutation is exactly a sample from the associated permuton.  From the present viewpoint, their work can be seen as an early study of classical permutation statistics for permutations sampled from a continuum measure.

\begin{figure}[t]
	\[
	\includegraphics[scale=0.395]{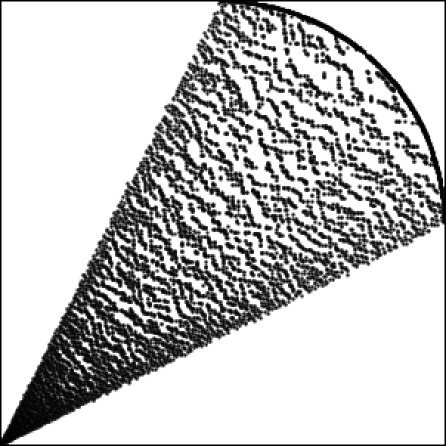}
	\hspace{0.2cm}
	\includegraphics[scale=0.35]{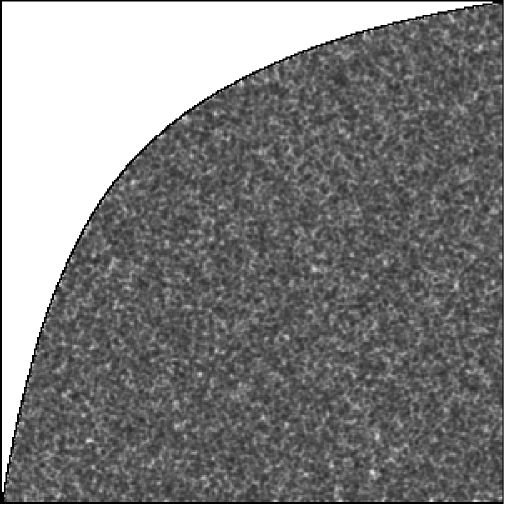}
	\hspace{0.2cm}
	\includegraphics[scale=0.465]{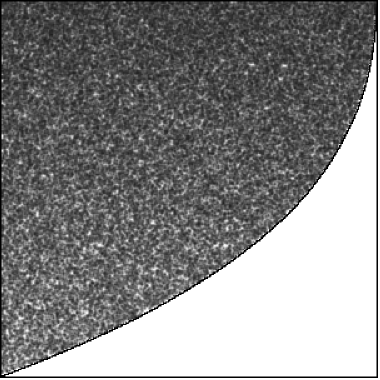}
	\hspace{0.2cm}
	\includegraphics[scale=0.35]{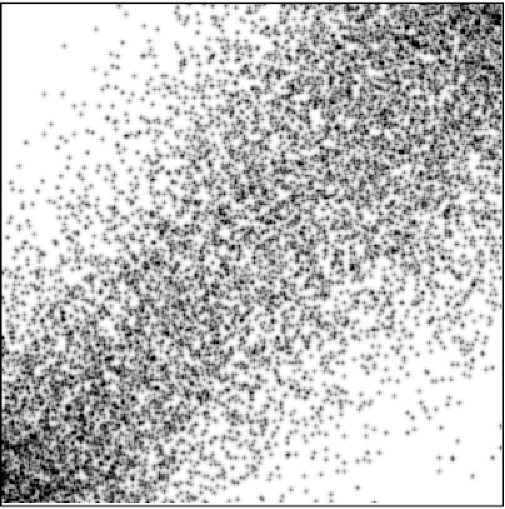}
	\hspace{0.2cm}
	\includegraphics[scale=0.465]{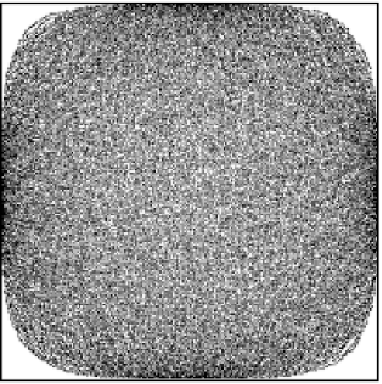}
	\]
	\caption{\textbf{From left to right}: {Permuton associated with a large} (1) {Grothendieck random permutation obtained from a uniformly random pipe dream}; (2) {record-biased random permutation};  (3) {uniform random permutation after applying the \emph{runsort} algorithm}; (4) {permutation sampled from the Spearman rank-correlation model}; (5) {uniform Erd\H{o}s--Szekeres extremal permutation.}
		\label{fig:fig-intro-2}}
\end{figure}

A third family of examples comes from what one may call \emph{algebraic probability theory}.  Some natural probability measures on algebraic or enumerative combinatorial objects produce random permutations whose diagrams have nontrivial macroscopic shapes.  Gnedin and Kerov~\cite{gnedin2002fibonacci}, in their study of the Young--Fibonacci graph, showed that the points of certain random involutions become uniformly spread in the square.  A more structured example arises from the Robinson--Schensted correspondence.  A permutation of size $n^2$ whose longest increasing and decreasing subsequences both have length $n$ is extremal for the Erd\H{o}s--Szekeres theorem and corresponds, under the Robinson--Schensted bijection, to a pair of square $n\times n$ standard Young tableaux.  Pittel and Romik proved a limit shape theorem for random square Young tableaux, and Romik used it to describe the typical diagram of such Erd\H{o}s--Szekeres extremal permutations \cite{pittel2007limit_shapes_square_tableaux,MR2266895}; see the fifth panel in~\Cref{fig:fig-intro-2}. This is an early example where a global constraint on the Robinson--Schensted shape produces a deterministic macroscopic geometry for the permutation diagram. Section~\ref{sect:algebraic-and-integrable} contains more examples of algebraic and integrable-probability limits.

Random sorting networks provide another important precursor.  A sorting network is a shortest path from the identity permutation to the reverse permutation using adjacent transpositions.  Angel, Holroyd, Romik and Virág \cite{angel2007sorting_networks} considered the permutations obtained at a fixed fraction of the time along a uniform random sorting network and gave conjectural formulas for their limiting shapes.  These limits are naturally permutons, and the conjectures were later proved by Dauvergne \cite{dauvergne2018archimedean}; {see Figure~\ref{fig:sorting-process}.}  In fact, sorting networks are best viewed not only as giving individual limiting permutons at fixed times, but also as giving a permuton-valued process. {Speaking of permuton processes, these objects
had also appeared in disguise in a completely different context
well before permutons themselves were introduced.
Indeed, already in 1989, Brenier~\cite{brenier1989least} studied 
{\em generalized incompressible flows}, which are nothing else than permutons evolving with time, where the uniformity of one-time marginals expresses incompressibility. 
He found an explicit sine-curve expression for such flows with minimal energy and given boundary conditions, which remarkably coincides, after rescaling, with the limiting process of random sorting networks; see Section~\ref{ssec:processes}.}

A fourth origin of permutons lies in \emph{extremal permutation theory}, in particular in the pattern packing problem.  Fix a pattern $\pi$.  Among permutations of size $n$, what is the largest possible density of occurrences of $\pi$, and what does this maximum converge to as $n\to\infty$?  Equivalently, what is the largest possible limiting density of $\pi$?  Presutti and Stromquist \cite{presutti2010packing} reformulated packing densities as optimization problems over probability measures on the square with uniform marginals and used this viewpoint to construct a lower bound for the packing density of the pattern $2413$.  Likewise, optimal constructions for several short patterns, including the layered optimizer in the $132$-packing problem, are most naturally understood as extremal permutons rather than as a single finite permutation \cite{price1997packing,albert2002packing}. These examples are conceptually important: permutons are not only limits of random models, but also optimizers of natural variational problems involving pattern densities.

Finally, \emph{large deviation theory} also encountered the same object before the terminology was standardized.  Trashorras \cite{trashorras2008large_deviations} proved a large deviation principle for the symmetrized empirical measure \[ \frac1n\sum_{i=1}^n \delta_{(X_i^n,X_{\sigma(i)}^n)}, \] where the $X_i^n$'s are random variables and $\sigma$ is a uniform random permutation.  For suitable choices of $X_i^n$, this empirical measure is essentially the permuton associated with $\sigma$.  The same large deviation principle was later rediscovered and developed further within the permuton literature.

\subsection{The modern theory and recent developments} The formal introduction of permutation limits by Hoppen, Kohayakawa, Moreira, Ráth and Sampaio \cite{hoppen-permuton} provided a unified language for these earlier ideas.  The starting point of the theory is the equivalence between two notions of convergence: a sequence of permutations converges as a sequence of measures on the unit square if and only if the density of every fixed pattern converges.  This makes the space of permutons a compact completion of the space of finite permutations.  Section~\ref{sec:theory} develops this basic theory. In particular, we revisit several of its foundational results and give new, often shorter or more direct proofs of most of the results presented there. We define permuton samples, prove that pattern densities determine the limit, and recall the algebraic structure generated by pattern density functions.

Once pattern densities are seen as continuous functions on a compact space, it becomes natural to ask which values they can take.  This leads to the feasible region of permutation pattern densities.  One-dimensional projections recover the packing problem; other projections encode quasirandomness questions.  The quasirandomness theorem of Král'--Pikhurko \cite{kral2013quasirandom}, the later refinements of Bergsma--Dassios \cite{bergsma2026quasirandomness_with_ties}, Chan--Král'--Noel--Pehova--Sharifzadeh--Volec \cite{kral2020quasirandom_pattern_sum}, and Crudele--Dukes--Noel \cite{crudele2023six_patterns_quasirandomness}, as well as the use of flag algebras in packing problems, all fit naturally into this framework.  We discuss these questions in Section~\ref{sect:feas-reg}.

Permutons also provide the right state space for large deviation and entropy methods.  For a uniform random permutation, the associated permuton converges to Lebesgue measure on the unit square, but deviations from this typical behavior occur at an exponential scale.  The rate function is the (negative) entropy of the limiting permuton.  This principle leads to variational problems for random permutations biased by pattern densities or by other statistics with continuous permuton analogues.  In particular, the Mallows model in the appropriate scaling regime can be understood by maximizing entropy plus an inversion-density term (see the third panel in~\Cref{fig:fig-intro}).  More generally, Mukherjee \cite{Mukherjee-estimation} and Kenyon--Král'--Radin--Winkler \cite{ Kenyon-Permutations-densities} showed how such variational problems control limit shapes and enumeration with prescribed pattern densities.  Section~\ref{sec:LDP} presents this circle of ideas.

A major source of recent activity has been the search for permuton limits of natural random permutation models.  Some limits are deterministic and analytic, such as the Mallows permuton or permutons arising from certain sequential sampling rules.  Others are \emph{random} and \emph{fractal}.  Uniform separable permutations (see the fourth panel in~\Cref{fig:fig-intro}) converge to the Brownian separable permuton \cite{bassino-separable-permuton}, and a much larger family of substitution-closed permutation classes has related Brownian separable limits \cite{Bassino-universal-perm}.  Baxter-type permutations (see the fifth panel in~\Cref{fig:fig-intro}) lead to the skew Brownian permutons \cite{ bm-baxter-permutation,borga-skew-permuton,borga2021permuton}.  Sorting networks lead to the Archimedean process \cite{dauvergne2018archimedean}, while models from algebraic combinatorics, such as Grothendieck random permutations (see the first panel in Figure~\ref{fig:fig-intro-2}), connect permuton limits with interacting particle systems and KPZ-type phenomena \cite{morales2025grothendieck}.  Section~\ref{sec:limit-shape} surveys these limit shape results.

The permuton topology is powerful, but it is deliberately macroscopic. It captures pattern densities, but it does not automatically control all classical permutation statistics.  Longest increasing subsequences, Robinson--Schensted shapes, fixed points, cycles, binary search trees, or the product of two permutations, depend on finer information than the weak limit of the permutation diagram. 
However, when the random permutations of interest have a nontrivial permuton limit,
a comparison with the samples of the limiting permuton often gives a good guess
for the limiting behavior of such statistics, and sometimes paves the way to a proof strategy by comparison arguments.
This is illustrated in Section~\ref{sec:samples}, 
and motivates the study of permuton samples.
Recent works study, among other things, their entropy \cite{maga2025entropy}, 
their Robinson--Schensted shapes \cite{sjostrand2023monotone}, 
the heights and shapes of their binary search trees \cite{corsini25-BST-samples}, 
and the minimal degrees in their inversion graphs~\cite{corsini26-min-inv-perm-sample}.

The final part of the survey focuses on a particularly rich family of random fractal permutons: the skew Brownian permutons.  They contain the Brownian separable permutons and the Baxter permuton as special cases, and they arise as universal limits for several families of constrained permutations.  They also admit a surprising connection with objects studied in \emph{random geometry}, such as Liouville quantum gravity surfaces and Schramm--Loewner evolution curves.  This connection reveals that a permuton limit can encode much more than a limiting density of points: in some cases, it remembers the structure of an entire random surface decorated by curves.  Section~\ref{sect:univ-perm-skew} is devoted to these objects, their construction, their pattern densities, their fractal supports, and the behavior of longest increasing subsequences sampled from them. We present this material in an accessible way, using only the intuition behind Liouville quantum gravity and Schramm--Loewner evolution, and assuming no prior knowledge of either theory.

We conclude in Section~\ref{sect:open-prob} with open problems.  Some are classical questions in a new language, such as the structure of feasible regions or the uniqueness of entropy-maximizing permutons.  Others are motivated by recent developments, such as proving the conjectural meandric permuton limit or understanding longest increasing subsequence exponents for general skew Brownian permutons.  We also briefly point to companion theories inspired by permutons, including limit theories for Latin squares~\cite{garbe2020limits}, local limits for permutations \cite{borga2020local} with their corresponding feasible region of consecutive patterns~\cite{borga2020feasible}, higher-dimensional permutations \cite{borga2024high_dimensional}, and more general integer tables \cite{tablons}.  All these new directions suggest that the theory of permutons is still young, and that the compact measure-theoretic viewpoint introduced by permutons will continue to interact fruitfully with many areas of combinatorics and probability. Of course, these developments raise many natural questions, and we believe that this is only the beginning of the story.

\section{The theory of permutons}\label{sec:theory}
Formally, a \textbf{permuton} is a probability measure on the unit square $[0,1]^2$ with both its marginals uniform, meaning that, for any Lebesgue-measurable subset $A$ of $[0,1]$:
\begin{equation}\label{eq:uniform_marginals}
    \mu(A\times[0,1])=
    \mu([0,1]\times A)=\lambda(A)
\end{equation}
where $\lambda$ is the Lebesgue measure on $[0,1]$.
Equivalently, if $(X,Y)$ has distribution $\mu$, then $X$ and $Y$ are both uniform
random variables in $[0,1]$, but they need not be independent, and the permuton $\mu$ precisely encodes this dependence.

The space $\mathcal P$ of permutons is equipped with the topology 
of weak convergence of measures.
In particular $\mathcal P$ is a closed subset of the space $\mathcal M_1([0,1]^2)$
of probability measures on $[0,1]^2$.
Since $\mathcal M_1([0,1]^2)$ is compact\footnote{The compactness
is a consequence of Prohorov's theorem, see, e.g., \cite[Theorem 5.1]{billingsley1999convergence}.}
 and metrizable,
the space $\mathcal P$ is compact and metrizable as well.

{A useful metric on $\mathcal P$ is the \textbf{rectangular distance}, introduced in \cite[Section~4]{hoppen-permuton} and defined by
\begin{equation}\label{eq:rectangular_distance}
	d_{\mathrm{rect}}(\mu,\nu)
	=\sup_{\substack{I,J\subseteq[0,1]\\ I,J\text{ intervals}}}
	\big|\mu(I\times J)-\nu(I\times J)\big|.
\end{equation}
This metric induces the topology of weak convergence on $\mathcal P$; see \cite[Lemma~5.3]{hoppen-permuton}. Indeed, every rectangle has boundary of measure zero under any permuton, since its marginals are uniform. Thus weak convergence implies convergence of the mass of each fixed rectangle. Approximating arbitrary rectangles by rectangles whose endpoints lie on a fixed finite grid, and using the uniform marginals to bound the mass of the intervening strips, makes this convergence uniform over all rectangles. Conversely, convergence in $d_{\mathrm{rect}}$ implies weak convergence by approximating continuous functions by step functions on a finite rectangular grid. 

In this survey, we mainly use the language of weak convergence to apply standard probabilistic tools, while the rectangular distance is particularly useful for quantitative estimates; see for instance Remark~\ref{rem:quantitative_convergence}.}

{The rectangular distance in~\eqref{eq:rectangular_distance} is analogous to the cut distance in graphon theory, with two important differences: graphon cut norms use arbitrary measurable subsets rather than intervals, and the cut distance also accounts for measure-preserving relabelings; see also Section~\ref{sect:open-prob-1} for further discussion.}

\subsection{Samples, moments and a criterion for convergence}
\label{ssec:samples_moments_convergence}
\begin{figure}[t]
\[\includegraphics{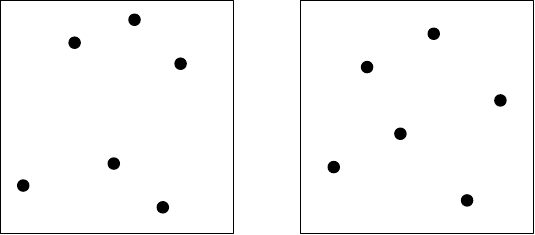}\]
\caption{\textbf{Left:} A set of points in $[0,1]^2$.
\textbf{Right:} The set of points $\{(i,\sigma(i))\}$,
for $\sigma=253614$. Both are in the same relative order.
\label{fig:sample}}
\end{figure}

Given a permuton $\mu$, we can associate, for each $n$, a random permutation
 $\sigma_n^\mu$ of size $n$ as follows.
 Let $(X_1,Y_1),\ldots,(X_n,Y_n)$ be a sequence of $n$ i.i.d.~random points\footnote{The pairs $(X_1,Y_1)$,\ldots,$(X_n,Y_n)$
 are i.i.d.~but, for a given $i$, the real numbers $X_i$ and $Y_i$ are in general not independent.} in  $[0,1]^2$ 
 with distribution $\mu$.
 We let $\pi$ and $\rho$ be the almost surely unique permutations such that
 \[X_{\pi(1)} < \dots < X_{\pi(n)}\quad\text{and}\quad Y_{\rho(1)} < \dots < Y_{\rho(n)}.\]
 Then we define $\sigma_n^\mu= \rho^{-1} \circ \pi$.
 Geometrically, $\sigma=\sigma_n^\mu$ is the unique
 permutation such that the sets of points $\{(X_1,Y_1),\ldots,(X_n,Y_n)\}$
and $\{(i,\sigma(i)), i \le n\}$ are in the same \enquote{relative order;} see \Cref{fig:sample}.
We say that the points $(X_1,Y_1),\ldots,(X_n,Y_n)$ form the permutation $\sigma$,
and
we will refer to $\sigma_n^\mu$ as a \textbf{permuton sample};
it is the random permutation formed from a sample of size $n$ of the measure $\mu$.

Let, in addition, $\tau$ be a permutation of size $k$. We write
\begin{equation}\label{eq:def_dtau_permuton}
d_\tau(\mu) = \mathbb P \big[ \sigma_k^\mu=\tau \big],
\end{equation}
and call $d_\tau$ the \textbf{density} of $\tau$ in $\mu$. 
Note that the permutation $\tau$ can be obtained from $k!$ pairs $(\rho,\pi)$ in
the above construction.
 By exchangeability of the sample $(X_1,Y_1),\ldots,(X_k,Y_k)$, all these pairs
 occur with the same probability.
 This means that, for any permutations $\pi$ and $\rho$ of size $k$, we have
\begin{equation}\label{eq:proba_pi_and_rho}
\mathbb P \big[ X_{\pi(1)} < \dots < X_{\pi(k)}, Y_{\rho(1)} < \dots < Y_{\rho(k)} \big]
= \tfrac{1}{k!} \, d_{\rho^{-1} \circ \pi} (\mu).
\end{equation}

These density functions are central elements in the theory of permutons.
We start with a simple result. We denote by $\Sym=\biguplus_{k \ge 1} \Sym_k$ the set of all permutations of all sizes.
\begin{prop}\label{prop:densities_separate}
The functions $d_\tau$ for $\tau \in \Sym$
separate points in $\mathcal P$, i.e.~any permuton
$\mu$ is determined by the values $(d_\tau(\mu))_{\tau \in \Sym}$.
\end{prop}
\begin{proof}
We will show that the moments of the measure $\mu$ can be expressed
in terms of density functions.
To do that, let us fix nonnegative integers $r,s \ge 0$,
and consider
 $r+s+1$ i.i.d.\ random points on $[0,1]^2$,
 \[P_0, \dots, P_{r+s},\] 
 with distribution $\mu$. We write $P_i=(x_{P_i},y_{P_i})$.
We consider the event
\[E= \left\{ \begin{array}{l}
x_{P_i} \le x_{P_0} \text{ for } 1 \le i \le r;\\
y_{P_i} \le y_{P_0} \,\text{ for } r+1 \le i \le r+s.
\end{array}\right\},\]
In words,
we want that $P_0$ has a higher $x$-coordinate than the next $r$ points in the sample,
and a higher $y$-coordinate than the last $s$ ones.
We will use a kind of double-counting argument, 
by computing $\prob(E)$ in two different ways.

A first method is to condition with respect to $P_0$.
Since $\mu$ has uniform marginals, for each $i \ge 1$,
\[ \prob\big[ x_{P_i} \le x_{P_0} \, \big| \, P_0 \big]=
\prob\big[ P_i \in [0, x_{P_0}] \times[0, 1] \big| P_0 \big]
= \mu\big( [0, x_{P_0}] \times[0, 1]\big) = x_{P_0},\]
and likewise,  $\prob\big[ y_{P_i} \le y_{P_0} \, \big| \, P_0 \big] =y_{P_0}$.
Since the $P_i$ are independent,
\[ \prob [E]= \expec \big[\expec \big[ \mathbf{1}_E | P_0 \big]\big] = \expec \big[x_{P_0}^r y_{P_0}^s \big]
=\int_{[0,1]^2} x^r y^s \mu(dx,dy),\]
where, in the last equation, we used the fact that $P_0$ has distribution $\mu$.
We note that $\prob [E]$ is nothing else than the joint moment of $\mu$ of order $(r,s)$.

The second method to compute $\prob [E]$ consists of looking at the (random)
permutations $\pi$ and $\rho$ of $\{0,\dots,r+s\}$ such that
 \[X_{\pi(0)} < \dots < X_{\pi(r+s)} \quad\text{and} \quad Y_{\rho(0)} < \dots < Y_{\rho(r+s)}.\]
The event $E$ is satisfied if and only if $\pi^{-1}(0)$ is greater than $\pi^{-1}(1)$,
\dots, $\pi^{-1}(r)$, and $\rho^{-1}(0)$ is greater than $\rho^{-1}(r+1)$,
\dots, $\rho^{-1}(r+s)$. Thus, using \eqref{eq:proba_pi_and_rho}, we obtain that
\[ \prob [E]= \sum_{\pi,\rho}  \frac{1}{(r+s+1)!} \, d_{\rho^{-1} \circ \pi} (\mu),\]
where the sum runs over permutations $\pi$ and $\rho$ satisfying the above conditions.

Hence, we see that all joint moments of $\mu$ can be expressed
as linear combinations of density functions. In particular, if two permutons $\mu$ and $\nu$
have the same evaluations by all density functions,
then they have the same joint moments and thus, they are equal
(measures on the compact set $[0,1]^2$ are determined by their moments thanks to the Stone--Weierstrass theorem and the Riesz representation theorem).
\end{proof}

This connection between density functions and moments can be further
used to prove that density functions not only separate points in $\mathcal P$,
but also generate the topology of $\mathcal P$.
The following theorem appeared first in \cite{hoppen-permuton}
(these authors use a different but equivalent notion of permutons, and simply call
them permutation limits) but the implication $(i) \Rightarrow (ii)$
had been proved earlier by Presutti and Stromquist \cite[Lemma 3.1]{presutti2010packing}.
\begin{thm}\label{thm:cv_permutons_densities}
Let $\mu_n$ be a sequence of permutons.
The following are equivalent:
\begin{enumerate}[(i)]
\item there exists a permuton $\mu$ such that $\mu_n$ converges to $\mu$;
\item for every $\tau$ in $\Sym$, there exists $c_\tau$ in $[0,1]$ 
such that $d_{\tau}(\mu_n)$ converges to $c_\tau$.
\end{enumerate}
Moreover, if (i) and (ii) hold, then, for every $\tau$ in $\Sym$, we have $d_{\tau}(\mu)=c_{\tau}$.
\end{thm}
\begin{proof}
Assume (i) and fix some $\tau$ in $\Sym_k$. 
For $n \ge 1$, we let $(X^n_1,Y^n_1),\ldots,(X^n_k,Y^n_k)$ be an i.i.d.~sample of law $\mu_n$;
we also let $(X_1,Y_1),\ldots,(X_k,Y_k)$ be an i.i.d.~sample of law $\mu$.
Since $\mu_n$ converges to $\mu$, the $k$-tuple $(X^n_1,Y^n_1),\ldots,(X^n_k,Y^n_k)$
converges in distribution to $(X_1,Y_1),\ldots,(X_k,Y_k)$ in $[0,1]^{2k}$.
We let $f_\tau$ be the following function on $[0,1]^{2k}$: 
\[f_\tau\big((x_1,y_1), \dots, (x_k,y_k)\big)= \begin{cases}
1 & \text{ if $(x_1,y_1), \dots, (x_k,y_k)$ induce the pattern $\tau$,}\\
0 & \text{ otherwise.} \end{cases}
\]
We note that $f_\tau$ is continuous on  $[0,1]^{2k} \setminus D$, where 
\[D= \{((x_1,y_1), \dots, (x_k,y_k)):  \exists i \ne j \text{ with }x_i=x_j\text{ or }y_i=y_j\} \]
is the set of points where at least two of the $x$-coordinates or two of the $y$-coordinates coincide.
Since $\mu$ is a permuton, the $X_i$'s (resp.~the $Y_i$'s) are i.i.d.~uniform random variables in $[0,1]$ and 
\[ \prob\Big[ \big((X_1,Y_1),\ldots,(X_k,Y_k)\big) \in D \Big] = 0.\]
By the mapping theorem (see, e.g., \cite[Theorem 29.2]{billingsley2012probability_and_measure}), 
 $f_\tau\big((X^n_1,Y^n_1),\ldots,(X^n_k,Y^n_k)\big)$
 converges in distribution to 
 $f_\tau\big((X_1,Y_1),\ldots,(X_k,Y_k)\big)$.
 Since $f_\tau$ is bounded, we also have convergence of expectation.
 But, by construction,
 \[\expec\Big[f_\tau\big((X^n_1,Y^n_1),\ldots,(X^n_k,Y^n_k)\big) \Big] =
 d_{\tau}(\mu_n),\]
 so (ii) is proved with 
 \[c_{\tau}=\expec\Big[f_\tau\big((X_1,Y_1),\ldots,(X_k,Y_k)\big) \Big] 
 =  d_{\tau}(\mu).\]
 Thus we have proved that (i) implies (ii), as well as the additional statement $d_{\tau}(\mu)=c_{\tau}$.
 
  Conversely, let us assume (ii). We recall from the proof of \Cref{prop:densities_separate}
 that moments of permutons $\nu$ can be expressed as linear combinations
 of density functions. Hence, (ii) implies that for any integers $r,s \ge 0$, there exists $c_{r,s}$
 such that 
 \begin{equation}\label{eq:conv_moments_mun}
  \lim_{n \to +\infty} \int_{[0,1]^2} x^r y^s \mu_n(dx,dy) =c_{r,s}.
  \end{equation}
Recalling that $\mathcal P$ is compact, one can extract from $\mu_n$
a weakly convergent subsequence $\mu_{n_k}$. Letting $\mu$ be its limit,
we  have
\[ \int_{[0,1]^2} x^r y^s \mu(dx,dy) = \lim_{k \to +\infty} \int_{[0,1]^2} x^r y^s \mu_{n_k}(dx,dy) =c_{r,s}.\]
Indeed, since we are working on $[0,1]^2$, weak convergence
and moment convergence are equivalent (convergence of all mixed monomial moments implies convergence against all polynomials, and polynomials are uniformly dense in $C([0,1]^2)$ by the Stone--Weierstrass theorem).
Going back to \eqref{eq:conv_moments_mun}, 
we conclude that $\mu_n$ converges weakly to $\mu$, i.e.~(i) holds.
The theorem is proved.
\end{proof}

\begin{remark}\label{rmk:trans-perm}
Let $\mu$ be a probability measure on $\bR^2$ whose projections on the axes do not have atoms.
Then we can construct a random permutation $\sigma_k^\mu$ as above.
We will, however, argue that  $\sigma_k^\mu$ has the same distribution
as  $\sigma_k^{\tilde \mu}$, where $\tilde \mu$ is a permuton associated with $\mu$.
To do this, we let $F_1$ (resp.~$F_2$) be the cumulative distribution function of $x_P$ (resp.~$y_P$), where $P=(x_P,y_P)$ is a random point in $\bR^2$ with distribution $\mu$.
Then we set $\tilde P=(F_1(x_P),F_2(y_P))$,
and let $\tilde \mu$ be the distribution of $\tilde P$. 
It is straightforward to check that $\tilde \mu$ has uniform marginals,
and that $\sigma_k^\mu$ and $\sigma_k^{\tilde \mu}$ have the same law
(since an i.i.d.~sample with distribution $\tilde \mu$ is obtained from a sample
with distribution $\mu$ by applying nondecreasing functions to the $x$- and $y$-coordinates).
Consequently,
\begin{itemize}
\item there is no loss of generality in considering only permutons $\mu$ when studying 
the random permutation $\sigma_k^\mu$ ;
\item \Cref{prop:densities_separate} only holds for permutons: if $\mu$ is not a permuton,
then the measures
$\mu$ and $\tilde \mu$ are different, but all density functions have the same
value on $\mu$ and $\tilde \mu$.
\end{itemize}
\end{remark}

\subsection{Embedding permutations into permutons}
The set $\Sym$ of all permutations can be embedded into the space of permutons as follows:
with a permutation $\sigma\in \Sym_n$, we associate the permuton $\mu_\sigma$ with density 
\begin{equation}\label{eq:defn-perm-ass}
    \mu_\sigma(dxdy) = n\One_{\sigma(\lceil xn \rceil) = \lceil yn \rceil} dxdy.
\end{equation}
Note that it amounts to replacing every point $(i,\sigma(i))$ in the diagram of $\sigma$ (normalized to the unit square) 
by a square of the form $[(i-1)/n,i/n]\times [(\sigma(i)-1)/n,\sigma(i)/n]$, 
with a mass $1/n$ uniformly distributed on it. 
This allows us, in particular, to define \textbf{convergent sequences} of permutations:
we say that $(\sigma_n)_n$ converges to a permuton $\mu$ when $\mu_{\sigma_n} \to \mu$ weakly. 
\smallskip

On the other hand, there is a natural notion of pattern densities in permutations.
If $\tau$ and $\sigma$ are permutations of size $k$ and $n$ respectively, with $k \le n$,
we let $d_\tau(\sigma)$ be the proportion of patterns of size $k$ in $\sigma$ that are equal to $\tau$;
in formula,
\begin{equation}\label{eq:def_dtau_permutation}
d_\tau(\sigma)= \binom{n}{k}^{-1} \, \big| \big\{ I \subseteq [n]: |I|=k \text{ and }\pat_I(\sigma)=\tau \big\}\big|,
\end{equation}
where $\pat_{\bm I}(\sigma)$ denotes the pattern induced by the points $(i,\sigma(i))_{i \in I}$ and $[n]=\{1,2,\dots,n\}$.
This has a probabilistic interpretation: if $\bm I$ is a uniform random subset of size $k$ of $[n]$, then
\[d_\tau(\sigma)= \prob\big[\pat_{\bm I}(\sigma)=\tau\big].\]

\smallskip

We warn the reader that $d_\tau(\sigma) \ne d_\tau(\mu_\sigma)$ in general, 
but they are close to each other for fixed $\tau$ and large permutations $\sigma$, as ensured by the following lemma
(see \cite[Lemma 3.5]{hoppen-permuton}).
\begin{lemma}\label{lem:densityofassociatedpermuton}
  If $\tau \in \Sym_k$ and $\sigma \in \Sym_n$, then 
	$$
	|d_\tau(\sigma)-d_{\tau}(\mu_\sigma)| \leq \frac 1 n \binom k 2.
	$$
\end{lemma}
\begin{proof}
  Let $(X_1,Y_1),\dots,(X_k,Y_k)$ be an i.i.d.~sample of distribution $\mu_\sigma$
  and let us call $\pi_k$ the random permutation formed by these $k$ points.
  For $1 \le j \le k$, set $I_j=\lceil nX_j \rceil$. We note that $I_1,\dots,I_k$ are i.i.d.~random variables in $\{1,\dots,n\}$.
  Let now $E=\{ \exists\, j < j': I_j=I_{j'} \}$ be the event that two of the $I_j$'s coincide.
  Using the union bound, we get
  \[\prob[E] \le \sum_{j < j'} \prob[I_{j}=I_{j'}] = n^{-1}\, \binom{k}2.\]
  On the other hand, conditionally on the complement $\bar{E}$ of $E$,
  the random set $\bm I=\{I_1,\dots,I_k\}$ is a uniform random subset of $\{1,\dots,n\}$
  and  the permutation $\pi_k$ formed by the points $(X_1,Y_1),\dots,(X_k,Y_k)$
  is the pattern $\pat_{\bm I}(\sigma)$ induced by $\sigma$ on the set $\bm I$ of positions.
  This implies
  \[ d_\tau(\sigma) = \prob\big[\pat_{\bm I}(\sigma)=\tau \big| \bar{E} \big].\]
  We now write
  \[ d_{\tau}(\mu_\sigma) = \prob\big[ \pi_k=\tau, E] + \prob\big[ \pi_k=\tau, \bar{E}],\]
  and observe that $\prob\big[ \pi_k=\tau, E] \le \prob[E]$ and
  \[\prob\big[ \pi_k=\tau, \bar{E}] = \prob[\bar{E}] \cdot \prob\big[\pat_{\bm I}(\sigma)=\tau \big| \bar{E} \big]= (1-\prob[E]) d_\tau(\sigma).\]
  We finally find
  \[ d_\tau(\sigma)- \prob[E] \le (1-\prob[E]) d_\tau(\sigma) \le d_\tau(\mu_\sigma) \le (1-\prob[E]) d_\tau(\sigma) +\prob[E] \le d_\tau(\sigma) +\prob[E].\]
  Recalling the bound $\prob[E] \le n^{-1}\, \binom{k}2$, this proves the lemma.
\end{proof}
As a consequence, if one considers a sequence $\mu_n=\mu_{\sigma_n}$ of permutons associated with permutations of size tending to $+\infty$,
Theorem~\ref{thm:cv_permutons_densities} has an immediate analog, where we replace $d_\tau(\mu_n)=d_\tau(\mu_{\sigma_n})$
by $d_{\tau}(\sigma_n)$ in item (ii). 
In particular, a sequence $\sigma_n$ of permutations of size tending to $+\infty$ is convergent if and only if, for any pattern $\tau$,
the proportion $d_\tau(\sigma_n)$ of occurrences of $\tau$ in $\sigma_n$ converges.
\begin{remark}
In \cite{hoppen-permuton}, the convergence of $d_\tau(\sigma_n)$ for all $\tau$ is taken
as the definition of {\bf convergent sequences of permutations},
and permutons are found as convenient analytic objects representing the possible limits,
in the spirit of graphon theory.
\end{remark}

A natural question is whether the image of $\Sym$ in $\mathcal P$ is dense, i.e.~if every permuton
can be obtained as limit of a sequence of finite permutations. 
We will see that this is indeed the case, using again the notion of permuton samples 
introduced in \Cref{ssec:samples_moments_convergence}.
\begin{prop}\label{prop:convergence_samples}
  Let $\mu$ be a permuton. Then its permuton samples $\sigma_n^\mu$ converge
  almost surely to $\mu$, as $n$ tends to $+\infty$.
\end{prop}
\begin{proof}
  Let $\tau$ be a permutation in $\Sym_k$. We want to show that the random variable $d_\tau(\sigma_n^{\mu})$
  converges almost surely to the deterministic number $d_\tau(\mu)$.
  For this, we recall that $\sigma_n^{\mu}$ is the permutation associated with a sample $(X_1,Y_1),\dots,(X_n,Y_n)$ of law $\mu$.

  Using the probabilistic interpretation of $d_\tau$, we have
  \[ \expec\big[ d_\tau(\sigma_n^{\mu}) \big]= \prob \big[ \pat_{\bm I}(\sigma_n^{\mu})=\tau \big],\]
  where $\bm I$ is a uniform random subset of size $k$ of $[n]$, independent from $\sigma_n^{\mu}$.
  Taking $\pat_{\bm I}(\sigma_n^{\mu})$ amounts to looking at the permutation formed by a uniform random
  subsample of size $k$ of the sample $\{(X_1,Y_1),\dots,(X_n,Y_n)\}$.
  But such a subsample has the same law as an i.i.d.~sample of size $k$.
  \[ \expec\big[ d_\tau(\sigma_n^{\mu}) \big] = \prob\big[ \Perm( (X_1,Y_1),\dots,(X_k,Y_k) )=\tau \big]=d_\tau(\mu),\]
  where $\Perm( (X_1,Y_1),\dots,(X_k,Y_k) )$ is the permutation formed by an i.i.d.~sample $(X_1,Y_1),\dots,(X_k,Y_k)$ of size $k$.
  In particular, the expectation of $d_\tau(\sigma_n^{\mu})$ converges (and is even equal) to $d_\tau(\mu)$.

  To prove almost sure convergence of $d_\tau(\sigma_n^{\mu})$ to $d_\tau(\mu)$, we further need some concentration
  inequality of $d_\tau(\sigma_n^{\mu})$ around its expectation. This will be provided by McDiarmid's inequality.
  With a slight abuse of notation, we see $d_\tau=d_\tau(\sigma_n^{\mu})$ as a function of $n$ pairs $(x_1,y_1)$, \dots, $(x_n,y_n)$.
  Changing a pair $(x_i,y_i)$ by another pair $(x'_i,y'_i)$ then changes at most $\binom{n-1}{k{-}1}$ patterns of size $k$
  in the associated permutations, so that, for all $i$, $(x_1,y_1)$, \dots, $(x_n,y_n)$, $(x'_i,y'_i)$
  \[
  \Big| d_\tau\big( (x_1,y_1),\dots,(x_i,y_i),\dots,(x_n,y_n) \big) - d_\tau\big( (x_1,y_1),\dots,(x'_i,y'_i),\dots,(x_n,y_n) \big) \Big| \le \frac{\binom{n-1}{k{-}1}}{\binom{n}{k}} = \frac{k}{n}.
  \]
  {Using McDiarmid's inequality \cite[Theorem~3.1]{McDiarmid1998Concentration}, we obtain, for every $n\ge k$ and $t>0$,
  \begin{equation}\label{eq:deviation_dtau}
  		\prob\Big[
  		\big|d_\tau(\sigma_n^\mu)-d_\tau(\mu)\big|\ge t
  		\Big]
  		\le 2\exp\!\left(-\frac{2nt^2}{k^2}\right).
  \end{equation}
  Here we used $\expec[d_\tau(\sigma_n^\mu)]=d_\tau(\mu)$, as proved above. In particular, taking $t=n^{-1/4}$ gives 
  \[ 
  	\prob\Big[ \big|d_\tau(\sigma_n^\mu)-d_\tau(\mu)\big|\ge n^{-1/4} \Big] \le 2e^{-2k^{-2}n^{1/2}}. 
  \] 
  Since the right-hand side is summable in $n$, the Borel--Cantelli lemma implies that, almost surely, $\big|d_\tau(\sigma_n^\mu)-d_\tau(\mu)\big|<n^{-1/4}$ for all sufficiently large $n$. Thus $d_\tau(\sigma_n^\mu)$ converges almost surely to $d_\tau(\mu)$.}

  In the preceding argument, $\tau$ was a fixed permutation in $\Sym$. 
  However, since $\Sym$ has countably many elements, this implies
  that almost surely, the convergence $d_\tau(\sigma_n^{\mu}) \to d_\tau(\mu)$
  holds for all $\tau$ simultaneously.
  Using \Cref{thm:cv_permutons_densities} (and the discussion after \Cref{lem:densityofassociatedpermuton}),
  we conclude that $\sigma_n^{\mu}$ converges almost surely to $\mu$.
\end{proof}

\begin{cor}\label{cor:permutations_dense_in_P}
  The set $\{\mu_\sigma,\sigma \in \Sym\}$  is dense 
  in the space $\mathcal P$ of permutons.
\end{cor}

{
	\begin{remark}\label{rem:quantitative_convergence}
		Eq.~\eqref{eq:deviation_dtau} gives in some sense a rate of convergence
		for the convergence of $\sigma_n^\mu$ to $\mu$.
		Indeed it implies that for every fixed pattern $\tau$,
		we have
			$d_\tau(\sigma_n^\mu)-d_\tau(\mu)=O_{\mathbb P}(n^{-1/2}).$
		Other quantitative convergence results have appeared in the literature,
		using the rectangular distance
		$d_{\mathrm{rect}}$ of~\eqref{eq:rectangular_distance} rather than pattern densities, see for example \cite[Lemma~4.2]{hoppen-permuton}.
		
		Instead of considering the permuton sample $\sigma_n^\mu$,
			one can  optimize the choice of the approximating permutation
			and study
			\[
			e_n(\mu)=\min_{\pi\in\Sym_n}
			d_{\mathrm{rect}}(\mu_\pi,\mu).
			\]
			This quantity need not have the same order as 
			$d_{\mathrm{rect}}(\mu_{\sigma_n^\mu},\mu)$. Maga~\cite{maga2026approximation} studies this optimal
			approximation problem and shows how the geometry of $\mu$ can
			affect the approximation rate.
	\end{remark}
}

{Proposition~\ref{prop:convergence_samples} can also be completed 
by a central limit theorem.}
{
	\begin{prop}\label{prop:sample_pattern_clt}
		Let $\mu$ be a fixed permuton, and let $\tau_1,\ldots,\tau_r$
		be patterns of respective sizes $k_1,\ldots,k_r$,
		and let
		$h_a(z_1,\ldots,z_{k_a})$ be the indicator that the points $z_1,\ldots,z_{k_a}$
			induce $\tau_a$.
		For i.i.d.\ points $Z_1,Z_2,\ldots$ of law $\mu$, define
		\[
		g_a(z)=\expec[h_a(z,Z_2,\ldots,Z_{k_a})].
		\]
		Then
		\[
		\sqrt n\,
		\bigl(d_{\tau_a}(\sigma_n^\mu)-d_{\tau_a}(\mu)\bigr)_{1\le a\le r}
		\ \xrightarrow{\mathrm d}\ \mathcal N_r(0,\Sigma),
		\qquad
		\Sigma_{ab}
		=k_a k_b\operatorname{Cov}\bigl(g_a(Z_1),g_b(Z_1)\bigr).
		\]
		The limiting Gaussian vector is allowed to be degenerate.
	\end{prop}
}

{
	\begin{proof}
		For each $a$, write $k=k_a$, $d_{\tau}(\mu)=d_{\tau_a}(\mu)$, $h=h_a$, and $g=g_a$. Define the centered kernel \[ \widetilde h(z_1,\ldots,z_k) =h(z_1,\ldots,z_k)-d_{\tau}(\mu) -\sum_{j=1}^k(g(z_j)-d_{\tau}(\mu)). \] By the definition of $g$ and the symmetry of $h$, we have $\expec[\widetilde h(Z_1,\ldots,Z_k)\mid Z_j]=0$ for every $j$. Since $d_{\tau_a}(\sigma_n^\mu)$ is the average of $h$ over all $k$-subsets of the sample, we obtain \[ d_{\tau_a}(\sigma_n^\mu)-d_{\tau}(\mu) =\frac{k}{n}\sum_{i=1}^n(g(Z_i)-d_{\tau}(\mu))+R_{n,a}, \] where \[ R_{n,a}=\binom nk^{-1} \sum_{\substack{I\subseteq[n]\\|I|=k}} \widetilde h((Z_i)_{i\in I}). \] 
		Since $\mathbb{E}[\widetilde h(Z_1,\ldots,Z_k)]=0$, we have
		\begin{align*}
			\mathbb{E}[R_{n,a}^2]
			&= \binom{n}{k}^{-2}
			\sum_{\substack{I,J\subseteq[n]\\ |I|=|J|=k}}
			\mathbb{E}\left[
			\widetilde h((Z_i)_{i\in I})
			\widetilde h((Z_j)_{j\in J})
			\right].
		\end{align*}
		The summands vanish whenever $|I\cap J|\leq 1$:
		if $I\cap J=\varnothing$, this follows from independence;
		if $|I\cap J|=1$, it follows by conditioning on the
		common variable and using the conditional expectation
		identity above. Therefore,
		\begin{align*}
			\mathbb{E}[R_{n,a}^2]
			= \binom{n}{k}^{-2}
			\sum_{\substack{I,J\subseteq[n]\\
					|I|=|J|=k,\ |I\cap J|\geq 2}}
			\mathbb{E}\left[
			\widetilde h((Z_i)_{i\in I})
			\widetilde h((Z_j)_{j\in J})
			\right] 
			= O\left(\frac{n^{2k-2}}{n^{2k}}\right)
			= O(n^{-2}),
		\end{align*}
		since $\widetilde h$ is bounded and there are
		$O(n^{2k-2})$ pairs $(I,J)$ with $|I\cap J|\geq 2$.
		In particular, $\sqrt n R_{n,a}\to0$ in $L^2$. The multivariate central limit theorem applied to the i.i.d.\ vectors $(k_a(g_a(Z_i)-d_{\tau_a}(\mu)))_{1\le a\le r}$ now gives the result.
	\end{proof}
}

{For a single pattern $\tau$, the limit is nondegenerate precisely
	when $\operatorname{Var}(g_\tau(Z_1))>0$.
	{Concretely $\operatorname{Var}(g_\tau(Z_1))=0$ means that the probability
	that $Z_1, \dots,Z_k$ form the pattern $\tau$ conditionally on the first point $Z_1$
	does not depend on $Z_1$.
	The permuton $\mu$ is then called {\em conditionally constant} (CC) with respect to~$\tau$;
	interestingly, this notion appeared in~\cite{borga2022large}
	in some large deviation study, a priori not directly related to Proposition~\ref{prop:sample_pattern_clt}.
	Being CC is of course a peculiar situation, but there are nontrivial examples. In particular, all $\tau$-avoiding or $\tau$-packing
	permutons\footnote{A permuton $\mu$ is $\tau$-avoiding if $d_\tau(\mu)=0$,
	and $\tau$-packing if $d_\tau(\mu)$ has the maximal possible value over the whole space of permutons.}  are CC w.r.t. $\tau$~\cite[Remark 1.18]{borga2022large}.
	A complete classification of CC permutons seems to be a difficult problem:
	\cite[Proposition 1.31]{borga2022large} provides necessary conditions in the case $\tau=21$.}
	
	If this variance vanishes, {i.e.~if $\mu$ is CC w.r.t.~$\tau$,
	then} the $\sqrt n$-scaled fluctuations
	converge to zero; finer fluctuations, when present, can have a
	different scale and need not be Gaussian.
	{In the uniform case (i.e.~$\mu$ is the Lebesgue measure on $[0,1]^2$, and $\sigma_n^\mu$ is a uniform random permutation),
	the asymptotic variance $\operatorname{Var}(g_\tau(Z_1))$
	never vanishes, but the asymptotic covariance matrix for all patterns in $S_n$
	has only rank $(n-1)^2$ (while it has size $n!$), as shown in \cite{janson2015CLT_patterns}.}
		Quantitative refinements and mod-Gaussian convergence in the
	nondegenerate case are established in
	\cite[Theorem~26]{feray2020moudli_spaces}.
	{Note that} these statements concern the sampling model $\sigma_n^\mu$:
	convergence to $\mu$ alone does not determine the fluctuations.}

\subsection{{The algebra of pattern densities and convergence for random permutons}}

The following proposition can be obtained from Razborov's general theory of flag algebras \cite{razborov2007flag}, together with the results of the previous section relating permutation limits to permutons. {However, no familiarity with flag algebras is needed here, as we give a direct proof in the permuton setting.} See also \cite{vargas2014hopf_permutations2} for a non-asymptotic version of the corresponding algebra.

\begin{prop}\label{prop:densities_span_algebra}
The linear span $\mathcal O:=\Span(d_\tau, \tau \in \Sym)$ is a dense subalgebra of the space
$C(\mathcal P)$ of real-valued continuous functions on $\mathcal P$ {with respect to the sup-norm $\|f\|_{\infty}=\sup_{\mu\in \mathcal P}|f(\mu)|$}.
\end{prop}
\begin{proof}
We first prove that $\mathcal O$ is a subalgebra of $C(\mathcal P)$, i.e.~that any product $d_{\tau_1} d_{\tau_2}$ is a linear combination of $d_{\sigma}$'s.
Let $\mu$ be a permuton, and $\tau_1$ and $\tau_2$ be permutations of respective sizes
$k_1$ and $k_2$. 
The product $d_{\tau_1}(\mu) d_{\tau_2}(\mu)$ has the following probabilistic interpretation.
We take $(X_1,Y_1) \dots (X_{k_1+k_2},Y_{k_1+k_2})$ i.i.d.~random points with law $\mu$.
Fix permutations $\pi_1$, $\rho_1$, $\pi_2$ and $\rho_2$ such that $\rho_i^{-1} \circ \pi_i=\tau_i$ for $i$ in $\{1,2\}$. Then, from \eqref{eq:proba_pi_and_rho}, we have
\[ d_{\tau_1}(\mu) = k_1! \,\prob\big[X_{\pi_1(1)} < \dots < X_{\pi_1(k_1)}, Y_{\rho_1(1)} < \dots < Y_{\rho_1(k_1)} \big].\]
Similarly, if we interpret $\pi_2$ and $\rho_2$ as permutations of $\{k_1+1,\dots,k_1+k_2\}$
by shifting everything by $k_1$, we get
\[ d_{\tau_2}(\mu) = k_2! \,
\prob\big[X_{\pi_2(k_1+1)} < \dots < X_{\pi_2(k_1+k_2)}, Y_{\rho_2(k_1+1)} < \dots < Y_{\rho_2(k_1+k_2)} \big].\]
Since $(X_1,Y_1), \dots, (X_{k_1},Y_{k_1})$ is independent from $(X_{k_1+1},Y_{k_1+1}), \dots, (X_{k_1+k_2},Y_{k_1+k_2})$, we get
\[ d_{\tau_1}(\mu) d_{\tau_2}(\mu) = k_1!k_2! 
\, \prob\left[ \begin{array}{c} X_{\pi_1(1)} < \dots < X_{\pi_1(k_1)}, X_{\pi_2(k_1+1)} < \dots < X_{\pi_2(k_1+k_2)}, \\ Y_{\rho_1(1)} < \dots < Y_{\rho_1(k_1)}, Y_{\rho_2(k_1+1)} < \dots < Y_{\rho_2(k_1+k_2)} \end{array}\right]. \]
As for $\prob[E]$ in the proof of \Cref{prop:densities_separate}, 
one can prove that the above probability is a linear combination of $d_\tau$'s,
for $\tau$ in $\Sym_{k_1+k_2}$.
Thus $\mathcal O$ is a subalgebra of $C(\mathcal P)$.

To prove that $\mathcal O$ is dense, observe that:
\begin{itemize}
\item $\mathcal O$ is a subalgebra of $C(\mathcal P)$;
\item $\mathcal O$ contains the constant functions (indeed, if $\tau$ is the unique size $1$ permutation, then $d_{\tau}(\mu)=1$ for any permuton $\mu$);
\item $\mathcal O$ separates points in $\mathcal P$ (by \Cref{prop:densities_separate}).
\end{itemize}
Using the Stone-Weierstrass theorem, see, e.g.,~\cite{rao2005StoneWeierstrass}
(recall that $\mathcal P$ is compact), we conclude that $\mathcal O$
is dense in $C(\mathcal P)$.
\end{proof}

We can use this algebra structure, and in particular the density result,
 to obtain a convergence theorem for {\em random} permutons,
 generalizing \Cref{thm:cv_permutons_densities}.
This theorem first appeared in \cite{Bassino-universal-perm}, 
but the proof given here is different.

\begin{thm}\label{thm:cv_random_permutons}
Let $\mu_n$ be a sequence of random permutons. Then the following are equivalent:
\begin{enumerate}[(i)]
\item There exists a random permuton $\mu$ such that $\mu_n$
 converges to $\mu$ in distribution.
\item For any pattern $\tau$, there exists a $[0,1]$-valued random variable $X_\tau$
such that $d_{\tau}(\mu_n)$ converges to $X_\tau$ in distribution.
\item For any pattern $\tau$, there exists $c_\tau$ in $[0,1]$ such that 
$\expec[d_{\tau}(\mu_n)]$ converges to $c_\tau$.
\end{enumerate}
Moreover, if these conditions are fulfilled, we have
\[X_\tau \stackrel{\mathrm{d}}{=} d_{\tau}(\mu) \text{ and }c_\tau= \expec[X_\tau]=\expec[d_{\tau}(\mu)].\]
\end{thm}
Before proving the theorem, we make a couple of remarks.
First, this theorem can in particular be applied to a sequence $\mu_n=\mu_{\sigma_n}$ of permutons
associated with \emph{random} permutations $\sigma_n$, whose size tends to $+\infty$.
In this case, $d_{\tau}(\mu_n)$ can be replaced by $d_{\tau}(\sigma_n)$ in items (ii) and (iii) above
(recall \Cref{lem:densityofassociatedpermuton}).
Second, the theorem implies that 
the distribution of a random permuton $\mu$ is uniquely
determined by its expected pattern densities $\expec[d_{\tau}(\mu)]$, for $\tau$ in $\Sym$
(improving \Cref{prop:densities_separate}).
\begin{proof}[Proof of \Cref{thm:cv_random_permutons}]
By \Cref{thm:cv_permutons_densities}, for any $\tau$ in $\Sym$,
the function $d_\tau$ from $\mathcal P$ to $[0,1]$ is continuous.
Hence, (i) implies that $d_{\tau}(\mu_n)$ converges in distribution to $d_{\tau}(\mu)$,
proving (ii) and the additional statement 
$X_\tau \stackrel{\mathrm{d}}{=} d_{\tau}(\mu)$ (which itself implies $\expec[X_\tau]=\expec[d_{\tau}(\mu)]$).

Now assume (ii). Since the random variable $d_{\tau}(\mu_n)$ is bounded,
convergence in distribution implies convergence of the expectation.
Namely, $\expec[d_{\tau}(\mu_n)]$ converges to $\expec[X_{\tau}]$,
proving (iii) and the additional statement $c_\tau= \expec[X_{\tau}]$.

It remains to prove that (iii) implies (i).
By linearity, (iii) implies that $\expec\big[f(\mu_n)\big]$ converges for any $f$ in 
 $\mathcal O=\Span(d_\tau, \tau \in \Sym)$. Using the density of $\mathcal O$ in $C(\mathcal P)$ (\Cref{prop:densities_span_algebra}), we claim that it is also the case for any $f$ in $C(\mathcal P)$.
 Indeed, for any $f$ in $C(\mathcal P)$ and $g$ in 
 $\mathcal O$, we have
\[ \limsup_{n,m \to +\infty} \big| \expec\big[f(\mu_n)\big] - \expec\big[f(\mu_m)\big] \big|
\le  2 \|g-f\|_{\infty} + \limsup_{n,m \to +\infty} \big| \expec\big[g(\mu_n)\big] - \expec\big[g(\mu_m)\big] \big|.\]
For any $g$ in  $\mathcal O$, the second term is equal to $0$
 since $\expec\big[g(\mu_n)\big]$ converges.
Moreover, since $\mathcal O$ is dense in $C(\mathcal P)$,
 the first term can be made arbitrarily small, implying 
 \[\limsup_{n,m \to +\infty} \big| \expec\big[f(\mu_n)\big] - \expec\big[f(\mu_m)\big]\big| =0.\]
 In other words, $ \expec\big[f(\mu_n)\big]$ is a Cauchy sequence and hence, converges.
 
 Now, by compactness of $\mathcal P$, and hence of the set of probability measures on $\mathcal P$ thanks to Prokhorov's theorem, 
 the sequence $\mu_n$ admits a convergent (in distribution) subsequence $\mu_{n_k}$.
 Let us denote $\mu$ its limit.
 Using the same argument as above for $d_\tau$, we have that for any function $f$ in $C(\mathcal P)$,
\[ \lim_{k \to +\infty} \expec\big[f(\mu_{n_k})\big] = \expec\big[f(\mu)\big]. \]
On the other hand, since $\expec\big[f(\mu_n)\big]$ converges, its subsequence $\expec\big[f(\mu_{n_k})\big]$
converges to the same limit, i.e.~we have
\[ \lim_{k \to +\infty} \expec\big[f(\mu_{n_k})\big] = \lim_{n \to +\infty} \expec\big[f(\mu_{n})\big].\]
Summarizing, we have that $ \expec\big[f(\mu_n)\big]$ converges to $\expec\big[f(\mu)\big]$
for any continuous function $f$, i.e.~$\mu_n$ converges to $\mu$ in distribution.
\end{proof}

\begin{remark}
The equivalence $(i) \Leftrightarrow (iii)$ can also be seen as a consequence of
 the following general result of Le Cam (see \cite{lecam1957convergence_determining} or \cite[Theorem 2.7]{loehr2013GromovBox}): 
Let $X$ be a completely regular Hausdorff space (in particular, a metric space), 
and let $\mathcal F \subseteq C_b(X)$ be a family of bounded continuous functions 
that is closed under multiplication and generates the topology of $X$, meaning that the topology of $X$ is the coarsest topology making every $f\in\mathcal F$ continuous.
Then $\mathcal F$ is convergence-determining for Radon probability measures on $X$: 
\begin{equation}\label{eq:convergence-determining}
 \nu_n \stackrel{d}{\to} \nu \ \text{ if and only if } \ 
\lim_{n \to +\infty}  \int f(x) \nu_n(dx) =  \int f(x) \nu(dx) \quad\text{for every  $f\in\mathcal F$}. 
\end{equation}
In our setting, we take $X=\mathcal P$, $\mathcal F=\mathcal O$, and $\nu_n$ and $\nu$ to be the distributions of $\mu_n$ and $\mu$. Then the above conditions are fulfilled, and the equivalence \eqref{eq:convergence-determining} corresponds to $(i) \Leftrightarrow (iii)$ with the following caveat. To prove $(iii) \Rightarrow (i)$, we need first to identify $c_\tau$
with $ \int d_\tau(x) \nu(dx)$ for some measure $\nu$ on $\mathcal P$---this can be done using a compactness argument as in the above proof.
\end{remark}

\subsection{Intensity measures} \label{subsec:intensity-measures} 
Given a random permuton $\mu$, its \textbf{intensity measure} is the deterministic measure $\mathbb{E}[\mu]$ defined by 
\[ 
	\mathbb{E}[\mu](A)=\mathbb{E}[\mu(A)], \qquad A\subseteq[0,1]^2\text{ Borel}. 
\] 
It is itself a (deterministic) permuton. For any bounded continuous function $f$ on $[0,1]^2$, we have 
\[ 
	\int f\,d\mathbb{E}[\mu] =\mathbb{E}\!\left[\int f\,d\mu\right]. 
\] 
Consequently, if $\mu_n$ converges in distribution to $\mu$, then $\mathbb{E}[\mu_n]$ converges weakly to $\mathbb{E}[\mu]$. Indeed, the map $\nu\mapsto\int f\,d\nu$ is bounded and continuous on $\mathcal{P}$, so convergence in distribution implies convergence of its expectation. See also \cite[Section~1.3.3]{feray2023recursive_separable} for this interpretation of intensity measures as limits of averaged permutation diagrams. 

For a random permutation $\pi$ of size $n$, the intensity measure $\mathbb E[\mu_\pi]$
 is determined by the one-point probabilities $\mathbb{P}(\pi(i)=j)$.
 More precisely, using in place of $\mu_\pi$ the empirical measure
 $\nu_\pi=\frac1n\sum_{i=1}^n\delta_{(i/n,\pi(i)/n)}$ for simplicity, we have
\begin{equation} \label{eq:intensity-one-point-probabilities}
	 \mathbb{E}\!\left[\int f\,d\nu_\pi\right] =\frac1n\sum_{i,j=1}^n f\!\left(\frac{i}{n},\frac{j}{n}\right) \mathbb{P}\bigl(\pi(i)=j\bigr). 
\end{equation} 
Recall that the permuton $\mu_\pi$ defined in~\eqref{eq:defn-perm-ass} replaces each atom of $\nu_\pi$ with Lebesgue measure over the corresponding square. Thus, replacing $\nu_\pi$ by $\mu_\pi$ in the left-hand side of~\eqref{eq:intensity-one-point-probabilities} changes it by at most $\omega_f(\sqrt{2}/n)$, where $\omega_f(\delta)=\sup_{z,z'\in[0,1]^2:\,\|z-z'\|\leq\delta} |f(z)-f(z')|$ is the modulus of continuity of $f$. In particular, the same formula~\eqref{eq:intensity-one-point-probabilities} holds for $\mu_\pi$ with an error $o(1)$ as $n\to\infty$. 

Estimation of the values of $\mathbb{P}\bigl(\pi(i)=j\bigr)$ by sampling $\pi$ is sometimes referred to as a \textbf{heatmap} of the random permutations $\pi$ in the permutation pattern literature; see, e.g.,~\cite[Figure 9]{combinatorial-exploration} for all heatmaps of uniform random permutations avoiding a pair of patterns of size 4. These heatmaps are approximations of the intensity measures of the potential limiting random permuton of $\pi$.

We warn the reader that the pattern densities of the intensity measure need not equal the expected pattern densities of the random permuton. If $\pi_n$ is a random permutation of size $n$ and $\mu_{\pi_n}$ converges in distribution to $\mu$, then Theorems~\ref{thm:cv_permutons_densities} and~\ref{thm:cv_random_permutons} (together with Lemma~\ref{lem:densityofassociatedpermuton}) give 
\[ 
	\begin{aligned} \mathbb{E}[d_\tau(\pi_n)] &\longrightarrow \mathbb{E}[d_\tau(\mu)],\\ d_\tau\bigl(\mathbb{E}[\mu_{\pi_n}]\bigr) &\longrightarrow d_\tau\bigl(\mathbb{E}[\mu]\bigr). \end{aligned} 
\] 
The two limits need not agree. In the first case, all sampled points are chosen conditionally independently from \emph{the same} realization of $\mu$; in the second, they are independent samples from $\mathbb{E}[\mu]$, which amounts to choosing an independent realization of $\mu$ for each point. Averaging the measure therefore loses the dependence between points sampled from a common realization. 

For example, let $A,U$ be independent uniform random variables on $[0,1]$, and let $\mu_A$ be the conditional distribution of $(U,(U+A)\bmod 1)$ given $A$. Then $\mathbb{E}[\mu_A]=\Leb_{[0,1]^2}$, whereas \[ \mathbb{E}[d_{21}(\mu_A)] =\int_0^1 2a(1-a)\,da =\frac13 \neq\frac12 =d_{21}(\Leb_{[0,1]^2}). \] Indeed, for a fixed $a$, an inversion occurs exactly when the two sampled $x$-coordinates lie on opposite sides of $1-a$. Alternating between copies of $\mu_A$ and the deterministic uniform permuton gives a sequence with constant intensity measures which does not converge in distribution. Thus convergence of intensity measures alone does not imply permuton convergence.

\section{Feasible region, quasi-randomness and admissible supports}\label{sect:feas-reg}

Given a permuton $\mu$, we defined in the previous section its pattern densities $d_\tau(\mu)$, indexed by permutations $\tau$ of all sizes. Since $\mu$ is a measure, we can also consider its support, which is a closed subset of $[0,1]^2$. In this section, we discuss some inverse problems: does there exist a permuton with some given pattern densities or some given support? In Section~\ref{sect:feas-reg-1}, we are interested in fixing a finite number of pattern densities, while Section~\ref{sect:packing-problem} considers the possible values of a single pattern density. Section~\ref{sect:quasirandomness} looks at uniqueness problems in this setting and searches for pattern density constraints for which the uniform permuton is the only solution. Finally, Section~\ref{sect:perm-support} addresses the existence or non-existence of permutons with given supports.

\subsection{The feasible region and its dimension}\label{sect:feas-reg-1}
An important question raised, e.g., by \Cref{thm:cv_permutons_densities}, is understanding the set 
\begin{equation}\label{eq:defn-feas-reg}
    \mathcal F=\big\{(d_\tau(\mu))_{\tau \in \Sym}: \mu \text{ is a permuton} \big\} \subseteq (\mathbb R_+)^{\Sym},
\end{equation}
{called the \textbf{feasible region} of pattern densities.}
Since this lives in an infinite dimensional space, we will rather ask questions about finite dimensional projections.
In terms of large permutations, this can be equivalently stated as follows 
(thanks to \Cref{thm:cv_permutons_densities} and \Cref{cor:permutations_dense_in_P}):
\[
\text{What are the possible limiting values of $d_\tau(\sigma_n)$ for a given pattern $\tau$, or a family of patterns $\tau$?}
\]
Classical problems of extremal combinatorics, such as the packing problem (see Section~\ref{sect:packing-problem}) or quasirandomness criteria (see Section~\ref{sect:quasirandomness}), fall into this framework.

\smallskip

There is another way to look at the feasible region $\mathcal F$.
Recall from \Cref{prop:densities_span_algebra} that the functions $d_\tau$ span an algebra linearly.
{In fact, looking at the proof of \Cref{prop:densities_span_algebra}, we can do a bit more: we can construct an associative
operation\footnote{The coefficients
$c_{\tau_1,\tau_2}^\tau$ have in fact some probabilistic interpretation as decomposition of $\tau$ as a disjoint union
of an occurrence of $\tau_1$ and an occurrence of $\tau_2$, so that $(\mathbb R\Sym, \ast)$ is a kind of shuffle algebra,
see, e.g., \cite[Proposition 17]{feray2020moudli_spaces}.} $\ast$ on the real vector space with basis $\Sym$ such that, if
\[
\tau_1 \ast \tau_2:=\sum_{\tau \in \Sym} c_{\tau_1,\tau_2}^\tau \, \tau,
\]
then we have
\[ d_{\tau_1} d_{\tau_2} = \sum_{\tau \in \Sym} c_{\tau_1,\tau_2}^\tau d_\tau.\]}
Then, every permuton $\mu$ defines a nonnegative homomorphism $\phi_\mu:\tau \mapsto d_\tau(\mu)$
from $(\mathbb R\Sym, \ast)$ to $\mathbb R$, where nonnegative means that $\phi_\mu(\tau) \ge 0$ for all $\tau$ in $\Sym$. We call such a homomorphism normalized if $\phi(1)=1$, where $1$ denotes the one-point permutation.
Conversely, Razborov \cite[Theorem 3.3]{razborov2007flag} showed in a rather general setting that any normalized nonnegative homomorphism from $\mathbb R \Sym$ to $\mathbb R$ corresponds to some $\mu$ in $\mathcal P$, i.e.~
\begin{equation}\label{eq:feasible_homomorphism}
\mathcal F=\big\{(\phi(\tau))_{\tau \in \Sym}, \phi \text{ normalized nonnegative homomorphism }(\mathbb R \Sym,\ast) \to \mathbb R \} \subseteq (\mathbb R_+)^{\Sym}. 
\end{equation}

\begin{figure}[t]
\[\includegraphics[width=0.48\textwidth]{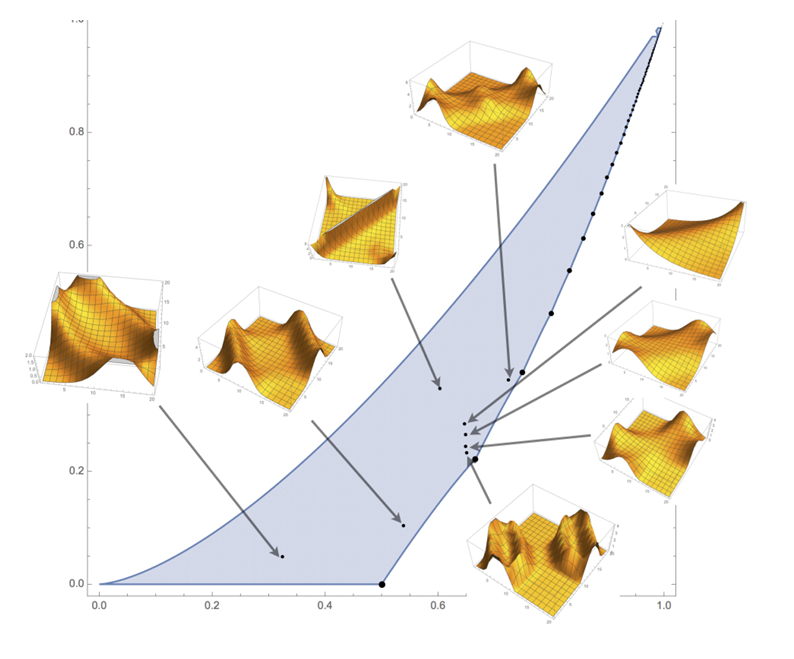}
\includegraphics[width=0.48\textwidth]{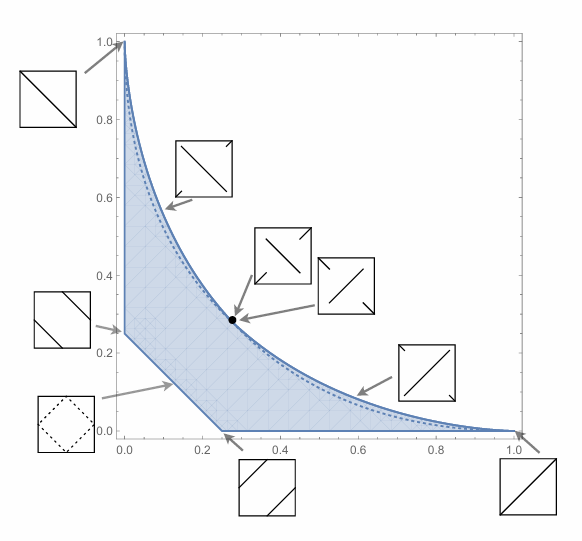}
\]
\caption{{These figures are reproduced from~\cite{ Kenyon-Permutations-densities}. \textbf{Left:} The feasible region for $(d_{12}(\mu),d_{123}(\mu))$, with corresponding permutons (computed numerically) at selected points.
\textbf{Right:} The feasible region for $(d_{123}(\mu),d_{321}(\mu))$. It is bounded above by the
parameterized curves $(1-3t^2+2t^3,t^3)$ and $(t^3,1-3t^2+2t^3)$ which intersect
at $(x,y)=(0.278\ldots,0.278\ldots)$. The lower boundaries consist of the axes and the
line $x+y=1/4$.}
\label{fig:feas-ref}}
\end{figure}

Despite this characterization, the feasible region $\mathcal F$ is a complicated object 
and a concrete description seems out of reach.
Looking at the one-dimensional projection is equivalent to the so-called pattern packing problem,
which we shall discuss in detail in the next section.
In \cite{ Kenyon-Permutations-densities}, building on results concerning feasible regions for subgraph densities,
Kenyon, Král', Radin, and Winkler consider the projection
of $\mathcal F$ on $\mathbb R^{\{12,123\}}$ and $\mathbb R^{\{123,321\}}$, 
i.e., the possible values of the pair $(d_{12}(\mu),d_{123}(\mu))$ or $(d_{123}(\mu),d_{321}(\mu))$.
The feasible regions for these patterns already have complicated structures, see Figure~\ref{fig:feas-ref}.

\smallskip

In another direction, for $n \ge 1$, let us denote $\Sym_{\le n}$ the set of patterns (or permutations) of size at most $n$,
and $\mathcal F_{\le n}$ for the projection of $\mathcal F$ to $(\mathbb R_+)^{\Sym_{\le n}}$.
Since one cannot hope for a complete description of $\mathcal F_{\le n}$, we may ask at least what its dimension is.
In \cite{vargas2014hopf_permutations2}, it is shown that the algebra $(\mathbb R\Sym, \ast)$ is freely generated
by a certain family of permutations called {\em Lyndon permutations}
(by analogy with Lyndon words, which generate freely the shuffle algebra \cite{radford1979basis_shuffle}).
Hence $\mathcal F_{\le n}$ is contained in an algebraic variety of dimension $\ell_{\le n}$, 
the number of nontrivial Lyndon permutations of size at most $n$.
On the other hand, if $\mathcal L_{\le n}$ is the set of nontrivial Lyndon permutations of size at most $n$,
it was proven by Garbe, Král', Malekshahian and Penaguiao \cite{garbe2023dimensionfeasibleregionpattern}
that the projection of $\mathcal F$ on $\mathbb R^{\mathcal L_{\le n}}$ has a nonempty interior.
Hence, we have the following nice result. 
{\begin{thm}[{\cite{garbe2023dimensionfeasibleregionpattern}}]\label{thm:fes-dim-lynd}
    The dimension of $\mathcal F_{\le n}$ is exactly the number  $\ell_{\le n}$ of nontrivial Lyndon permutations of size at most $n$.
\end{thm}}
{The number $\ell_n$ of Lyndon permutations of size $n$ can be computed recursively.
Indeed, by a result of Vargas~\cite{vargas2014hopf_permutations2}, Lyndon permutations are free generators of a commutative algebra over permutations (using a {\em bi-shuffle} product). This implies via standard Hilbert series computation that
\[\sum n! t^n = \prod_{n \ge 1} \frac{1}{(1-t^n)^{\ell_n}}.\]
It follows that $\ell_n$ corresponds to the sequence A112354 in OEIS~\cite{oeis}. In particular, asymptotically $\ell_n \sim n!$, meaning that $\mathcal F_{\le n}$ has essentially full dimension.}

{\begin{remark}
	A related feasible-region problem concerns probability distributions on $\Sym_n$ induced by ranking independent, not necessarily identically distributed, real-valued random variables. More precisely, let $X_1,\ldots,X_n$ be independent random variables with almost surely no ties, and let $\sigma\in\Sym_n$ be the random permutation determined by $X_{\sigma(1)}<\cdots<X_{\sigma(n)}$.
	Writing $p_\tau=\mathbb P(\sigma=\tau)$, consider the set of achievable ranking laws
	\[
	\mathcal A_n=\left\{(p_\tau)_{\tau\in\Sym_n}: X_1,\ldots,X_n \text{ are independent}\right\}.
	\]
	Turvey~\cite{Tur26} recently proved that the dimension of $\mathcal A_n$ is $\sum_{k=2}^n\binom nk(k{-}1)!$, resolving a conjecture of Babson et al.~\cite{BDIPTT26}. As in the dimension problem for the feasible region of permuton pattern densities (Theorem~\ref{thm:fes-dim-lynd}), the proof exploits shuffle-algebra identities and Lyndon decompositions.
\end{remark}}

\subsection{The packing problem}\label{sect:packing-problem}
In this section, we consider one-dimensional projections of the feasible region,
i.e.~we fix a pattern $\tau$ of size $k\ge 2$, and want to find the set
\[\mathcal F_\tau=\{d_{\tau}(\mu), \mu \in \mathcal P\}.\]
It is easy to see that $0$ is always in $\mathcal F_\tau$: if $\tau$ is different from the identity,
take $\mu$ to be the Lebesgue measure on the main diagonal $\{y=x\}$ of the square,
and if $\tau$ is the identity, we take $\mu$ to be the Lebesgue measure on the other diagonal $\{y=1-x\}$.
One can also easily see that $\mathcal F_\tau$ is a closed interval, and thus $\mathcal F_\tau=[0,p(\tau)]$,
for some number $p(\tau)$ in $[0,1]$.
By definition, we have
\begin{equation}\label{eq:def_packing}
p(\tau)= \max \{d_{\tau}(\mu), \mu \in \mathcal P\} = {\sup}_{\sigma_n:|\sigma_n| \to +\infty} \big\{{\limsup}_{n\to \infty} d_{\tau}(\sigma_n)\big\},
\end{equation}
where the {supremum} is taken over sequences of permutations whose sizes tend to $+\infty$.
In words, we try to pack as many patterns $\tau$ as possible into larger and larger permutations,
whence the name of \enquote{packing problem.}
Under this form, this problem has been around in the combinatorics literature for some time,
see e.g.~\cite{price1997packing,albert2002packing}.

Earlier results on packing problems considered small patterns and used ad-hoc arguments.
In particular, when $\tau$ is layered, meaning that it can be written as a direct sum of decreasing permutations, the packing density $p(\tau)$ can be achieved by a sequence of layered permutations; see, e.g., \cite[Theorem 5.36]{bona_book_permutations}---and optimizing over the set of layered permutations is sometimes more tractable than over all permutations.
{The structure of such optimizing layered permutations can nevertheless be rather subtle. Kabela, Král', Noel, and Pierron \cite{Kabela-22-layred} showed that, for some layered patterns, every asymptotically maximizing sequence of layered permutations must have a number of layers tending to infinity; they also gave sufficient conditions ensuring the existence of asymptotic or exact maximizers with a bounded number of layers.}
The reformulation as an optimization problem on measures on the unit square was first found by Presutti and Stromquist \cite{presutti2010packing} (before permutons were introduced by Hoppen and co-authors). They utilized it to find a new lower bound for the packing density of $2413$ by constructing a measure with a high density of occurrences of $2413$.

However, the strongest improvements in the bounds in this area were achieved through the use of flag algebras. Flag algebras, originating in graph theory,
were introduced by Razborov in the seminal paper~\cite{razborov2007flag}.
He uses the description in \eqref{eq:feasible_homomorphism} of $\mathcal F$ as the image of normalized nonnegative homomorphisms from 
$(\mathbb R \Sym,\ast)$ to $\mathbb R$ to obtain some inequalities fulfilled by the pattern densities in $\mathcal F$.
The power of the method is that it recasts the problem as a semidefinite programming problem for which
there exist some efficiently implemented methods.
We will not explain flag algebras in detail here, since this is not the topic of this survey,
and since there are already excellent references: {see~\cite{silva2016survey_flags}
for a survey of applications} in the context of graph theory,
and \cite{balogh2015monotone_subsequences,sliacan2017packing_flag,crudele2023six_patterns_quasirandomness} 
for references in the framework of permutations.

Returning to the packing problem, flag algebras have yielded improved upper bounds for patterns of length 4, as well as bounds for longer patterns.
For all patterns of length $4$, and for some of length 5 and 6,
one can construct a measure $\mu$ so that $d_\tau(\mu)$ is close to the proved
upper bound for $p(\tau)$, so that the value of $p(\tau)$ is known up to a small error.
All these results are presented in the article by Slia{\v{c}}an and Stromquist \cite{sliacan2017packing_flag}.

\subsection{Quasi-randomness}\label{sect:quasirandomness}
Another classical question related to the feasible region is that of quasirandomness.
A sequence of permutations $\sigma_n$, {with $|\sigma_n|\to\infty$,} is said to be \textbf{quasirandom} if, for any pattern $\tau$,
the density $d_\tau(\sigma_n)$ converges to $1/|\tau|!$---which is the limiting value for a sequence
of uniform random permutations $\sigma_n$.
In terms of permutons, this means that $\sigma_n$ converges to the uniform permuton $\Leb([0,1]^2)$, where $\Leb$ denotes the Lebesgue measure.
Such a quasirandomness notion can be defined for various combinatorial objects,
and it turns out that, in many cases, only the convergence of a {\em finite number} of substructure densities is enough to ensure quasirandomness,
see~\cite{chung1989quasi-random,chung1990QRhypergraphs,chung1991QRtournaments,chung1991QRSetSystems,cooper2022quasirandomLatin}.
For permutations, one can prove that $\sigma_n$ is quasirandom if and only if $d_\tau(\sigma_n)$ converges to $1/24$ for all patterns of size 4.
This result was first announced\footnote{There seem, however, to be some missing arguments in the proof.} by Yanagimoto \cite{yanagimoto1970quasi_randomness} 
and then later rediscovered by Král' and Pikhurko \cite{kral2013quasirandom}
{(in fact, the proof of \cite{kral2013quasirandom}
shows that specifying the values of 3 linear combinations of the 24 pattern densities
of size 4 is enough to characterize quasirandomness)}.
In terms of feasible region, this means that the full infinite-dimensional feasible region $\mathcal F$
only contains a single point $(\alpha_\tau)_{\tau \in \Sym}$ such that $\alpha_{\tau}=1/24$ for all patterns $\tau$ of size $4$.

This is of importance in statistics, for non-parametric independence tests.\footnote{The above-cited article by Yanagimoto \cite{yanagimoto1970quasi_randomness}
and the article by Bergsma and Dassios \cite{bergsma2014test_patterns}, cited below, are written in the language of statistical tests.}
Indeed, let us consider an i.i.d.~sample $(X_i,Y_i)_{i \le n}$ of an unknown random variable $(X,Y)$ in $\mathbb R^2$.
We want to test the independence of $X$ and $Y$. 
We have no a priori knowledge of the law of $X$ and $Y$, except that they are continuous random variables, in
the sense that their cumulative distribution functions $F_X$ and $F_Y$ are continuous
(we do not want to assume that they belong to some specific family; this is a non-parametric test).
The idea, suggested by Hoeffding \cite{hoeffding1948test}, is to consider the permutation $\sigma_n$ associated with the points
$(X_i,Y_i)_{i \le n}$ in the sense of \Cref{ssec:samples_moments_convergence}.
If $X$ and $Y$ have uniform marginal distributions, then the distribution of $(X,Y)$ is a permuton $\mu$,
and $\sigma_n$ follows the law of the permuton sample $\sigma_n^\mu$.
In the general case, as explained in Remark~\ref{rmk:trans-perm}, $\sigma_n$ still has the law of some permuton sample $\sigma_n^\mu$,
but $\mu$ is the law of
\begin{equation}\label{eq:meas-to-perm}
    (F_X(X),F_Y(Y)).
\end{equation}
In particular, $X$ and $Y$ are independent if and only if $\mu$ is the uniform permuton,
which occurs if and only if $\sigma_n$ is quasirandom.
Thanks to the above result, quasirandomness can be tested by only looking at a finite number of pattern densities
(in fact, a single linear combination suffices).
We refer to \cite{evenZohar2021counting_patterns} for details, numerous references,
and for algorithmic results on the complexity of computing these pattern densities.

In view of this application, it is natural to try to find smaller sets of patterns 
or simpler linear combinations that characterize quasirandomness---in which case they are called {\bf quasirandom-forcing}. The first result of this kind is due to Bergsma and Dassios \cite{bergsma2014test_patterns}, who found a quasirandom-forcing linear combination of 8 patterns of size 4, but only for a restricted class of permutons (mixtures of densities and atoms, but no other types of singular parts). Their result was later generalized to all permutons~\cite{drton2020quasirandom_tau}, and very recently even to measures whose marginals may have atoms~\cite{bergsma2026quasirandomness_with_ties}.

In a later paper, Chan, Král', Noel, Pehova, Sharifzadeh, and Volec \cite{kral2020quasirandom_pattern_sum} completely characterize the sets of size-4 patterns whose sum is quasirandom-forcing.  More recently, Crudele, Dukes, and Noel \cite{crudele2023six_patterns_quasirandomness} found a nonnegative linear combination
of only six permutations (of sizes 3 and 4) which characterizes quasirandomness.
These results go in pairs with necessary condition results:
any quasirandom-forcing set of patterns of any size must have size at least 4~\cite{kurecka2022lower_bound_quasirandom};  any quasirandom-forcing set of patterns of size four must have size at least 5~\cite{kral2024forcing};
and a nonnegative quasirandom-forcing linear combination must involve at least 6 patterns~\cite{crudele2023six_patterns_quasirandomness}
(so that the result of Crudele, Dukes, and Noel is optimal when we restrict to nonnegative coefficients).
As for packing density results, the framework of flag algebras 
has been an important tool in this line of research.\medskip

\subsubsection*{Finite forcibility}
We conclude this section by mentioning briefly a concept related to quasi-randomness.
A permuton $\mu$ is called {\em finitely forcible} if there exists a finite set of patterns $S$ such that, for any permuton $\nu$, we have
\begin{equation}\label{eq:def-finitely-forcible} 
d_{\tau}(\mu) = d_{\tau}(\nu) \text{ for all $\tau$ in $S$} \ \Longrightarrow \ \nu=\mu.
\end{equation}
Recall that, by Proposition~\ref{prop:densities_separate}, it is always the case that $\mu=\nu$ whenever $d_{\tau}(\mu) = d_{\tau}(\nu)$ for {\em all} patterns $\tau$. The challenge in the above definition is therefore to
find a {\em finite} subset of patterns 
for which the implication \eqref{eq:def-finitely-forcible} holds.
Král'--Pikhurko's result above implies that the uniform permuton is finitely forcible
(taking $S$ as the set of patterns of size $4$, and further results found smaller and smaller sets $S$ such that \eqref{eq:def-finitely-forcible} holds).
Further examples of finitely forcible permutons are given in~\cite{glebov15-forcible}.
On the other hand, an example of a permuton that is not finitely forcible can be found
in~\cite{glebov17-permuton-testing}. 
Completely characterizing finitely forcible permutons remains an open problem.

\subsection{Sets that can support a permuton}\label{sect:perm-support}
\begin{figure}[t]
\includegraphics[width=0.7\textwidth]{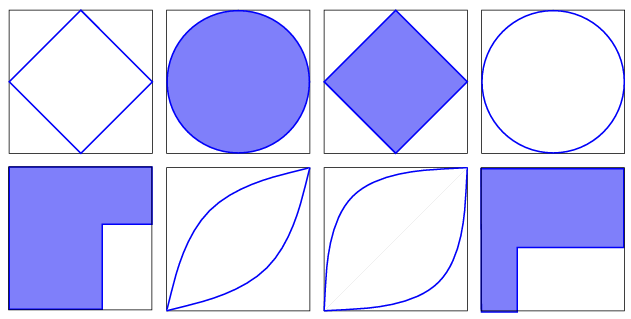}
\caption{Which of these subsets of $[0,1]^2$ are the support of a permuton? 
\label{fig:Support-permutons}}
\end{figure}

Recall that the support of a permuton $\mu$ is the smallest closed subset $S$ of $[0,1]^2$ such that $\mu(S)=1$.
This leads to the following question, illustrated by Figure~\ref{fig:Support-permutons}: which closed sets $S \subseteq [0,1]^2$ can be the support of a permuton?
This question was addressed in the recent preprint by Coscia, Tassy, and Winkler ~\cite{CTW26}, whose results are briefly summarized below.

{Let $X=Y=[0,1]$ and $S\subseteq X \times Y = [0,1]^2$ be closed.
For a Borel set $U\subseteq X$, write
\[
N(U):=\bigl\{y\in[0,1]:(x,y)\in S\text{ for some }x\in U\bigr\}.
\]
We say that $U$ is \textbf{compliant} if
\[
\lambda\bigl(N(U)\bigr)\geq\lambda(U),
\]
where, as earlier, $\lambda$ denotes Lebesgue measure on $[0,1]$.
The set $S$ itself is said to be compliant if every Borel set $U \subset X$ is compliant.
It is not hard to see that $S$ is compliant if and only if there exists a permuton $\mu$ such that
\[
\operatorname{supp}(\mu)\subseteq S;
\]
see~\cite[Corollary~3]{CTW26}. Figuring out under what conditions the support can be {\em equal} to $S$ is more delicate.
A Borel set $U \subset X$ is said to be \textbf{expansive} if, in the event that $\lambda\bigl(N(U)\bigr) = \lambda(U)$,
there are no open intervals $A \subset X \setminus U$ and $B \subset Y$ such that $\lambda(B\setminus N(U))=0$ and $B\cap N(A)\neq \emptyset$.
The set $S$ itself is termed expansive if every Borel $U \subset X$ is expansive.

\begin{thm}[{\cite[Theorem~11]{CTW26}}]
	If $S$ is regular (i.e., $S$ is the closure of its interior), then $S$ is the support of a permuton if and only if
	it is compliant and expansive.
\end{thm}

With help from this result, one can see that the four sets on the left of Figure~\ref{fig:Support-permutons}
are the support of a permuton, while the four sets on the right are not. See \cite[Section~4]{CTW26} for details.
For arbitrary closed sets, compliance and expansiveness are not sufficient in general; an exact characterization
(of sorts) is given in~\cite[Theorem~10]{CTW26} in terms of a certain local \emph{shavability} condition.}

\section{Permutons and large deviation theory}
\label{sec:LDP}

\subsection{Large deviation principle for uniform random permutations}\label{sect:LDP-sect}
Consider a uniform random permutation $\sigma_n$ in $\Sym_n$.
It is a simple exercise to show that $\sigma_n$ converges almost surely to
the uniform permuton $\mu=\Leb([0,1]^2)$. 
This can also be seen as a special case of \Cref{prop:convergence_samples},
since for $\mu=\Leb([0,1]^2)$, the permutation $\sigma_n^\mu$ is uniformly distributed on $\Sym_n$.

Hence, if $A$ is a subset of the space $\mathcal P$ of permutons whose closure $\overline A$ does not contain the 
permuton $\Leb([0,1]^2)$, then the quantity $\mathbb P(\mu_{\sigma_n} \in A)$ tends to $0$.
Large deviation theory is interested in the following question: 
\[
\text{How rare is the event $\{\mu_{\sigma_n} \in A\}$?}
\]
In many settings, such as for i.i.d.\ random variables, the analogous events have exponentially small probability, and we will see that the same phenomenon occurs here for permutations.
The question is to find the exponential rate of decay of $\mathbb P(\mu_{\sigma_n} \in A)$,
i.e.~to find the limit (if it exists) of 
\[\tfrac1n \log\big(  \mathbb P(\mu_{\sigma_n} \in A) \big).\]

To answer this question, it is useful to define the entropy of a permuton $\mu$ as
\[H(\mu) = \begin{cases} - \int_{[0,1]^2} g(x,y) \log(g(x,y)) dx dy & \text{if $\mu$ has a density $g(x,y)$ with respect to $\Leb([0,1]^2)$,}\\
-\infty& \text{otherwise.}
\end{cases}\]
(By convention $0 \log(0)$ is interpreted as $0$.)
Using Jensen's inequality, we see that $H(\mu) \le 0$ for any permuton $\mu$
with equality if and only if $\mu=\Leb([0,1]^2)$.
Moreover, it is known that $H$ is a strictly concave 
and upper-semicontinuous function on $\mathcal P$; see, e.g.~\cite[Section III, Theorem 1]{entropy1975}.
The following result is due to Trashorras\footnote{In fact, Trashorras associates with a permutation $\sigma$ the measure $\nu_\sigma=\frac1n \sum_{i=1}^n \delta_{(i/n,\sigma(i)/n)}$ and not the permuton $\mu_\sigma$ considered here, but since these measures are asymptotically close to each other, {large deviation principles for one and for the other are equivalent.}}~\cite{trashorras2008large_deviations}; but see also~\cite[Theorem 4.1 in the arXiv version]{Mukherjee-estimation} and~\cite[Theorem 1]{ Kenyon-Permutations-densities}. 

\begin{thm}\label{thm:large_deviation}
Let $A$ be a subset of the space $\mathcal P$ of permutons and denote by $\Int(A)$ the interior
of the set $A$. Then
\[\sup_{\mu \in \Int(A)} H(\mu) \le \liminf_{n \to +\infty} \tfrac1n \log\big( \mathbb P(\mu_{\sigma_n} \in A) \big)
\le  \limsup_{n \to +\infty} \tfrac1n \log\big( \mathbb P(\mu_{\sigma_n} \in A) \big) \le \sup_{\mu \in \overline A} H(\mu). \]
\end{thm}
Such a result is called a {\bf large deviation principle} (LDP) for $\sigma_n$, with
$n$ being the speed of the LDP and {$-H$} the rate function.
Large deviation principles are ubiquitous in probability theory:
a classical reference is
 the book by Dembo and Zeitouni~\cite{dembo-ld}.
As we shall see, this LDP is useful for studying large random permutations constrained to
satisfy some rare events (having a small number of inversions, for example),
or whose distribution is biased with respect to some statistics.

\begin{remark}\label{rmk:further_LDP}
Since many models of non-uniform random permutations are discussed in this survey,
it is natural to ask about large deviation principles for some of these models.
Considering the exponentially biased random permutation (Section~\ref{sect:exponentially-biased-permutations}), as in the Mallows model,
will lead to a simple adjustment of the rate function.
Another family of non-uniform random permutation models is the family of permuton
samples, discussed in \Cref{sec:theory}.
A large deviation principle for such models, generalizing \Cref{thm:large_deviation},
 is given by Borga, Das, Mukherjee and Winkler in \cite{borga2022large}.
 
  In another direction, large deviation principles for 231 and 321-avoiding permutations
 (and thus, by symmetry, for permutations avoiding any pattern of size $3$)
 have been given by Budzinski, Dubach, F\'eray, Kammoun, and Ma\"ida in the recent preprint \cite{budzinski2026constrained_mallows}.
 It would be interesting to find large deviation principles for some further models of random
 pattern-avoiding permutations, such as those mentioned in \Cref{sect:pattern-av-perm}.
 Some general questions can also be considered: e.g., should we expect to have an LDP at speed $n$
 for all exponentially growing permutation classes, or can other speeds appear?
\end{remark}

\subsection{Exponentially biased permutations}\label{sect:exponentially-biased-permutations}
This section largely follows \cite{Mukherjee-estimation}.
 We consider non-uniform random permutation models,
defined through an exponential change of probability measures
with respect to the uniform distribution.
Namely, if $f:\Sym \to \mathbb R$ is a permutation statistic,
and $q=q_n>0$ is a positive real parameter (possibly depending on $n$),
then we let $\sigma_{n,q}^f$ be the random permutation taking values in $\Sym_n$ 
with the following distribution:
\begin{equation}\label{eq:def_exp_biased}
\prob \big(\sigma_{n,q}^f=\tau \big) = \frac{q^{f(\tau)}}{\sum_{\pi \in \Sym_n} q^{f(\pi)}}, 
\text{ for $\tau$ in $\Sym_n$.}
\end{equation}
Given a function $f:\Sym \to \mathbb R$, we call $f^o$
a \textbf{continuous version of $f$} on $\mathcal P$ if $f^o$ is continuous and
\begin{equation}\label{eq:def_continuous_version}
\lim_{n \to +\infty} \sup_{\sigma \in \Sym_n} | f(\sigma) -f^o(\mu_\sigma)| =0.
\end{equation}
For example, the function $d_\tau$ defined by \eqref{eq:def_dtau_permuton} on $\mathcal P$ is a continuous version of the function
$d_\tau$ defined by \eqref{eq:def_dtau_permutation} on $\Sym$. In that case, we made a slight abuse of notation and used the same notation for the function on $\Sym$ and its continuous version on $\mathcal P$.
We are interested in finding a limiting permuton for the sequence of permutations
$\sigma_{n,q}^f$ when $f$, or one of its normalizations, 
has a continuous version on $\mathcal P$.
Then the large deviation principle turns such limit shape problems
into analytic optimization problems, as shown by the following result.
\begin{prop}\label{prop:limit_sigma_qf}
Let $\sigma_{n,q_n}^f$ be the exponentially biased permutation defined in \eqref{eq:def_exp_biased}. Assume that there exists a sequence {$a_n\in\mathbb{R}_{>0}$}
 such that:
\begin{enumerate}[(i)]
\item the function $\tilde f:\Sym \to \mathbb R$ defined by $\tilde f(\sigma)=f(\sigma)/a_{|\sigma|}$ has a continuous version $\tilde f^o$ on $\mathcal P$;
\item the sequence $a_n \log(q_n)/n$ has a limit $\beta \in (-\infty,+\infty)$.
\end{enumerate}
Then the normalization constant in \eqref{eq:def_exp_biased} satisfies
\begin{equation}\label{eq:estimate_logZ}
 \lim_{n \to +\infty} \frac1n \log\left(\frac{1}{n!} \sum_{\pi \in \Sym_n} q_n^{f(\pi)}\right)= \max_{\mu \in \mathcal P} \ H(\mu) + \beta \tilde f^o(\mu),
 \end{equation}
and the support of any distributional limit of a subsequence of $\sigma_{n,q_n}^f$
is included in the set of permutons $\mu$ achieving the maximum above. 
\end{prop}
Naturally, this proposition is particularly meaningful when the maximum is achieved by a single element $\mu_0$ in $\mathcal P$, in which case it implies that $\sigma_{n,q}^f$
converges weakly to the deterministic permuton $\mu_0$.

Results such as \Cref{prop:limit_sigma_qf} are standard in large deviation theory.
As far as we are aware, 
\Cref{prop:limit_sigma_qf} does not appear as is in the literature,
but it has been stated for some special cases where the maximizer is unique
(see examples below).

\begin{proof}
The proof can be done as in \cite[Theorem 1.4, proofs of items (a) and (b)]{Mukherjee-estimation}, using general results from large deviation theory
(in particular \cite[Theorem 4.3.1 and Lemma 4.3.6]{dembo-ld}).
Nevertheless, to give the reader some intuition for the connection between large deviation theory and such limit statements for biased random permutations, we provide a proof based only on \Cref{thm:large_deviation}. Thanks to the compactness of $\mathcal P$, the argument is slightly simpler than those in \cite[Theorem 4.3.1 and Lemma 4.3.6]{dembo-ld}, although the underlying idea is the same.

\medskip

Let us set $Z_n(q_n,f)=\sum_{\pi \in \Sym_n} q_n^{f(\pi)}$.
We have
\[Z_n(q_n,f)=\sum_{\pi \in \Sym_n}  q_n^{a_n(\tilde f^o(\mu_\pi) +o(1))}=e^{o(a_n \log(q_n))} \sum_{\pi \in \Sym_n}  e^{a_n \log(q_n)\, \tilde f^o(\mu_\pi) }
= e^{o(n)} \sum_{\pi \in \Sym_n} e^{n \beta \tilde f^o(\mu_\pi)},\]
where, in the middle equality, we used the fact that the error in the approximation ${f(\pi)/a_n}=\tilde f^o(\mu_\pi)+o(1)$ is uniform in $\pi$ (see \eqref{eq:def_continuous_version}) and, in the last equality,
we used assumption (ii) and the boundedness of $\tilde f^o(\mu_\pi)$.
In particular
\begin{equation}\label{eq:approxA}
\frac1n \log\Big(\frac{1}{n!}Z_n(q_n,f) \Big) = \frac1n \log\Big(\expec \big(e^{n \beta \tilde f^o(\mu_{\pi_n})}\big) \Big) +o(1),
\end{equation}
where $\pi_n$ is uniformly distributed on $\Sym_n$. {We assume that $\beta\geq 0$, the case $\beta\leq 0$ is similar.}

\medskip

\noindent{\emph{\textbf{Step 1: The lower-bound in~\eqref{eq:estimate_logZ}.}}} Fix $\eps>0$.
Let $\mu_\star$ be a maximizer of  $H(\mu) + \beta \tilde f^o(\mu)$ on $\mathcal P$
(an upper-semicontinuous function on a compact space always has a maximizer).
Since $\tilde f^o$ is continuous, we can find an open neighborhood $\mathcal P_\star$
of $\mu_\star$ such that 
\begin{equation}\label{eq:def_Pstar}
 \mu \in \mathcal P_\star \ \Rightarrow\   \tilde f^o(\mu) \ge  \tilde f^o(\mu_\star) - \eps.
\end{equation}
Since $\mathcal P_\star$ is open, \Cref{thm:large_deviation} implies that, for $n$ large enough,
\begin{equation}\label{eq:bewibweio}
    \prob\big[\mu_{\pi_n} \in \mathcal P_\star \big] \ge \exp\Big(n \, \big( \sup\nolimits_{\mu \in \mathcal P_\star} H(\mu) -\eps \big) \Big) \ge \exp\Big(n \big(H(\mu_\star) -\eps\big)\Big).
\end{equation}
Using the two inequalities above, we have, again for $n$ large enough,
\begin{align}
    \expec \big(e^{n \beta \tilde f^o(\mu_{\pi_n})}\big)\geq \expec \big(e^{n \beta \tilde f^o(\mu_{\pi_n})}\mathbf{1}_{\{\mu_{\pi_n} \in \mathcal P_\star\}}\big) &\stackrel{\eqref{eq:def_Pstar}}{\ge}
    \prob\big[\mu_{\pi_n} \in \mathcal P_\star \big] \exp\big(n \beta (\tilde f^o(\mu_\star) - \eps) \big)\\ 
    &\stackrel{\eqref{eq:bewibweio}}{\ge}
    \exp \big(n  (H(\mu_\star) +  \beta \tilde f^o(\mu_\star) -  (\beta+1) \eps) \big).
\end{align}
Since $\eps>0$ is arbitrary, we have proved 
 \begin{equation}\label{eq:lowerA}
 \liminf_{n \to +\infty} \frac1n \log\Big(\expec \big(e^{n \beta \tilde f^o(\mu_{\pi_n})}\big) \Big)
\ge H(\mu_\star) + \beta \tilde f^o(\mu_\star) =  \max_{\mu \in \mathcal P} (H(\mu) + \beta \tilde f^o(\mu)).
\end{equation}

\medskip

\noindent{\emph{\textbf{Step 2: The upper-bound in~\eqref{eq:estimate_logZ}.}}} We use a 
 finite covering 
$\mathcal P=\bigcup_{i \in I} \mathcal P_i$ by closed sets such that, for all $i$ in $I$,
\begin{equation}\label{eq:def_Pi}
 \mu,\nu \in \mathcal P_i \ \Rightarrow\  |\tilde f^o(\mu)-\tilde f^o(\nu)|\le \eps.
\end{equation}
Such a finite closed cover exists since $\tilde f^o$ is continuous and $\mathcal P$ is compact.
We define 
\[
H_i=\sup_{\mu \in  \mathcal P_i} H(\mu)
\qquad\text{ and }\qquad 
F_i= \sup_{\mu \in  \mathcal P_i} \tilde f^o(\mu).
\]
Then we have 
\[\expec \big(e^{n \beta \tilde f^o(\mu_{\pi_n})}\big) \le \sum_{i \in I}
\prob\big[\mu_{\pi_n} \in \mathcal P_i \big] \exp(n \beta F_i) \le  \sum_{i \in I} \exp(n  (H_i+\beta F_i+o(1))),\]
where the second inequality follows from \Cref{thm:large_deviation} recalling that $\mathcal P_i$ is closed.
On the other hand, using \eqref{eq:def_Pi},  for all $\mu$ in $\mathcal P_i$, we have 
$H(\mu)+\beta \tilde f^o(\mu) \ge H(\mu)+\beta F_i -\beta \eps$. 
Maximizing over $\mu$ in $\mathcal P_i$
yields:
\[H_i +\beta F_i \le \sup_{\mu \in \mathcal P_i} (H(\mu)+\beta \tilde f^o(\mu))  +\beta \eps.\]
We finally obtain  
\begin{equation}\label{eq:upperA}
\limsup_{n \to +\infty}
 \frac1n \log\Big(\expec \big(e^{n \beta \tilde f^o(\mu_{\pi_n})}\big) \Big)
\le \max_{i \in I} (H_i+\beta F_i) \le  \max_{\mu \in \mathcal P}
 (H(\mu) + \beta \tilde f^o(\mu)) + \beta \eps.
\end{equation}
Since~\eqref{eq:upperA} holds for any $\eps>0$, and comparing with~\eqref{eq:lowerA} 
we conclude that
\[ \lim_{n \to +\infty} \frac1n \log\Big(\expec \big(e^{n \beta \tilde f^o(\mu_{\pi_n})}\big) \Big)
=  \max_{\mu \in \mathcal P} H(\mu) + \beta \tilde f^o(\mu).\]
Together with~\eqref{eq:approxA}, this proves~\eqref{eq:estimate_logZ}.

\medskip

\noindent{\emph{\textbf{Step 3: Possible limits.}}}
Consider now the second statement in the proposition, regarding potential limits of subsequences of $\sigma_{n,q_n}^f$. 
For this, let us call $E$ the set of $\mu$ realizing the above maximum, and let $C$ be a closed subset of $\mathcal P$ disjoint from $E$.
We have
\[ \prob\big[{\mu_{\sigma_{n,q_n}^f}} \in C\big] = \frac{1}{Z_n(q_n,f)} \sum_{\pi \in \Sym_n
\atop \mu_\pi \in C} q_n^{f(\pi)} .\]
The same proof as above shows that 
\[ \limsup_{n \to +\infty} \frac1n \log\left(\frac{1}{n!} \sum_{\pi \in \Sym_n
\atop \mu_\pi \in C} q_n^{f(\pi)}\right) \le  \max_{\mu \in C} H(\mu) + \beta \tilde f^o(\mu) <  \max_{\mu \in \mathcal P} H(\mu) + \beta \tilde f^o(\mu).\]
Hence, $\prob\big[{\mu_{\sigma_{n,q_n}^f}} \in C\big]$ tends to $0$. 
Since this holds for any closed subset $C$ of $\mathcal P$ disjoint from $E$,
this concludes the proof of the proposition.
\end{proof}

\begin{figure}[t]
    \centering
    \includegraphics[width=0.24\textwidth]{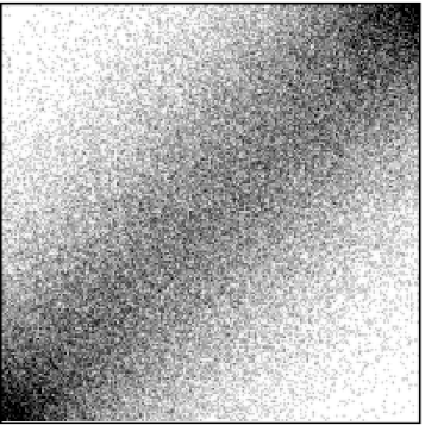}
    \includegraphics[width=0.24\textwidth]{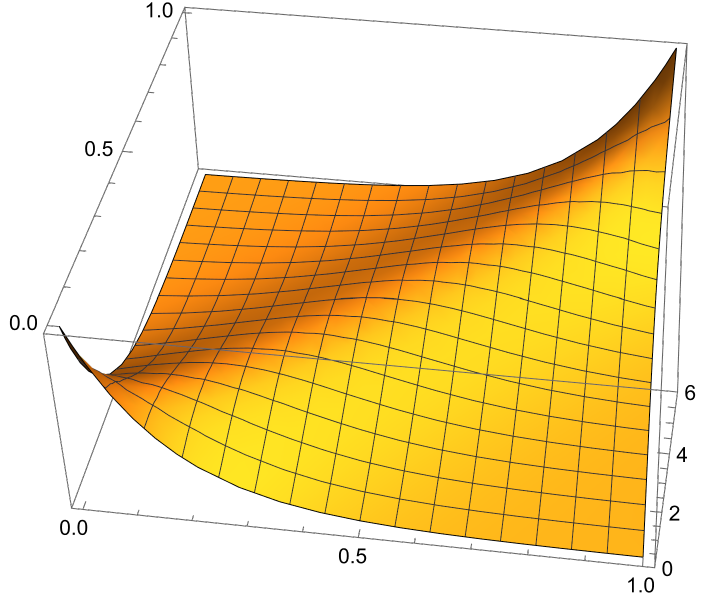}
    \includegraphics[width=0.24\textwidth]{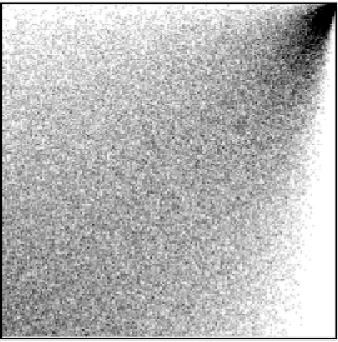}
    \includegraphics[width=0.24\textwidth]{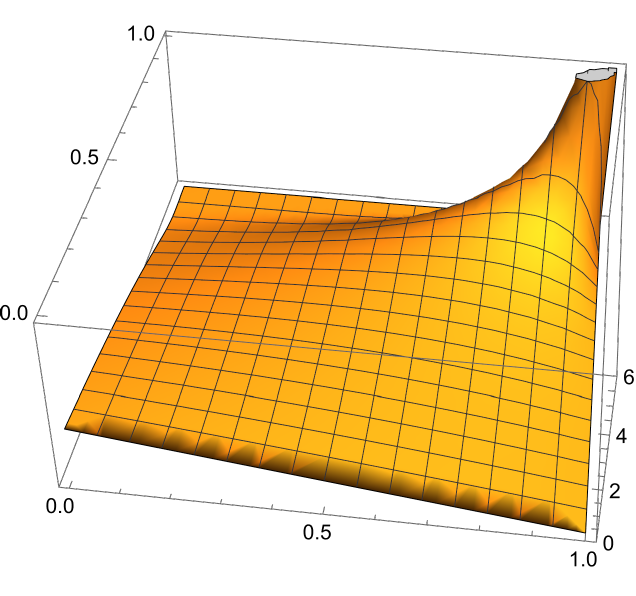}
    \caption{\textbf{From left to right:} 
    (1) A Mallows random permutation of size {$n=100,000$} and parameter {$q_n=0.99994$}. Note that $\frac{n}{2} \log(q_n) \approx -3$; 
    (2) The 3D graph of the function $\rho_\beta$ in \eqref{eq:Mallows-dens} for $\beta=-3$. Starr's convergence result predicts that the density of points on the left is approximately given by the function on the right. We see in particular that the densest area of points on the left is around the main diagonal and especially {near its endpoints}, where $\rho_\beta$ achieves its highest values. 
    (3) The diagram of a Luce-distributed permutation---{discussed later in Section~\ref{ss:analytic-limit-shapes}}---with Sukhatme weights of size $100,000$, see~\cite{borga2025luce} for more details; (4) The density $\rho(x,y)$ of the permuton limit of Luce-distributed permutations with Sukhatme weights.}
    \label{fig:mallows_luce}
\end{figure}

\begin{example}[Mallows permutations]
\label{ex:Mallows}
Random permutations defined by \eqref{eq:def_exp_biased} in the special case
when $f(\pi)$ is the number of inversions in $\pi$ 
are often referred to as \textbf{Mallows random permutations} and
have been widely studied in the last twenty years
 (see, e.g., \cite{dubach2024mallows} for recent results on the topic and multiple references).
 In this case, $f/\binom{n}{2} $ is the $21$ pattern density $d_{21}$
 and has a continuous version on $\mathcal P$, which we also denote by $d_{21}$.
 Therefore, if $\frac12 n \log(q_n)$ has a finite limit $\beta\in\mathbb{R}$, 
 \Cref{prop:limit_sigma_qf} applies.
 It turns out that the variational problem on the right-hand side of \eqref{eq:estimate_logZ} can be explicitly solved in this case,
 and that there is a unique maximizing permuton $\mu_\beta$, which has the following density with respect to Lebesgue measure:
\begin{equation}\label{eq:Mallows-dens}
    \rho_\beta(x,y)=\frac{\beta \sinh(\beta)}{\big[e^{-\beta/2} \cosh(\beta(x-y))-e^{\beta/2} \cosh(\beta(x+y-1))\big]^2}.
\end{equation}
By  \Cref{prop:limit_sigma_qf}, the Mallows random permutation $\sigma_{n,q_n}^f$
 converges in probability to the permuton $\mu_\beta$,
 as can be seen {in the first two panels of~\Cref{fig:mallows_luce}}.
  
 This convergence result was first obtained by Starr in \cite{starr2009thermodynamic} (see also~\cite{StarrWalters-uniqueness} for complementary work),
 who found and solved the variational problem above, but without going through the 
 large deviation principle of \Cref{thm:large_deviation}.
 More recently, Mallows random permutations with restricted positions were considered 
 by Gorin and Kenyon in connection with the 6-vertex model \cite{gorin2024restrictedMallows}, 
 leading to variants of the variational principle above and new limit shape results.
\end{example}

\begin{example}[Simply indexed permutation statistics]\label{exemp:simply-indexed}
A {\bf simply indexed permutation statistic} is a statistic of the form
\[f(\sigma) = \sum_{i=1}^n a_{n,i,\sigma(i)},\]
where, for each $n$, $(a_{n,i,j})_{1 \le i,j \le n}$ is a matrix of real numbers.
The study of such statistics dates back to Hoeffding \cite{hoeffding1951combinatorial}, 
who gave asymptotic normality conditions for $f(\sigma_n)$,
with $\sigma_n$ uniformly distributed on the symmetric group.
Here, following Mukherjee, we will assume that \[a_{n,i,j}=a(i/n,j/n)\]
 for some continuous function $a$ on the unit square.
 Then the resulting function $f$, divided by $a_n=n$, has a continuous version
 $\tilde f^o(\mu)=\int a(x,y) \mu(dx,dy)$.
 This is a linear function of $\mu$, so the function $H+\beta \tilde f^o$ is strictly concave.
 Hence, it reaches its maximum on $\mathcal P$ at a unique point $\mu^\star_{f,\beta}$.
 From \Cref{prop:limit_sigma_qf}, taking $q_n=q=e^\beta$,
 the random permutation $\sigma_{n,q}^f$ converges to $\mu^\star_{f,\beta}$; see~\cite[Theorem 1.4 (b)]{Mukherjee-estimation} for the original reference. 
 For concrete functions $a$, finding this maximum can be intractable,
 but Mukherjee provides an approximation algorithm for it \cite[Theorem 1.8]{Mukherjee-estimation}.
Moreover, if $a$ is twice differentiable and symmetric, it is possible to refine \eqref{eq:estimate_logZ} to find an asymptotic
for the partition function itself (and not only for the logarithm); the subexponential term then involves
a Fredholm determinant, see \cite{pal2024partition_Mallows,tripathi2026partition_Mallows}. {The fourth panel in Figure~\ref{fig:fig-intro-2} illustrates the case $a(x,y)=xy$ and $\beta=20$, known as the Spearman rank-correlation model, by showing the permuton associated with a permutation of size $n=10,000$ obtained by approximate MCMC sampling.}
{For $\beta>0$, taking $a(x,y)=-|x-y|^p$, with $p=1$ or $2$, gives Spearman's
footrule and squared-$L^2$ Mallows models, respectively, with probabilities proportional to
\[
	\exp\!\left(-\frac{\beta}{n^p}
	\sum_{i=1}^n|\sigma(i)-i|^p\right).
\]
For $p=2$, this is {equivalent to taking} $a(x,y)=xy$ as above with parameter $2\beta$. Zhong~\cite[Theorems~3.4 and~3.2]{zhong2021hitrun} analyzes hit and run samplers for these models and proves $O(\log n)$ total-variation mixing times at any fixed accuracy for every fixed $0<\beta<1/16$.}
\end{example}

\begin{figure}[t]
\includegraphics[height=5cm]{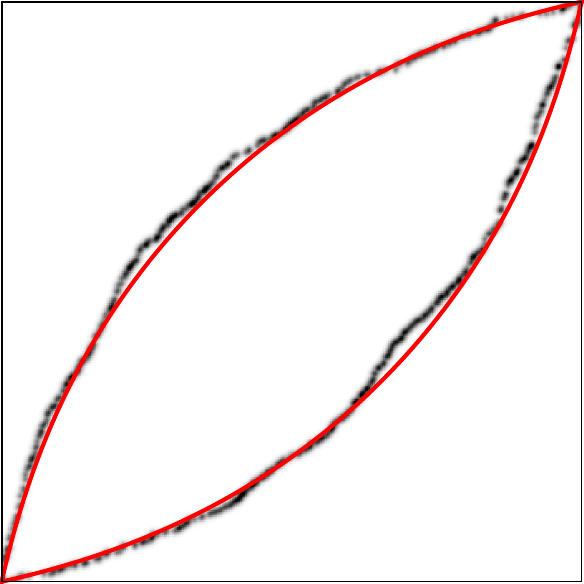}
\hspace{1cm}
\includegraphics[height=5cm]{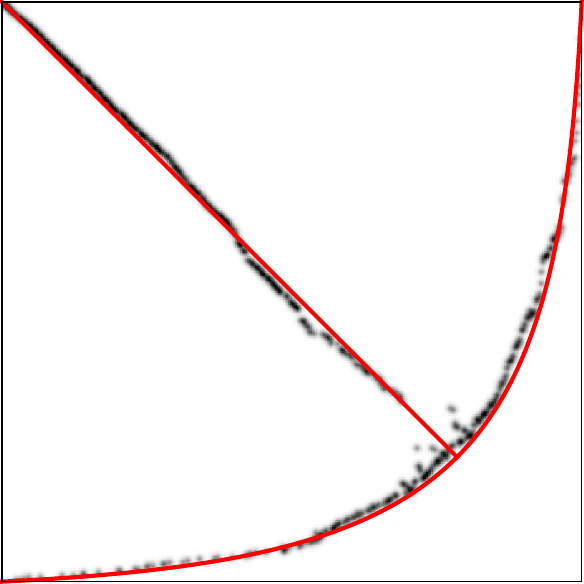}
\caption{Simulations of Mallows random permutations conditioned to avoid a pattern of length 3 ($321$ on the left, $231$ on the right) with $n=600$, $q_n=\exp(\beta/n)$ and $\beta=3$.
The red curves show the supports of the corresponding limiting permutons.}
\label{fig:ConstrainedMallows}
\end{figure}

\begin{example}[{Constrained Mallows permutations}]\label{expem:cond-mallows}
The same strategy---that is, using a large deviation principle to turn a limit shape problem for a biased measure into an analytic optimization problem---can be used for Mallows permutations constrained to avoid a pattern of length 3, thanks to the large deviation principle mentioned in Remark~\ref{rmk:further_LDP}.
Again here, we obtain a nontrivial permuton limit when the bias parameter is of the form
$q_n=\exp(\beta/n)$; see Fig.~\ref{fig:ConstrainedMallows}
and the paper~\cite{budzinski2026constrained_mallows} for more details.
\end{example}

\begin{example}[Major index distribution]
    {Another natural exponential bias, which falls outside the direct scope of Proposition~\ref{prop:limit_sigma_qf}, is given by the \emph{major index}
    \[
        \maj(\sigma)
        := \sum_{\substack{1\leq i<n\\ \sigma(i)>\sigma(i+1)}} i.
    \]
    For fixed $q\in(0,1)$, the major-index distribution assigns to $\sigma\in\mathfrak S_n$ probability proportional to $q^{\maj(\sigma)}$. Coopman~\cite{coopman2025major} gives a simple sampler for this model based on i.i.d.\ geometric random variables and identifies its limiting permuton. The support of this permuton is the closure of the union, over $i\geq 0$, of the line segments joining $(q^{i+1},0)$ to $(q^i,1)$. The major index depends on consecutive descents and is not continuous in the permuton topology, illustrating that nontrivial permuton limits can also arise for exponential biases not covered by the variational framework above.}
\end{example}

\subsection{Entropy-maximizing permutons, and an enumeration result}\label{sect:Entropy-maximizing-permutons}

Let us work in the framework of \Cref{prop:limit_sigma_qf}.
In the special case where $\tilde f^o$ is a pattern density $d_\tau$,
 any permuton $\mu_\star$ maximizing $H(\mu) +\beta \tilde f^o$ is, in particular, a maximizer of $H$
over permutons with the same density $d_\tau$. 
In formula, 
\[H(\mu_\star) = \max \{H(\mu), \mu \in \mathcal P \text{ s.t. }d_\tau(\mu)=d_\tau(\mu_\star)\}.\]
More generally, if $\tilde f^o$ is a linear combination of the pattern densities $d_{\tau_1}$, \dots, $d_{\tau_r}$,
then one needs to maximize $H$ over permutons with fixed values of the pattern densities $d_{\tau_1}$, \dots, $d_{\tau_r}$.

This maximizing problem also occurs in the following problem of asymptotic enumeration,
as shown by Kenyon, Král', Radin and Winkler \cite{ Kenyon-Permutations-densities}.

\begin{thm}[{\cite[Theorem 2]{Kenyon-Permutations-densities}}]\label{thm-entr-ken}
  Fix a collection $\tau_1$, \dots, $\tau_r$ of patterns and some desired densities $\alpha_1$, \dots, $\alpha_r$ for these patterns.
  For $\eps>0$ and $n \ge 1$, we let
  \[S_n^{ \bm \alpha,\eps}=\big\{\sigma \in S_n, |d_{\tau_i}(\sigma) - \alpha_i| < \eps \text{ for each } i \le r\big\}.\]
  Then
  \begin{equation}\label{eq:fixed_pattern_densities_asymptotic_enumeration}
 \lim_{\eps \to 0} \lim_{n \to +\infty} \frac1n \log \frac{\#S_n^{ \bm \alpha,\eps}}{n!} = \max_{\mu \in \mathcal P^{\bm \alpha}} H(\mu),
  \end{equation}
  where $\mathcal P^{\bm \alpha}=\big\{\mu \in \mathcal P,\ d_{\tau_i}(\mu) = \alpha_i \text{ for each } i \le r\big\}$.
  \label{thm:fixed_pattern_densities_asymptotic_enumeration}
\end{thm}
By the method of Lagrange multipliers, in the interior of the feasible region,
any solution $\mu_\star$ of the constrained maximization
problem on the right-hand side of \eqref{eq:fixed_pattern_densities_asymptotic_enumeration} is a critical point of
\[H(\mu) + \beta_1 d_{\tau_1}(\mu) + \dots + \beta_r d_{\tau_r}(\mu),\]
for some real numbers $\beta_1$, \dots, $\beta_r$.
Hence, the maximizing problems in \Cref{prop:limit_sigma_qf} and \Cref{thm:fixed_pattern_densities_asymptotic_enumeration}
are closely tied (with a nontrivial relation between $\alpha$ and $\beta$).
For $\tau=21$, as explained in \Cref{ex:Mallows}, the problem can be solved explicitly, and
the maximizer has an analytic density whenever the target proportion of $21$ is in $(0,1)$.
For larger patterns (or families of patterns), the problem is generally difficult.

We note that the theorem is meaningful only when $(\alpha_i)_{i \le r}$ belongs to the projection of the feasible region $\mathcal F$ on $\mathbb R^{\{\tau_1,\dots\tau_r\}}$ (recall its definition from \eqref{eq:defn-feas-reg});
otherwise, the set $S_n^{ \bm \alpha,\eps}$ is empty for large $n$ and small $\eps$,
and the set $\mathcal P^{\bm \alpha}$ is empty (so that the maximum in \eqref{eq:fixed_pattern_densities_asymptotic_enumeration}
is $-\infty$ by convention).

Numerical simulations performed by Kenyon, Král', Radin, and Winkler~\cite{ Kenyon-Permutations-densities} suggest that entropy-maximizing permutons at interior points of the feasible regions have analytic densities. Uniqueness, however, fails in general. Indeed, returning to the right-hand side of Figure~\ref{fig:feas-ref}, the inward corner of the feasible region (marked by a large dot) corresponds to two different singular permutons. \cite[Lemma 18]{ Kenyon-Permutations-densities} implies that along a non-empty line segment extending from the dot toward the origin, the maximum-entropy permuton looks like one of these two and in particular is not symmetric under the reflection $(x,y)\mapsto(x,1-y)$.  But this reflection exchanges the densities of $123$ and $321$ and preserves entropy; hence there are at least two maximum-entropy permutons at each point of the segment.

On the (relative) boundary of the feasible region, however, {most permutons seem to have} entropy $-\infty$,
so that \eqref{eq:fixed_pattern_densities_asymptotic_enumeration} roughly means
that the number of permutations with such pattern densities is super-exponentially smaller than $n!$.
This is consistent with the famous Stanley-Wilf conjecture---now Marcus-Tardos theorem \cite{marcus2004StanleyWilf}---which states that the number of permutations avoiding a given pattern grows at most exponentially
(avoiding a pattern is a special case of fixed pattern densities on the boundary of the feasible region),
but \Cref{thm:fixed_pattern_densities_asymptotic_enumeration} does not help much in this setting.

{
\subsection{Fluctuations beyond the law of large numbers}
	\label{ssec:fluctuations}
	The large deviation principles discussed above describe the
	exponential probabilities of macroscopic deviations. They do not,
	by themselves, determine the scale or distribution of typical
	fluctuations. Proposition~\ref{prop:sample_pattern_clt} gives a
	central limit theorem for samples from a fixed permuton; for
	other permutation models,} {fluctuation theorems have appeared 
	sporadically in the literature, and we review some results of this kind below.
We are however not aware of any general theory for fluctuations 
in the space of permutons;
see Open Problem \eqref{pb:fluctuations} in Section~\ref{ss:open_problems}.}
\medskip

{\emph{Mallows permutations.}
	For the model of Example~\ref{ex:Mallows}, Bufetov and
	Chen~\cite{bufetov2024local} prove a local central limit theorem
	for the height function
	\[
	H_n(x,y)=
	\#\bigl\{i\leq\lfloor nx\rfloor:
	\sigma_n(i)\leq\lfloor ny\rfloor\bigr\},
	\qquad x,y\in(0,1),
	\]
	when $q_n=1-c/n$ with $c>0$.
	In this regime, the permutations have a nontrivial permuton limit $\mu_{-c/2}$ (see Example~\ref{ex:Mallows}), and hence $H_n(x,y)$ is equivalent at first order
	to $n\mu_{-c/2}([0,x]\times[0,y])$.
	Bufetov and Chen proved that its fluctuations around
	$n\mu_{-c/2}([0,x]\times[0,y])$ are Gaussian on the
	$\sqrt n$ scale.}
	{A generalization of this result to random permutations constructed
	by recursive insertions of new values has been announced 
	in~\cite[Remark 2.15]{jasinska2025permuton_kcard_mallows}.
	This further result would also imply fluctuation results of pattern densities.
	In other regimes than $q_n=1-c/n$,
	 namely when $q_n$ in $(0,1)$ satisfies $n^{3/2}(1-q_n)\to0$ or 
	when $q_n=q\in(0,1)$ is fixed,
	Dubach~\cite{dubach2024mallows} establishes joint central limit
	theorems for pattern densities.}
\medskip

{\emph{Simply indexed Gibbs models.}
	In the setting of Example~\ref{exemp:simply-indexed}, fix $\beta\in\mathbb R$, let
	$\sigma_n=\sigma^f_{n,e^\beta}$, and write
	$\mu_\star=\mu^\star_{f,\beta}$.
	For a continuous function $g:[0,1]^2\to\mathbb R$, set
	$L_n(g)=\sum_{i=1}^n g(i/n,\sigma_n(i)/n)$.
	Duan, Lyu, and Powell~\cite[Corollary~2.15(ii)]{duan2026poisson}
	prove that
	\begin{align*}
		\frac{L_n(g)-\expec[L_n(g)]}{\sqrt n}
		\ \xrightarrow{\mathrm d}\ \mathcal N(0,V_{\mu_\star}(g)),\qquad
		V_{\mu_\star}(g)
		=\inf_{u,v\in L^2([0,1])}
		\int_{[0,1]^2}
		\bigl(g(x,y)-u(x)-v(y)\bigr)^2\,\mu_\star(dx,dy).
	\end{align*}
	The convergence holds jointly for finitely many test functions.
}
\medskip

{
{\em Pattern-avoiding and Grothendieck permutations.}
Anticipating a bit the results from next section, 
let us also mention fluctuation results for pattern avoiding permutations:
see e.g.~\cite{hoffman2017patternI} for shape fluctuations and 
\cite{janson2017av132,janson2019av321} for pattern density fluctuations
for random permutations avoiding patterns of length 3.
Remarkably, there exists also a model of random permutations,
called Grothendieck random permutations (see Section~\ref{ss:algebraic_limits})
for which fluctuations of Tracy-Widom type have been established for the height function
\cite{morales2025grothendieck}.
No fluctuation results for pattern counts seem to be known in this case.
}\medskip

\section{Some limit shape results}
\label{sec:limit-shape}

One of the main advantages of the language of permutons introduced in \Cref{sec:theory} is that it turns the vague question
\[
\text{\enquote{What does a large permutation look like?}}
\]
into a well-posed mathematical question about probability measures on the unit square.  In this section, we survey some of the answers that have been discovered in the last three decades for many different models of random permutations.  We use \emph{limit shape} in a broad sense: sometimes the limit is a deterministic permuton, sometimes it is a random permuton, and sometimes the permuton limit is only the first layer of a finer scaling picture. We present in this section a number of random permutation models, coming mostly from the enumerative/algebraic combinatorics literature, for which the permuton limit has been described. While uniform random permutations have the ``boring'' Lebesgue measure as the limiting permuton, the models from this section yield a remarkable collection of limiting objects with various properties: they can be random or deterministic, supported on a full square, or on (semi-)algebraic subsets of full or lower dimension. Some have self-similarity properties while others are simply unions of line segments, and so on.

\subsection{Recursive construction of random permutations, and deterministic analytic shapes}
\label{ss:analytic-limit-shapes}

The first explicit deterministic permuton limit for a widely studied non-uniform random permutation model is Starr's limit for the Mallows model; recall~\Cref{ex:Mallows}.

Recent work has considerably enlarged the list of deterministic analytic or semi-explicit limits arising from sequential biases.  Jasi\'nska and R\'ath introduced a model that interpolates between Mallows permutations and Travers's $k$-card-minimum model.  At each step, one removes a card from the remaining deck according to a prescribed distribution; under natural hypotheses, the associated permutations converge in probability to an explicit deterministic permuton~\cite{jasinska2025permuton_kcard_mallows}.  Another sequential model is the Luce model, in which cards are sampled without replacement with probabilities proportional to prescribed positive weights. Borga, Chatterjee, and Diaconis proved, under minimal assumptions on the weights, convergence to a permuton with an explicit density (see the right-hand side of~\Cref{fig:mallows_luce}), and they also obtained limiting pattern densities and local limits \cite{borga2025luce}.  These models show that entropy-tilting is not the only mechanism leading to tractable deterministic permutons: sequential sampling rules can also have explicit permuton limits.

A different bias, closer in spirit to records and algorithms, is the record-biased model. Here the probability of a permutation is tilted by a constant $\theta$ raised to the number of records; see the second panel in Figure~\ref{fig:fig-intro-2}. As in the above example, such random permutations can also be constructed by successive random \enquote{insertions} in the permutation. Bouvel, Nicaud, and Pivoteau show that in the regime where the bias creates a linear number of records, the random permutation converges to a deterministic permuton which they characterize fully \cite{bouvel2026record_biased}.

Some interesting limit shapes may also arise not by changing the probability measure directly, but by applying a deterministic transformation to a uniform random permutation.  A beautiful example is the runsort operator.  Given a permutation $\pi$, split it into its ascending runs and then sort these runs lexicographically, or equivalently, by their first entries.  This elementary operation turns a uniform random input into a highly non-uniform output; see the third panel in Figure~\ref{fig:fig-intro-2}. 
\begin{thm}[The runsort permuton \cite{alon2021runsort}] \label{thm:runsort_permuton} 
Let $\pi_n$ be uniform in $\Sym_n$.  Then
\[
        \mu_{\operatorname{runsort}(\pi_n)}
        \xrightarrow[n\to\infty]{a.s.}
        \mu^{\operatorname{runsort}},
\]
where $\mu^{\operatorname{runsort}}$ is an explicit deterministic permuton. In particular, its support is
\[
        \big\{(x,y)\in[0,1]^2: x\le y e^{1-y}\big\}.
\]
\end{thm}
The runsort limit is a useful reminder that deterministic algorithms on permutations can produce nontrivial continuum limits for the rescaled permutation diagram.

\subsection{Pattern avoidance: from collapsed shapes to random fractal measures}\label{sect:pattern-av-perm}
{For a set $B$ of patterns, $\operatorname{Av}(B)$ consists of the permutations containing no subsequence order-isomorphic to a member of $B$. A classical algorithmic example is Knuth's characterization of permutations sortable into increasing order by a single stack: when the input is read from left to right, these are exactly the $231$-avoiding permutations \cite{knuth1997fundamental}. All six patterns $\tau$ of size three have the same avoidance counts, $|\operatorname{Av}(\tau)\cap\Sym_n|=\frac{1}{n+1}\binom{2n}{n}$. Equality of avoidance counts for every $n$ is called \emph{Wilf equivalence}; it does not assert equality of the distributions of the numbers of occurrences. Patterns of size four already split into three Wilf classes, represented by $1234$, $1342$, and $1324$. {See \cite{bona_book_permutations,vatter2015permutation_classes} for surveys of enumerative results on permutation classes}. For permuton theory, {uniform random permutations in} avoidance classes provide natural constrained models whose limiting geometries depend strongly on the forbidden patterns: the examples below range from collapsed deterministic shapes to random fractal measures. These geometric questions complement the classical enumeration of the classes.}

Pattern avoidance provides perhaps the richest collection of fascinating permuton limits.  The first global limit shapes in this direction were obtained for permutations avoiding a single pattern of size three.  Miner and Pak initiated the study of the shape of random pattern-avoiding permutations, describing the asymptotic positions of individual entries and the corresponding geometric pictures~\cite{miner-pak-permuton}.  Madras and Pehlivan proved refined structural estimates for random $312$-avoiding permutations, and in particular, exponential concentration of points near the diagonal after applying symmetries~\cite{madras2016structure}.  Hoffman, Rizzolo and Slivken then identified the finer fluctuation scale in terms of Brownian excursion~\cite{hoffman2017patternI}.  In permuton language, the global limit for these single length-three avoidance classes is typically degenerate: after a symmetry, all the mass collapses to the diagonal.  The interesting randomness is only visible at smaller scales.

This collapsed behavior persists in several monotone-avoidance models. Hoffman, Rizzolo, and Slivken studied permutations avoiding a long decreasing sequence, equivalently permutations with a bounded longest decreasing subsequence, and related the fluctuations of the increasing layers\footnote{Every permutation avoiding the pattern $k\,k{-}1\dots1$ can be naturally decomposed into $k{-}1$ increasing subsequences.} to the eigenvalue process of a traceless Hermitian Brownian bridge~\cite{hoffman2019scaling}.  Hohmeier and Slivken recently used the language of permutons, together with Backelin--West--Xin type bijections, to prove anti-diagonal limits for further classes avoiding a single pattern \cite{hohmeier2025single_pattern}.  These results are valuable even when the limiting permuton is simple: they identify complicated mechanisms that force such collapse and clarify which information is lost in the permuton limit. Finally, affine and periodic variants provide another direction.  Madras and Troyka studied bounded affine permutations avoiding a fixed decreasing pattern.  Their scaling limit result is not literally a permuton theorem on the unit square, because the model has periodic boundary conditions, but it is permuton-like: a typical permutation avoiding the monotone decreasing pattern of fixed size $k$ converges to a union of $k{-}1$ random lines of slope $1$, satisfying a global sum constraint~\cite{madras2021bounded_affine}.

{
These examples fit into a general rigidity theorem for pattern-avoiding permutons. Garbe, Hladk\'y, Kun and Pek\'arkov\'a proved that if $d_\tau(\mu)=0$ for a pattern $\tau$ of size $k$, then, after disintegrating $\mu$ with respect to the first coordinate, almost every fiber is supported on at most $k{-}1$ points; moreover, this bound is sharp \cite{garbe2024avoiding}. In particular, every pattern-avoiding permuton is singular with respect to two-dimensional Lebesgue measure. 
}

\begin{figure}[t]
    \centering
    \includegraphics[width=0.19\textwidth]{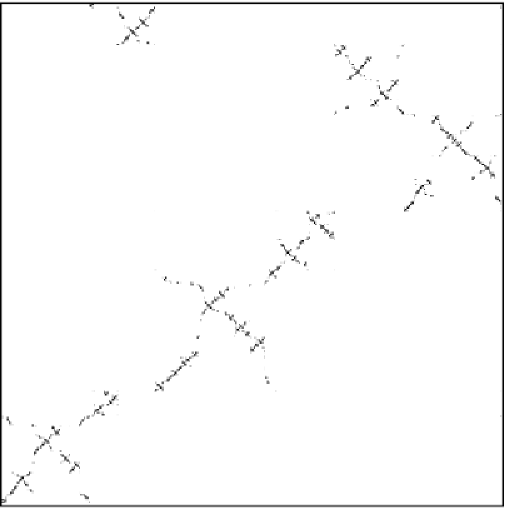}
    \includegraphics[width=0.19\textwidth]{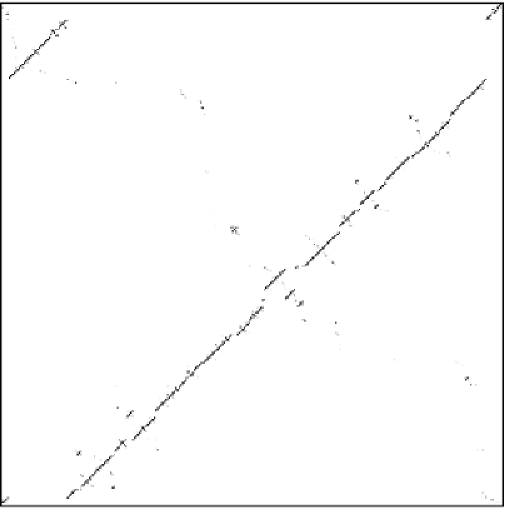}
    \includegraphics[width=0.19\textwidth]{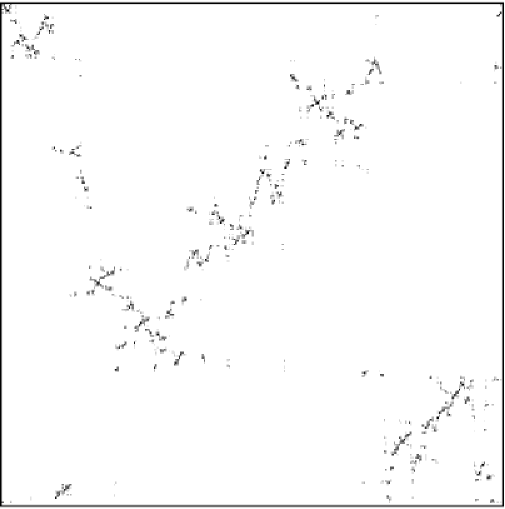}
    \includegraphics[width=0.19\textwidth]{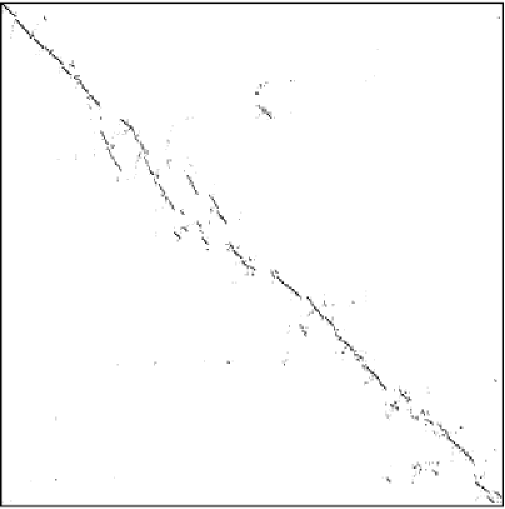}
    \includegraphics[width=0.19\textwidth]{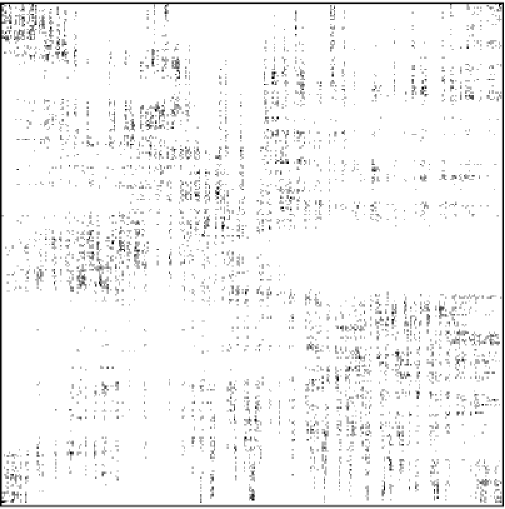}
    \caption{\textbf{From left to right:} Permutation samples $\sigma_n^{\bm{\mu}_{\rho,q}}$ of size $n=10,000$ from the skew Brownian permutons $\bm{\mu}_{\rho,q}$ with parameters $(\rho,q)=(1,0.5),(1,0.2),(-0.5,0.5),(-0.5,0.95),(-0.995,0.5)$. In particular, the first sample corresponds to the permuton limit of uniform separable permutations and the third sample corresponds to the permuton limit of uniform Baxter permutations. Recall that when $\rho=1$ (this is the case in the first two simulations), then $\bm{\mu}_{1,q}$ is equal to the Brownian separable permuton $\bm{\mu}_{1-q}$. When $(\rho,q)=(-1/2,1/2)$ (this is the case in the third simulation), then $\bm{\mu}_{\rho,q}$ is the Baxter permuton, i.e.\ the permuton limit of uniform Baxter permutations.}
    \label{fig:pat-av-permutons}
\end{figure}

{Despite this rigidity result, pattern avoidance does not force the global geometry of the limiting permuton to be simple or deterministic. When the corresponding permutation class admits a recursive substitution structure, highly nontrivial random fractal limits can arise. A fundamental example is} the class of separable permutations $\operatorname{Av}(2413,3142)$; see for a beautiful account~\cite{ghys2017singular}.
{Equivalently, separable permutations are generated from the singleton permutation $1$ by repeated direct sums and skew sums. For permutations $\alpha$ and $\beta$, the \emph{direct sum} $\alpha\oplus\beta$ places the diagram of $\beta$ to the right of and above that of $\alpha$, while the \emph{skew sum} $\alpha\ominus\beta$ places it to the right of and below that of $\alpha$. We return to separable permutations and their Brownian limits in \Cref{sect:univ-perm-skew}.}
Bassino, Bouvel, F\'eray, Gerin and Pierrot proved that a uniform separable permutation converges in distribution to the \textbf{Brownian separable permuton}, a random permuton constructed from a Brownian continuum random tree (CRT) decorated by independent $\oplus/\ominus$ signs \cite{bassino-separable-permuton}.  This was one of the first genuinely \emph{random} permuton limits.  It is fractal, singular with respect to Lebesgue measure (see the first diagram in \Cref{fig:pat-av-permutons}), but nevertheless sufficiently explicit to compute expected pattern densities and intensity measures (Section~\ref{subsec:intensity-measures}).

{More generally, for $\sigma\in\Sym_k$ and nonempty permutations $\pi_1,\ldots,\pi_k$, the \emph{inflation} $\sigma[\pi_1,\ldots,\pi_k]$ replaces each entry $\sigma(i)$ by a block order-isomorphic to $\pi_i$. Each block occupies consecutive positions and consecutive values, and the relative order of the blocks is prescribed by $\sigma$. Thus $\alpha\oplus\beta=12[\alpha,\beta]$ and $\alpha\ominus\beta=21[\alpha,\beta]$. A permutation class is \emph{substitution closed} if it contains all such inflations whenever $\sigma$ and $\pi_1,\ldots,\pi_k$ belong to the class \cite[Definitions~1.2--1.3]{Bassino-universal-perm}.}

The Brownian separable permuton $\bm{\mu}_{1/2}$ is just an instance of a larger universality class of permutons. Bassino, Bouvel, F\'eray, Gerin, Maazoun, and Pierrot proved that  a one-parameter biased version of the Brownian separable permuton, denoted by $\bm{\mu}_{p}$ with $p\in(0,1)$, is the permuton limit for a large family of substitution-closed permutation classes \cite{Bassino-universal-perm}; see the first two diagrams in~\Cref{fig:pat-av-permutons}.  For classes with finite combinatorial specifications, a later dichotomy theorem of the same authors shows that, under some connectivity hypothesis of the underlying specification, the limiting permuton is either a deterministic $X$-shaped permuton or a Brownian separable permuton \cite{bassino2019scaling}. Thus, the same recursive decomposition can produce either a deterministic macroscopic ``skeleton measure'' or a random Brownian fractal measure, depending on the analytic type of the specification.

{In a related heavy-tailed substitution framework, Stufler~\cite{stufler-26-dirichlet} introduced Poisson--Dirichlet permutons (and their graphon analogues). These permutons are constructed by ``assembling'' random permuton blocks whose masses are governed by a two-parameter Poisson--Dirichlet process. He proved a general invariance principle showing that they arise as universal limits and obtained concrete examples from iterated constructions based on separable permutations.}

There are also random recursive variants.  F\'eray and Rivera-Lopez introduced a Markovian model of random separable permutations and proved almost sure convergence to a random recursive separable permuton \cite{feray2023recursive_separable}.  This limit satisfies a distributional fixed-point equation, its expected pattern densities and intensity measure are given by explicit formulas, and its distribution is singular with respect to that of the Brownian separable permuton.  It is therefore a useful warning that similar recursive combinatorics may lead to quite different continuum random measures. Moreover, we will see shortly that the recursive separable permuton belongs to the closure of a class of permutons generated by random walks.

Baxter permutations form a bridge between pattern avoiding permutations and random planar geometry.  Dokos and Pak computed expected shapes for doubly alternating Baxter permutations and raised questions about the shape of large Baxter permutations \cite{ dokos-pak-bipolar-permuton}.  Borga and Maazoun later constructed the Baxter permuton and proved that it is the permuton limit of uniform Baxter permutations \cite{bm-baxter-permutation}; see the third diagram in~\Cref{fig:pat-av-permutons}.  Borga then introduced the two-parameter family of \textbf{skew Brownian permutons} $\bm{\mu}_{\rho,q}$, with $(\rho,q)\in(-1,1]\times(0,1)$, and proved that strong-Baxter and semi-Baxter permutations converge to random permutons of this family \cite{borga-skew-permuton,borga2021permuton}. The skew Brownian permutons contain both the Brownian separable permutons $\bm{\mu}_{p}$ (when $\rho=1$ and $q=1-p$) and the Baxter permuton (when $\rho=-1/2$ and $q=1/2$) as special cases, and they are discussed in much more detail in the later~\Cref{sect:univ-perm-skew}.  The skew Brownian  permutons are closely related to \emph{Liouville quantum gravity} surfaces and \emph{Schramm--Loewner evolution} curves, as discussed later in~\Cref{sect:space-filling-constr}.

\medskip

{To conclude this section on pattern-avoiding permutations, let us discuss permutation classes defined by geometric constraints. A large family is given by {\em monotone grid classes}: we are given a finite rectangular table of monotonicity constraints (increasing, decreasing, or empty), and a permutation belongs to the class if its diagram can be split into a corresponding rectangular grid satisfying these constraints. For some monotone grid classes, \cite{alshammari2025grid} proves convergence to a deterministic permuton supported on finitely many line segments, with masses determined by an analytic optimization problem.
	
Related to monotone grid classes, but more closely associated with permutons, are the {\em geometric grid classes} which require the points to be drawable on prescribed diagonal line segments~\cite{albert2013geometric}. Thus a geometric grid class consists, for each $n$, of the support of {permuton samples} $\sigma^\mu_n$ from certain permutons $\mu$.  There is a natural correspondence between {monotone and geometric} grid classes (each being defined by a $\{0,+1,-1\}$-valued matrix), and indeed for some matrices the classes agree; generally the geometric grid class is smaller.
	
Members of the ``permutation patterns" research community are naturally interested in the {\em asymptotic growth rate} $gr(C)$ of a permutation class $C$, that is, $\lim_{n \to \infty} |C \cap \Sym_n|^{1/n}$; this number (which exists for geometric grid classes, as a consequence of results in ~\cite{albert2013geometric}) {might be} related to the Shannon entropy $h(\sigma^\mu_n)$ of {permuton samples}.  The entropy $H$ defined in Section~\ref{sec:LDP}, which yields powerful large deviation theorems, is useless for singular permutons; instead we define the \textbf{catropy}\footnote{From the Greek root ``cata" as in ``catacomb," meaning ``under".} of a permuton $\mu$ by $C(\mu) = \lim_{n \to \infty} \frac1n h(\sigma^\mu_n)$, when this limit exists.  (See~\cite{maga2025entropy} for a general discussion of the behavior of $h(\sigma^\mu_n)$, including, e.g., a theorem giving the value of $C(\mu)$ when $\mu$ is supported on the graph of a piecewise-continuously-differentiable measure-preserving function with finitely many pieces.)	
	
Since the Shannon entropy of a probability distribution supported on
$k$ objects is at most $\log_2 k$, we have $h(\sigma^\mu_n)\leq \log_2|\mathcal C\cap\Sym_n|$, and hence $C(\mu)\leq \log_2\operatorname{gr}(\mathcal C)$ whenever the catropy exists. Computing and optimizing catropy for geometric grid permutons, and understanding when this bound is sharp, leads to the questions formulated in Problem~\eqref{pb:catropy}.}

A related model is that of square permutations, i.e.~permutations in which every point is a record from one of the four directions. Here the points of the permutations must live on four monotone subsequences, but there is no grid structure. The limiting permuton turns out to be a \emph{random} rectangle with sides of slope $\pm1$ \cite{MR4149526, borga2021almost}. In particular, we see a fundamental difference from connected monotone grid classes, which have deterministic permutons. The limiting object is nevertheless much simpler than the Brownian fractal permutons discussed above, but it illustrates that the randomness of the finite model may survive in the permuton limit.

\subsection{Permutations generated by random walk}
\label{ssec:random-walk}

A natural source of random permutations is the symmetric random walk.  Let $\mu$ be a symmetric, atomless
probability distribution on the real line, and let $0 = X_0,X_1,X_2,\dots,X_n$ be a random walk
of $n$ steps with step distribution $\mu$.  Let $\pi^\mu$ be the permutation associated with the points $\{(i,X_i),\, 0\le i\le n\}$ as in Section~\ref{ssec:samples_moments_convergence},
that is, the permutation of $\{0,1,\dots,n\}$ for which $X_{\pi^{-1}_0} < X_{\pi^{-1}_1} < \dots < X_{\pi^{-1}_n}$.
A classical theorem of Darling \cite{Dar51} states the surprising and not-so-well-known fact that the
one-point probabilities $\mathbb P(\pi^\mu(i) = j)$ are independent of $\mu$.

\begin{figure}[t!]
	\begin{center}
		\includegraphics[width=0.30\textwidth]{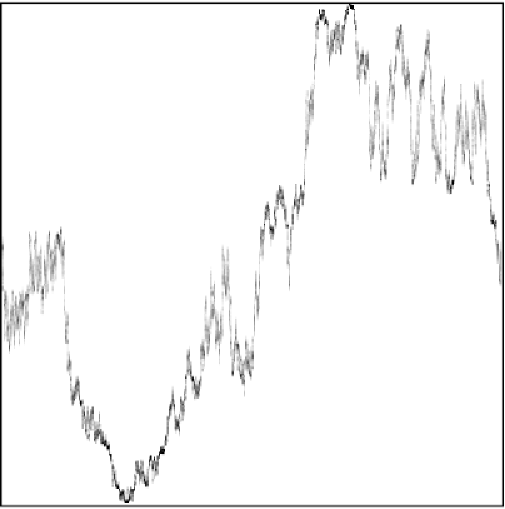}  \qquad \includegraphics[width=0.35\textwidth]{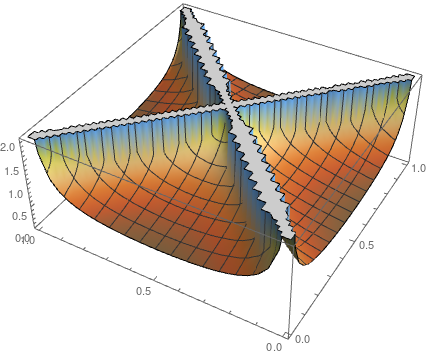} 
		\caption{\textbf{Left:} Sample of the permuton associated with the random permutation $\pi^\mu$ associated with a random walk of size $n=250,000$.
		\textbf{Right:} The density of the intensity measure of the Wiener permuton and of the recursive separable permuton. The small oscillations and truncation in the right panel are numerical artifacts of the plot and do not reflect genuine features of the density.\label{fig:wiener}}
	\end{center}
\end{figure}

Francis and Winkler~\cite{FW26} recently generalized this result as follows.
Let $A=\{a_1,\dots,a_n\}$ be a finite set of real numbers with distinct subset sums.
Let $(\sigma,\eps_1, \dots,\eps_n)$ be a uniform random
permutation of size $n$ together with $n$ uniform i.i.d.\ symmetric signs $\eps_i\in\{-1,1\}$, and set 
\[
    X_k=\sum_{i=1}^k  \eps_i a_{\sigma(i)}.
\]
Then we consider the random permutation $\pi^A$ associated with $(X_0,\dots,X_n)$ as above.
Again the one-point probabilities $\mathbb P(\pi^A(i) = j)$ are independent of $A$.
Taking $A$ to be an $n$-sample of law $\mu$ recovers the theorem of Darling stated above.
When $A = \{1,2,4,\dots,2^{n-1}\}$,
the random permutation $\pi^A$ has the same law as
 F\'eray and Rivera-Lopez's recursive separable permutations \cite{feray2023recursive_separable}.
 In particular, in this case $\mu_{\pi^A}$ converges to the so-called recursive separable permuton.  Since one-point probabilities are universal, and, recalling \eqref{eq:intensity-one-point-probabilities} in Section~\ref{subsec:intensity-measures}, for any bounded continuous $f$,
\begin{equation}\label{eq:rel-int-meas-proba}
	\mathbb{E}\!\left[\int f\,d\mu_{\pi^A}\right]
	=
	\frac{1}{n}\sum_{i,j}
	f\!\left(\frac{i}{n},\frac{j}{n}\right)
	\mathbb{P}\bigl(\pi^A(i)=j\bigr)\, +\, o(1),
\end{equation}
we deduce that for any $A$ (and thus, in particular, for any symmetric random walk), the intensity measure $\mathbb{E}[\mu_{\pi^A}(\cdot)]$ of $\mu_{\pi^A}$ converges to that of
the recursive separable permuton described in \cite{feray2023recursive_separable} and drawn on Fig.~\ref{fig:wiener} (right); this extends the classical \emph{arcsine law} of \cite{And49,EK47}.

Coming back to the random walk setting, if the step distribution has finite variance (symmetry is no longer needed here, just zero mean), \cite{FW26}
shows that the resulting permutation approaches a random permuton $\mathbb{W}$ called the \emph{Wiener permuton}
(see Fig.~\ref{fig:wiener}, left), which is nothing more than the rank-uniformized standard Brownian motion.  In fact,
the standard construction of Brownian motion can be made to produce a sequence of random walk
permutations that converge, not merely in distribution, to that permuton.  Also constructed is an interpolating
family of random permutons, all with the same intensity measure, for
which the zeroeth's limit is $\mathbb{W}$ and the ``infinityeth's" is the recursive separable permuton; see~\cite[Section 5]{FW26}.

Permutations generated by random walk when $\mu$ is 
the (symmetric) Laplace distribution play a special role in this work; 
in this case, we have the remarkable identity 
\[
\mathbb P[\pi^\mu=\tau]=\mathbb E[d_{\tau}(\mathbb{W})], \quad\text{for any pattern $\tau$ of $\{0,1,\dots,n\}$}.
\] 
A curious consequence is that the limiting \emph{consecutive} pattern
densities for the Laplace random walk permutations coincide with their expected limiting ``classical" pattern
densities.  As far as we know, this is the first nontrivial connection that
has been observed between the densities of the same pattern interpreted in both ways.

\subsection{Permutation-valued processes}
\label{ssec:processes}

The permuton topology also makes it possible to discuss limits of entire permutation-valued processes. The fundamental example is the random sorting network. A sorting network of size $n$ is a shortest path from the identity permutation to the reverse permutation using adjacent transpositions; its length is $N=\binom n2$. Angel, Holroyd, Romik and Vir\'ag discovered a remarkable conjectural picture: the values should move approximately along sine curves, and the permutation matrix at time $1/2$ should converge to the projection of the surface measure on the sphere~\cite{angel2007sorting_networks}. These conjectures were later proved by Dauvergne~\cite{dauvergne2018archimedean}.

\begin{figure}[t]
	\[
	\includegraphics[width=0.19\textwidth]{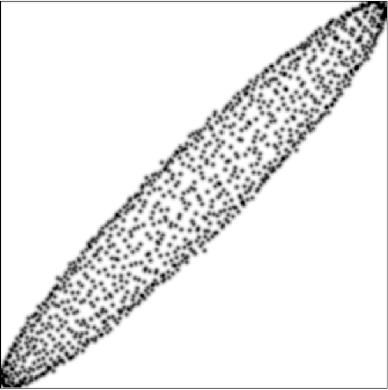}
	\hspace{0.2cm}
	\includegraphics[width=0.19\textwidth]{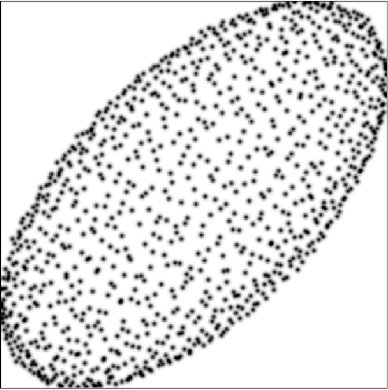}
	\hspace{0.2cm}
	\includegraphics[width=0.19\textwidth]{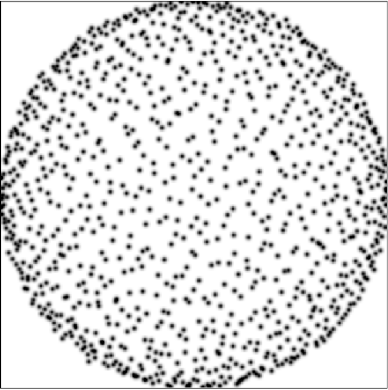}
	\hspace{0.2cm}
	\includegraphics[width=0.19\textwidth]{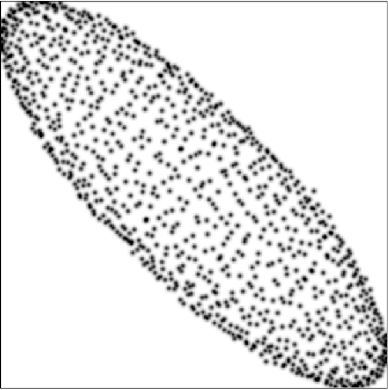}
	\hspace{0.2cm}
	\includegraphics[width=0.19\textwidth]{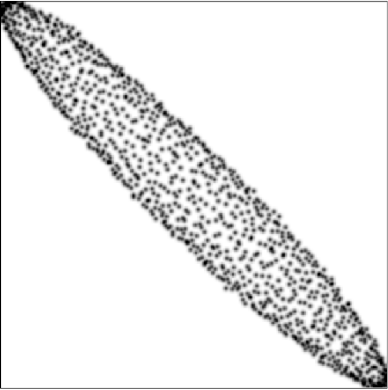}
	\]
	\caption{{\textbf{From left to right}: The permutons associated with the permutations $\sigma^n_{\lfloor tN\rfloor}$ along a uniform random sorting network of size $n=1000$, for $t=0.1,0.3,0.5,0.8,0.9$, where $N=\binom{n}{2}$. The corresponding deterministic limits $\mu_t^{\mathrm{Arch}}$ are described in Theorem~\ref{thm-sorting}; at $t=1/2$, the limit is the Archimedean permuton.} \label{fig:sorting-process}}
\end{figure}

Let $(A,B)$ have the \textbf{Archimedean law}, the projection of the uniform surface measure on the unit sphere to the first two coordinates. For $t\in[0,1]$, define the permuton
\[
\mu_t^{\operatorname{Arch}}
=\mathcal L\left(
\frac{1+A}{2},
\frac{1+A\cos(\pi t)+B\sin(\pi t)}{2}
\right).
\]
At $t=0$ this is the identity permuton, at $t=1$ the reverse permuton, and at intermediate times its support is an ellipse (with its interior); at $t=1/2$ one obtains the spherical, or Archimedean, permuton; see also Figure~\ref{fig:sorting-process}.

\begin{thm}[Archimedean sorting-network limit \cite{dauvergne2018archimedean}]\label{thm-sorting}
	Let $(\sigma^n_k)_{0\le k\le N}$ be a uniform random sorting network of size $n$, viewed as a process of permutations. Then, for each fixed $t\in[0,1]$,
	\[
	\mu_{\sigma^n_{\lfloor tN\rfloor}}
	\xrightarrow[n\to\infty]{\mathbb{P}}
	\mu_t^{\operatorname{Arch}}.
	\]
	Moreover, the trajectories of the permutation values are, uniformly with high probability, close to sine curves.
\end{thm}

{Remarkably, the same sine-curve process mentioned above had already appeared in Brenier's 1989 work on generalized incompressible flows~\cite[Proposition~6.3]{brenier1989least}. After rescaling, his explicit generalized solution of the incompressible Euler equations is the process
	\[
	X_t=\frac{1+A\cos(\pi t)+B\sin(\pi t)}{2},
	\qquad t\in[0,1],
	\]
with $(A,B)$ distributed according to the Archimedean law. This process minimizes the expected kinetic action subject to having uniform one-time marginals and satisfying the endpoint condition $X_1=1-X_0$. The uniform-marginal condition expresses incompressibility, so Brenier's work can be viewed as an early study of permuton processes, long before this terminology was introduced. His later paper~\cite[Figures~3 and~6]{brenier2008generalized} presents numerical approximations of this flow whose plots closely resemble random sorting networks. As Dauvergne explains in~\cite[Section~1.5]{dauvergne2018archimedean}, this connection was brought to his attention only after the first version of his sorting-network proof had appeared.}

The result in Theorem~\ref{thm-sorting} should be considered together with the variational theory of permutation-valued paths developed by Rahman, Vir\'ag and Vizer~\cite{rahman2019geometry}. They showed that the Archimedean path $(\mu_t^{\operatorname{Arch}})_{t\in[0,1]}$ is the unique minimum-energy path from the identity permuton to the reverse permuton in a natural metric on permutons; related work of Dauvergne and Vir\'ag established strong elliptic support estimates on fixed-time limits~\cite{dauvergne2020circular}. Kotowski and Vir\'ag~\cite{KV22} further connect this variational picture with probability and fluid dynamics{, building on Brenier's theory of generalized incompressible flows and its refinement by Ambrosio and Figalli~\cite{ambrosio2009geodesics}}. They prove large-deviation bounds for the interchange process on an interval in the space of permuton processes, with Dirichlet energy governing the rate, and use them to derive the Archimedean limit and an asymptotic enumeration of {\em relaxed} sorting networks.

{Another process-level analogue of the sorting-network picture is provided by the continuous-time Mallows process studied by Adamczak and Kotowski~\cite{AK25}. They prove both a global limit, in which suitably rescaled particle trajectories converge to an explicit process described in terms of Starr's Mallows permuton, and a local limit given by a permutation-valued process on $\mathbb Z$ with Mallows marginals.}

{A different process-level direction comes from up-down chains.}  An up-down chain first enlarges a combinatorial object and then projects back to the original size.  F\'eray and Rivera-Lopez~\cite{feray2025up_down} identified integrable up-down chains on permutations and proved convergence, after time scaling, to permuton-valued Feller diffusions.  These limiting diffusions have explicit diagonalizable semigroups, are ergodic, and have the random recursive separable permutons as stationary distributions.  This is a different kind of limit shape theorem: not a single static permuton and not a deterministic path, but a Markov process whose states are themselves \emph{random} permutons.

\subsection{Algebraic and integrable-probability limits}\label{sect:algebraic-and-integrable}
\label{ss:algebraic_limits}
A number of interesting models of random permutations also emerge from the algebraic combinatorics literature. The Erd\H{o}s--Szekeres extremal permutations, discussed in the introduction (Section~\ref{ssec:before}), which can be obtained by applying the Robinson-Schensted correspondence to a pair of independent uniform random tableaux of square shape, are one example~\cite{pittel2007limit_shapes_square_tableaux,MR2266895}. Their limit shape is a permuton with an analytic density supported on the interior of an algebraic curve of degree $4$; see the fifth panel in~\Cref{fig:fig-intro-2}. More generally, it would be interesting to have a systematic understanding of the link between the shape of the  Robinson-Schensted tableaux of a permutation, and its permuton limit; see the open problem section for further discussion on this topic.

The most recent additions to the list of explicit permuton limits come from algebraic combinatorics, especially Schubert and Grothendieck calculus.  The combinatorial objects here are pipe dreams, or related diagrams encoding polynomials representing cohomology and $K$-theory classes of flag varieties.  A pipe dream determines a permutation, so a natural random pipe dream produces a random permutation.  What is striking is that the macroscopic limit can be analyzed using interacting particle systems.

Morales, Panova, Petrov, and Yeliussizov~\cite{morales2025grothendieck} studied Grothendieck random permutations, where the probability of a permutation is proportional to the corresponding Grothendieck polynomial at parameter $\beta=1$, with all variables specialized to the same value.
  By mapping the model to a version of TASEP, they obtained the limiting permuton and proved $n^{1/3}$ Tracy--Widom GUE fluctuations {for the height function}, placing this algebraic model in the Kardar--Parisi--Zhang (KPZ) universality class; see the first panel in Figure~\ref{fig:fig-intro-2}.
  {A related model, based on principal specializations of Schubert polynomials
  instead of Grothendieck polynomials, was later studied
  by Anderson, Panova and Petrov in~\cite{anderson2026SchubertSpecializations}. 
 The pictures share some similarities with that of Grothendieck random permutations,
 but determining the permuton limit for such random permutations is an open problem.}

Defant~\cite{defant2025demazure} subsequently constructed several new families of permutons from the Demazure product.  One part of the work extends the Grothendieck pipe-dream limits from staircase shapes to arbitrary order-convex shapes; another studies deterministic sequences of bubble-sort operators applied to random initial permutations.  A key conceptual point is a min-plus, or tropical, formulation of the Demazure product, which makes it possible to define a Demazure product directly on permutons.  These models suggest a new (more) algebraic way of studying permuton limits: start with operations on finite permutations coming from Coxeter or Schubert-theoretic structures, tropicalize them, and then pass to continuum measures.

{Complementing these results, Tikhonov~\cite{Tikhonov26} studies a $q$-deformation of the random Demazure-product model introduced in~\cite{morales2025grothendieck}. For the standard reduced decomposition of the longest permutation, he proves convergence to a deterministic permuton and shows that the limit coincides with the $q=0$ limit after a certain parameter substitution, resolving a conjecture of Morales, Panova, Petrov, and Yeliussizov.}

\subsection{Conjectural and emerging limits}

We close with a few directions where the picture is still developing.  The most intriguing conjectural random limit is the meandric permuton.  A meandric permutation records the relative order in which a simple loop and a line visit their intersection points, see~\Cref{fig-meander}~and~\Cref{sect:univ-perm-skew} for further explanations.  Borga, Gwynne, and Sun~\cite{bgs-meander} constructed a random permuton from two independent space-filling SLE$_8$ curves on a suitable Liouville quantum gravity surface and conjectured that it is the scaling limit of uniform meandric permutations.  This conjecture is discussed in the later \Cref{sect:scaling-limit-results} on skew Brownian and SLE/LQG permutons; at the level of this section, its significance is that it would place a century-old enumerative problem about meanders into the same continuum framework as Baxter and separable permutations. 

Other examples are less settled but are worth keeping in mind. Returning to pattern-avoiding permutations, there is now an efficient algorithmic framework, called \emph{Combinatorial Exploration}, to find constructive descriptions (i.e.~combinatorial specifications) of permutation classes, leading to many enumerative results and simulations of uniform random permutations \cite[Section 6]{combinatorial-exploration}. A number of them seem to have nontrivial permuton limits (see the \enquote{heatmaps} of the permutation classes defined by avoiding 2 patterns of length 4 in \cite[Section 6]{combinatorial-exploration}), and it would be very desirable to have some algorithmic way to derive such limits. This is a rather challenging objective since the specifications generated by the {\em Combinatorial Exploration} framework frequently consist of hundreds of equations.

\section{Permuton samples and classical statistics}\label{sec:samples}

As seen in~\Cref{thm:cv_permutons_densities}, the pattern density functionals $d_{\tau}(\cdot)$ are continuous functionals on the space of permutons (and actually characterize permuton convergence).
This is very useful. If one proves the convergence of a sequence of permutations
$\sigma_n$ to some permuton $\mu$,
we get as an immediate corollary the convergence of $d_{\tau}(\sigma_n)$
to $d_{\tau}(\mu)$. 
In particular, the case $\tau=21$ corresponds to the (rescaled) number of inversions,
which is one of the most studied statistics in many random permutation models.
We refer to \cite[Section 1.2.2]{dubach2024mallows} and \cite[Proposition 1.6]{morales2025grothendieck},
where the asymptotic behavior of the number of inversions is deduced from permuton
convergence for Mallows and Grothendieck random permutations, respectively.

Yet, the number of inversions is just one of the many statistics that are typically studied for random permutations. Other standard statistics, such as cycle counts or the length of the longest increasing subsequence (LIS), are not continuous in the permuton topology.
Nevertheless, as we shall see in this section, 
permuton convergence can be a helpful tool to study these statistics as well.

Say we are interested in the asymptotic behavior of $f(\sigma_n)$
where  $\sigma_n$ is a sequence of random permutations of interest,
and $f$ is our favorite statistic (length of LIS, number of fixed points, \dots).
The idea is to compare it 
to $f(\sigma_n^\mu)$, where $\mu$ is the limiting permuton of $\sigma_n$.
A reasonable guess, though this can be false in general,
is that $f(\sigma_n)$ and $f(\sigma_n^\mu)$ have the same asymptotic behavior.
In some cases, this can be proven via some comparison argument,
either for some specific model (see the discussion on LIS of random permutations in Section~\ref{sect:lis-permutons} below)
or by finding regularity conditions on $\sigma_n$ so that the equivalence holds 
(as done for short cycle counts in~\cite{mukherjee2016fixed}, see \Cref{sec:fixed-points-cycles}).

All this motivates the study of three interconnected questions: Given a permutation statistic $f$, 
\begin{enumerate}[(i)]
    \item is (a rescaled version of) $f$ continuous, or upper/lower semicontinuous, with respect to the permuton topology or a stronger topology?
    \item how does $f$ behave on permuton samples?
    \item are there conditions on a permutation sequence $\sigma_n$ converging to a permuton $\mu$ that ensure $f(\sigma_n)$ and $f(\sigma^\mu_n)$ behave similarly?
\end{enumerate}
In Sections~\ref{sect:lis-permutons} and~\ref{sec:fixed-points-cycles}, we discuss the literature on such questions for the LIS and short cycle counts, respectively. Section~\ref{ss:other-stats}
briefly discusses results on other statistics.
Also Table~\ref{table:continuity} recapitulates such continuity properties.

\begin{table}
\[{\begin{tabular}{|p{4cm}|p{10.5cm}|} \hline
\hspace{8mm} {\bf Continuity} & \hspace{4cm} {\bf Statistics} \\ \hline
Continuous & Pattern densities (Theorem \ref{thm:cv_permutons_densities}), proportion of points in a rectangle, simply indexed permutation statistics (Example~\ref{exemp:simply-indexed})\\ \hline
Semi-continuous & LIS at linear order (Equation~\eqref{eq:LIS_semicontinuous})\\ \hline
Continuous under stronger hypotheses
& Short cycle counts (Section~\ref{sec:fixed-points-cycles}), degree process of inversion graphs~\cite{bhat17-deg-inv-graphs}, composition of permutations~\cite{gorski2026compositionpermutons}\\ 
\hline 
No known continuity properties & Total number of cycles, length of the longest cycles, order of the permutation, number of records, \dots \\ \hline
\end{tabular}}\]
\caption{Continuity properties of permutation statistics}
\label{table:continuity}
\end{table}

\subsection{Longest increasing subsequences and Robinson--Schensted shapes}
\label{sect:lis-permutons}

Recall that for a permutation
$\sigma\in\Sym_n$, an \textbf{increasing subsequence} is a sequence of
indices
\[
    1\le i_1<i_2<\dots<i_k\le n
    \qquad\text{such that}\qquad
    \sigma(i_1)<\sigma(i_2)<\dots<\sigma(i_k).
\]
We denote by $\LIS{\sigma}$ the maximal possible value of $k$.  The study
of $\LIS{\sigma}$ for a uniform random permutation is Ulam's problem~\cite{ulam1961monte}.  The first-order asymptotics were obtained by
Logan--Shepp~\cite{Logan-lis}~and~Vershik--Kerov~\cite{Vervsik-lis}, who proved that
\[
    \frac{\LIS{\sigma_n}}{\sqrt n}\xrightarrow[n\to\infty]{\prob}2
\]
for $\sigma_n$ uniform in $\Sym_n$ .
The fluctuations were later identified by Baik, Deift, and Johansson~\cite{Baik-lis}:
after centering by $2\sqrt n$ and scaling by $n^{1/6}$, the limiting law is
the Tracy--Widom GUE distribution.  This problem
is also one of the classical entry points to the Robinson--Schensted correspondence, since $\LIS{\sigma}$ is the length of the first row of the
Young diagrams associated with $\sigma$. See~Romik's book~\cite{Romik-lis} for a beautiful account.

There is now a large literature on longest increasing subsequences in
non-uniform random permutations.  For instance, Mallows permutations have
been studied in several parameter regimes~\cite{ms-mallows-lis,Bhatnagar-lis-mallows,Basu-lis-mallows} (see also~\cite{zhong2023length} for Spearman's
footrule and squared-$L^2$ Mallows models); random
pattern-avoiding permutations give rise to many different behaviors~\cite{Deutsch-lis,Madras-lis,Mansour-lis,bassino2022linear}; and
locally uniform random permutations, obtained from i.i.d.\ points in the plane
with non-constant density, connect the problem to geometric variational
principles~\cite{deuschel1995limiting_curves,deuschel1999increasing_subseq,sjostrand2023monotone}.
Permutons provide a natural language in which to organize these latter examples.

If one wants to give a meaning to the LIS of a permuton, there are two natural approaches.
The first is probabilistic: given a permuton $\mu$, one studies the sequence
$\LIS{\sigma_n^\mu}$, where $\sigma_n^\mu$ is a permutation of size $n$ sampled from
$\mu$.  The second is more geometric: one asks how much mass of $\mu$ can be
captured by an increasing subset of the square.  More precisely, call a set
$A\subseteq[0,1]^2$ \textbf{increasing} if any two points
$(x_1,y_1),(x_2,y_2)\in A$ are comparable for the product order, i.e.\
\begin{equation}\label{eq:defn-inc-set-perm}
    x_1<x_2 \quad\Longrightarrow\quad y_1\le y_2
\end{equation}
up to exchanging the two points.  Define
\[
    \widetilde{\mathrm{LIS}}(\mu)
    :=
    \sup\big\{\mu(A): A\subseteq[0,1]^2
          \text{ is closed and increasing}\big\}.
\]
(One can similarly define $\widetilde{\mathrm{LDS}}(\mu)$ by using
decreasing sets.)
This functional is a natural continuum analogue of the LIS.  More generally,
let $\mathrm{LIS}_k(\sigma)$ be the maximal size of a union of $k$
increasing subsequences of $\sigma$, and define
\[
    \widetilde{\mathrm{LIS}}_k(\mu)
    :=
    \sup\big\{\mu(A): A\subseteq[0,1]^2
          \text{ is a closed union of }k\text{ increasing sets}\big\}.
\]
Define $\mathrm{LDS}_k$ and $\widetilde{\mathrm{LDS}}_k$ analogously. The following theorem of Dubach connects the two approaches.

\begin{thm}[{\cite{dubach2023LIS-linear}}]
\label{thm:dubach-linear-rs}
Let $\mu$ be a permuton, and let $\sigma_n^\mu$ be its permuton samples.
Then, for every fixed $k\ge1$,
\[
    \frac{1}{n}\mathrm{LIS}_k(\sigma_n^\mu)
    \xrightarrow[n\to\infty]{\mathrm{a.s.}}
    \widetilde{\mathrm{LIS}}_k(\mu),
    \qquad
    \frac{1}{n}\mathrm{LDS}_k(\sigma_n^\mu)
    \xrightarrow[n\to\infty]{\mathrm{a.s.}}
    \widetilde{\mathrm{LDS}}_k(\mu).
\]
In particular,
\[
    \frac{\LIS{\sigma_n^\mu}}{n}
    \xrightarrow[n\to\infty]{\mathrm{a.s.}}
    \widetilde{\mathrm{LIS}}(\mu).
\]\end{thm}
Recall that, by Greene’s theorem, the differences 
\[
    \mathrm{LIS}_k(\sigma_n^\mu)
    -
    \mathrm{LIS}_{k{-}1}(\sigma_n^\mu),
    \qquad k\ge1,
\]
and their decreasing analogs are the row and column lengths of
the Robinson--Schensted (RS) shape  associated with $\sigma_n^\mu$.
This can be taken as a definition of the RS shape of a permutation,
and we will not recall here the usual definition of the RS shape through a combinatorial algorithm.
The above theorem gives the limiting RS shape
of permuton samples at a linear scale.
It defines, in some sense, a (linear) RS shape associated with a permuton $\mu$.
Having an equivalent definition via a continuous analog of the combinatorial
algorithm would be of great interest.

Coming back to the LIS,
\Cref{thm:dubach-linear-rs} tells us that when $\widetilde{\mathrm{LIS}}(\mu)>0$, the samples $\sigma_n^\mu$ from
$\mu$ have increasing subsequences of order $n$.  When
$\widetilde{\mathrm{LIS}}(\mu)=0$, the theorem only gives
$\LIS{\sigma_n^\mu}=o(n)$, and the interesting question becomes understanding
 the sublinear scale---the same is true more generally for the RS shape.  This distinction is important: the
permuton topology alone is too coarse to determine sublinear LIS behavior. 

Let us make this last point explicit.  If two permutation diagrams differ
only in $m_n=o(n)$ points, then their associated permutons have
distance tending to zero since the two empirical measures differ by a total
mass of at most $O(m_n/n)$.  However, the LIS can change by as much as $m_n$.  For
example, let $m_n=\lceil n^{\alpha}\rceil$ with $\alpha\in(0,1)$, and let
\[
    \tau_n = 12\cdots m_n \oplus \pi_{n-m_n},
\]
where $\pi_{n-m_n}$ is a uniform random permutation of size $n-m_n$.  Then
$\mu_{\tau_n}$ still converges in probability to the uniform permuton
$\Leb([0,1]^2)$, because the first block has vanishing mass and the second
block occupies all but a vanishing strip of the square.  Nevertheless,
\[
    \LIS{\tau_n}\ge m_n\geq n^\alpha,
\]
whereas the LIS of a uniform permutation is of order $\sqrt n$. Thus, two sequences with the same permuton limit can have very different LIS behavior below the linear scale.  What is robust under permuton convergence is the linear-scale upper bound
\begin{equation}\label{eq:LIS_semicontinuous}
    \limsup_{n\to\infty}\frac{\LIS{\sigma_n}}{n}
    \le
    \widetilde{\mathrm{LIS}}(\mu)
\end{equation}
whenever $\mu_{\sigma_n}\to\mu$ (note that here $\sigma_n\in\Sym_n$ is a general sequence of permutations, not the samples $\sigma_n^\mu$); this follows from the upper semicontinuity of $\widetilde{\mathrm{LIS}}$; see~\cite[Proposition 2.1]{dubach2023LIS-linear}. This is, however, in general a rather weak result. Finding a topology that makes LIS continuous or regularity conditions allowing one to infer information on LIS from permuton convergence seems to be a nontrivial task; see, for instance, Question~\eqref{quest:lis-exp}~in~\Cref{sect:open-prob}.
\medskip

Nevertheless, as explained above, the comparison with permuton samples can be
a good starting point to guess and prove the asymptotic behavior of $\LIS{\sigma_n}$.
The first question is then to find the asymptotic behavior of $\LIS{\sigma^\mu_n}$ 
for large families of permutons $\mu$.

For absolutely continuous permutons, the functional $\widetilde{\mathrm{LIS}}(\cdot)$ and
the RS-shape defined above are trivial.   The correct scale is then typically $\sqrt n$.  This is the setting of Deuschel--Zeitouni for the LIS statistics~\cite{deuschel1995limiting_curves} and, more recently, Sj\"ostrand for the RS shape~\cite{sjostrand2023monotone}.  Given a probability density $\rho$ on $[0,1]^2$, one samples $n$ i.i.d.\ points with density $\rho$ and forms the associated permutation. These permutations are sometimes called \textbf{locally uniform}.  After applying the two marginal distribution functions, as in \eqref{eq:meas-to-perm}, this is precisely a permuton-sample model, which we will therefore denote as $\sigma_n^{\rho}$.   For sufficiently regular $\rho$, that is $\rho>0$ everywhere and $\rho\in C^1([0,1]^2)$, the LIS has order $\sqrt n$, and the leading constant is described by a variational problem over increasing curves.
Namely, we have
\begin{equation}\label{eq:DZ-asymptotics}
   \frac{\LIS{\sigma_n^{\rho}}}{\sqrt n}\xrightarrow[]{\mathbb P} 
    2 \sup_\gamma\left(\int_0^1
    \sqrt{{\rho}(\gamma(t))\,\gamma_x'(t)\,\gamma_y'(t)}\,dt \right),
\end{equation}
where the supremum is taken over increasing curves $\gamma=(\gamma_x,\gamma_y)$ in the unit square. Sj\"ostrand~\cite{sjostrand2023monotone} extends this picture from one curve to surfaces associated with unions of order-$\sqrt n$ monotone subsequences, obtaining limits for the RS shape of locally uniform random permutations.  A later paper~\cite{sjostrand2024continuity} proves continuity and existence properties for the associated variational maximizers under boundedness assumptions on the density and its support.

Let us return to the analysis of the LIS of permuton samples.
Between the regular $\sqrt n$ regime and the linear regime, there are intermediate behaviors. In the paper~\cite{dubach2023LIS-inermediate}, Dubach showed that, for locally uniform permutations whose sampling densities have suitable singularities, the LIS can grow like essentially any prescribed power $n^\alpha$ with $\alpha\in(1/2,1)$, up to logarithmic factors.  For example, densities that diverge near an increasing curve force many sampled points to accumulate close to that curve, producing longer increasing subsequences than in the smooth case, but still not a linear proportion of all points.

As suggested above, these results on permuton samples
can be helpful for analyzing the LIS of other permutation models.
In this spirit,
Deuschel--Zeitouni's result has been utilized by Mueller and Starr \cite{ms-mallows-lis} to find the limiting behavior
of $\LIS{\sigma_n}$, where $\sigma_n$ is a Mallows random permutation in the regime $q=\exp(2\beta/n)$---this is the regime where $\sigma_n$ has a nontrivial permuton limit $\mu_\beta$; recall Example~\ref{ex:Mallows}. It turns out that $\LIS{\sigma_n}$ has the same asymptotic behavior as $\LIS{\sigma^{\mu_\beta}_n}$, which can be found from Deuschel--Zeitouni's result \eqref{eq:DZ-asymptotics}. Remarkably, in this case, the variational problem is solved by the diagonal curve $\gamma(t)=(t,t)$, and the corresponding value can be computed explicitly; see~\cite{ms-mallows-lis}.
The same strategy---first finding a limiting permuton and then studying the LIS by comparison---has later been used by Jin on products of independent Mallows permutations;
see~\cite{jin19-perm-lim-prod-Mallows,Ke-lis-mallows}.
Extending such a strategy---for a single Mallows permutation or for a product of two independent ones---to the full RS shape using the work of Sj\"ostrand seems quite challenging.

{The Mallows model is not the only model for which this strategy can be applied. In recent work, Borga and Dubach~\cite{BD26-LIS-Luce} study the LIS of the so-called \emph{Luce model} of random permutations. In this model, weights $\theta_1,\dots,\theta_n$ are given, and the entries of the permutation are chosen successively without replacement, with a remaining entry $i$ selected with probability proportional to its weight $\theta_i$. As explained in Section~\ref{ss:analytic-limit-shapes} (see also the right-hand side of~Figure~\ref{fig:mallows_luce}), the permuton limit of this model was recently identified in \cite{borga2025luce}. 

Borga and Dubach show that the LIS of Luce permutations can exhibit a wide range of asymptotic behaviors. For arbitrary nonincreasing weights, they give an estimate of the expected LIS in terms of the distribution of the weights across logarithmic scales; in particular, suitable choices of weights yield LIS of order $n^\beta$ for any $\beta\in[1/2,1)$. In the regular $\sqrt n$ regime, when the weights are obtained by sampling a sufficiently regular profile $f$, they further prove that the LIS has the same first-order asymptotics as that of a sample from the limiting permuton $\mu_f$. The limiting constant is therefore given by the Deuschel--Zeitouni variational formula applied to the density of $\mu_f$ (but, in this setting, solving the variational problem seems a difficult task).}

In conclusion, the picture that emerges is as follows.  The permuton topology partly controls the linear part of the longest increasing subsequence through the upper-semicontinuous functional $\widetilde{\mathrm{LIS}}$.  Absolutely continuous permutons have no linear increasing subsequence, and the classical scale is $\sqrt n$, governed by variational curves and surfaces. Singular permutons can produce linear behavior when an increasing set carries positive mass, or intermediate polynomial behavior when their mass is arranged along a random fractal geometry.  The skew Brownian permutons provide the first examples where these intermediate exponents are now known exactly; we will come back to this in ~\Cref{sect:univ-perm-skew}.

\subsection{Fixed points and cycle structure}
\label{sec:fixed-points-cycles}

We now turn to other classical permutation statistics: fixed points and cycles. For a permutation $\sigma\in\Sym_n$, let
\[
    F(\sigma):=\big|\{i\in[n]:\sigma(i)=i\}\big|
\]
be the number of fixed points, and let $C_\ell(\sigma)$ be the number of cycles of length $\ell$ in the cycle decomposition of $\sigma$.  For a uniform random permutation, the classical Goncharov--Kolchin theorem says that, for each fixed $L$,
\[
    \big(C_1(\sigma_n),\ldots,C_L(\sigma_n)\big)\xrightarrow[n\to\infty]{d}\big(Z_1,\ldots,Z_L\big),
\]
where $Z_1,\ldots,Z_L$ are independent and $Z_\ell\sim\Poisson(1/\ell)$. In particular, $F(\sigma_n)=C_1(\sigma_n)$ converges to $\Poisson(1)$.  See, for instance, \cite{arratia2003logarithmic} for a broad account of these classical Poisson approximation phenomena.

Unlike pattern densities, fixed-point and cycle counts are not continuous with respect to permuton convergence.  The reason is that they depend on the exact ``lattice-scale'' behavior, while the permuton topology only remembers the macroscopic behavior of the permutation diagram.  For example, let $e_n$ be the identity permutation and let $s_n$ be the cyclic shift $s_n(i)=i+1$ for $i<n$ and $s_n(n)=1$.  Then both $\mu_{e_n}$ and $\mu_{s_n}$ converge to the diagonal permuton, but
\[
    F(e_n)=n,\qquad F(s_n)=0.
\]
Thus even the normalized statistic $F(\sigma)/n$ is not continuous at the diagonal permuton.

A theorem of Mukherjee nevertheless gives a positive result when permuton convergence is supplemented by a local regularity assumption \cite{mukherjee2016fixed}.  Roughly speaking, the idea is that if a random permutation has a smooth limiting density $\rho$, and if the ``assignment probabilities'' are locally uniform at scale $1/n$, then fixed points and finite cycles have universal Poisson limits whose parameters are computed from $\rho$. We now explain this rigorously.

For $r\ge1$, set
\[
    S(n,r):=\big\{\bm p=(p_1,\ldots,p_r)\in[n]^r:p_a\ne p_b\text{ for }a\ne b\big\}.
\]
If $\pi_n$ is a random permutation, write
\[
    \pi_n(\bm p):=(\pi_n(p_1),\ldots,\pi_n(p_r)).
\]
Let $\rho$ be a strictly positive continuous density on $[0,1]^2$ with uniform marginals.  We say that $\pi_n$ satisfies the \textbf{local density condition} with density $\rho$ if, for every fixed $r\ge1$,
\begin{equation}
\label{eq:local-density-condition}
    \lim_{n\to\infty}\sup_{\bm p,\bm q\in S(n,r)}\left|\frac{n^r\,\prob\big(\pi_n(\bm p)=\bm q\big)}{\prod_{a=1}^r\rho\left(\frac{p_a}{n},\frac{q_a}{n}\right)}-1\right|=0.
\end{equation}
In words, the images of any fixed number of distinct positions behave, at the lattice scale, as if they were sampled independently from the limiting density $\rho$.  This is much stronger than the weak convergence of $\mu_{\pi_n}$ to the permuton with density $\rho$.\footnote{Mukherjee proves in~\cite[Theorem 1.1]{mukherjee2016fixed} that \eqref{eq:local-density-condition} follows from convergence to this permuton together with the joint equicontinuity condition in~\cite[Definition 1.2]{mukherjee2016fixed}. His Theorems 1.3 and 1.4 are stated under the resulting local-density condition (1.5), corresponding to~\eqref{eq:local-density-condition} here.}

We first record the fixed-point theorem in a slightly more general form. If $\pi,\sigma\in\Sym_n$, define their \textbf{overlap}
\[
    N_n(\pi,\sigma)
    :=
    \sum_{i=1}^n \One_{\{\pi(i)=\sigma(i)\}}.
\]
Thus $F(\pi)=N_n(\pi,e_n)$, where $e_n$ is the identity permutation.

\begin{thm}[{\cite[Theorem 1.3]{mukherjee2016fixed}}]
\label{thm:mukherjee-overlaps}
Let $\pi_n$ be random permutations satisfying \eqref{eq:local-density-condition} for a strictly positive continuous density $\rho$ with uniform marginals.  Let $\sigma_n\in\Sym_n$ be a deterministic sequence such that $\mu_{\sigma_n}$ converges to a permuton $\nu$.  Then
\[
    N_n(\pi_n,\sigma_n)\xrightarrow[n\to\infty]{d}\Poisson\big(\nu[\rho]\big),\qquad\nu[\rho]:=\int_{[0,1]^2}\rho(x,y)\,\nu(dx,dy).
\]
The convergence also holds in moments.
\end{thm}

Taking $\sigma_n=e_n$, whose permuton limit is the diagonal permuton, gives the limiting number of fixed points:
\[
    F(\pi_n)\xrightarrow[n\to\infty]{d}\Poisson\left(\int_0^1 \rho(x,x)\,dx\right).
\]
This formula has a simple interpretation.  Although the limiting permuton with density $\rho$ gives zero mass to the diagonal, the density of the permuton \emph{along} the diagonal determines the expected number of fixed points at the microscopic scale.

The cycle theorem has a similarly explicit form.  For $\ell\ge1$, define
\begin{equation}
\label{eq:c-rho-def}
    c_\rho(\ell):=\frac1\ell\int_{[0,1]^\ell}\rho(x_1,x_2)\rho(x_2,x_3)\cdots\rho(x_{\ell-1},x_\ell)\rho(x_\ell,x_1)\,dx_1\cdots dx_\ell.
\end{equation}
Equivalently, if $T_\rho$ is the integral operator on $L^2[0,1]$ with kernel $\rho$, then $c_\rho(\ell)=\ell^{-1}\operatorname{Tr}(T_\rho^\ell)$.\footnote{Strictly speaking, for $\ell=1$ the operator $T_\rho$ need not be trace class under our assumptions, so $\operatorname{Tr}(T_\rho)$ is not necessarily defined. In this case, $c_\rho(1)=\int_0^1 \rho(x,x)\,dx$ by definition.}

\begin{thm}[Mukherjee {\cite[Theorem 1.4]{mukherjee2016fixed}}]
\label{thm:mukherjee-cycles}
Let $\pi_n$ be random permutations satisfying \eqref{eq:local-density-condition} for a strictly positive continuous density $\rho$ with uniform marginals.  Then, for every fixed $L\ge1$,
\[
    \big(C_1(\pi_n),\ldots,C_L(\pi_n)\big)\xrightarrow[n\to\infty]{d}(Z_1,\ldots,Z_L),
\]
where $Z_1,\ldots,Z_L$ are independent and
\[
    Z_\ell\sim \Poisson(c_\rho(\ell)).
\]
The convergence also holds in all mixed moments.
\end{thm}

For the uniform density $\rho\equiv1$, we have $c_\rho(\ell)=1/\ell$, recovering the classical Poisson limit for finite cycles of a uniform random permutation.  For a non-uniform model, the same Poisson structure survives, but the means are deformed by the limiting density $\rho$.  The proof follows the usual rare-event philosophy.  A fixed point is the event $\pi_n(i)=i$, which has probability of order $1/n$; an $\ell$-cycle is the event that $\pi_n$ sends a fixed $\ell$-tuple cyclically to itself, which has probability of order $n^{-\ell}$.  The local density condition gives the correct asymptotics for these probabilities, and Stein's method or the method of moments controls the weak dependence arising from the fact that $\pi_n$ is a permutation.

\medskip

\noindent\textbf{Examples.}
Mukherjee verifies the hypotheses above for several important classes of non-uniform random permutations.
First, consider Mallows permutations with parameter $q_n$ in the regime
\[
    \frac{n}{2}\log q_n\longrightarrow \beta\in\mathbb R.
\]
In this regime, $\pi_n$ converges to Starr's Mallows permuton with density $\rho_\beta$; see \Cref{ex:Mallows}.  Mukherjee proves that the local density condition holds, and therefore
\[
    F(\pi_n)\xrightarrow[n\to\infty]{d}\Poisson\left(\int_0^1\rho_\beta(x,x)\,dx\right),
\]
while, for every fixed $L$,
\[
    \big(C_1(\pi_n),\ldots,C_L(\pi_n)\big) \xrightarrow[n\to\infty]{d} \big(Z_1,\ldots,Z_L\big), \qquad Z_\ell\sim \Poisson(c_{\rho_\beta}(\ell)),
\]
with the $Z_\ell$'s independent.  

Second, Mukherjee considers one-parameter exponential families of permutations of the form
\[
    \prob(\pi_n=\pi)
    \propto
    \exp\left\{
        n\theta\,\mu_\pi[f]
    \right\},
\]
where $f$ is a continuous function on $[0,1]^2$.  This family includes scaled versions of classical rank statistics such as Spearman's footrule and Spearman's rank correlation.  The associated random permutations converge to a deterministic permuton with a continuous density, say $g_{f,\theta}$, and the same Poisson formulas hold with $\rho$ replaced by $g_{f,\theta}$.

Third, if $\mu$ itself is a permuton with strictly positive continuous density $\rho$, and if $\pi_n=\sigma_n^\mu$ is sampled from $\mu$, then the local density condition holds.  Hence permutations sampled from a smooth permuton have asymptotically independent Poisson finite-cycle counts with parameters given by \eqref{eq:c-rho-def}.  This is a useful complement to the basic convergence $\mu_{\sigma_n^\mu}\to\mu$: the weak permuton limit identifies the macroscopic shape, while the density along cyclic configurations identifies the microscopic short cycle statistics.
\medskip

The regularity hypothesis is essential.  Mukherjee gives a particularly simple warning example.  Let $R_{n,\theta}$ be the probability measure on $\Sym_n$ defined by
\[
    R_{n,\theta}(\pi)\propto\exp\big\{\theta F(\pi)\big\}.
\]
For every fixed $\theta \neq 0$, a random permutation with law $R_{n,\theta}$ still converges to the uniform permuton.  Nevertheless
\[
    F(\pi_n)\xrightarrow[n\to\infty]{d}\Poisson(e^\theta),
\]
not $\Poisson(1)$.  Thus, two random permutation models may have the same permuton limit but different limiting numbers of fixed points.  The missing ingredient is precisely the local equicontinuity: the biasing factor $\exp\{\theta F(\pi)\}$ changes probabilities at the lattice scale along the diagonal without changing the macroscopic empirical measure.

{Long cycles in Spearman's
footrule model have also been analyzed. Let $\pi_n$ have distribution proportional to $\exp(-\beta_n\sum_{i=1}^n|\pi_n(i)-i|)$, where $\beta_n>0$. Zhong~\cite[Theorems~1.1 and~1.3]{zhong2023cycles} shows that the expected length of the cycle containing a given index is of order $\min\{\max\{\beta_n^{-2},1\},n\}$, uniformly in the chosen index. Moreover, if $\sqrt n\,\beta_n\to0$, the decreasing sequence of cycle lengths, divided by $n$, converges in distribution to the Poisson--Dirichlet law with parameter $1$. In particular, this latter conclusion applies when $\beta_n=\beta/n$ for any fixed $\beta>0$}.

\subsection{Other statistics}\label{ss:other-stats}
Other statistics of permuton samples have been considered in recent years. For example, the paper \cite{corsini25-BST-samples} analyzes the height of the binary search tree (BST) associated with the permuton sample $\sigma_n^\mu$ and gives sufficient conditions on $\mu$ so that it behaves essentially as that of  a uniform random permutation.
In another direction, the recent paper \cite{maga2025entropy} analyzes the Shannon entropy\footnote{Note that this entropy of permuton samples is a priori unrelated to the entropy $H$ of the permuton, as defined in the context of large deviation (see Section~\ref{sec:LDP}).} of permuton samples. Various behaviors (including a surprising log-periodicity phenomenon) appear depending on the structure of  $\mu$ (permuton with density, or supported by the graph $\{(x,f(x)), x\in [0,1]\}$ of a function, \dots). Let us also mention two papers studying degrees in the inversion graph of permuton samples: \cite{bhat17-deg-inv-graphs} is interested in the global degree distribution, while \cite{corsini26-min-inv-perm-sample} looks at the minimal degree. Interestingly, \cite{bhat17-deg-inv-graphs} also provides a general condition (the convergence of the ``permutation process'') which implies the convergence of the degree distribution. Finding such ``continuity criteria'' for the other statistics mentioned above---height of the BST, entropy, minimal degree of the inversion graphs---would be interesting.
{Let us finally mention the recent preprint \cite{gorski2026compositionpermutons},
which introduces a notion of compositions of permutons (extending that of permutations),
and studies its continuity properties.
In particular, the author shows that, in certain cases, one can deduce convergence of a composition of permutations from the convergence of the factors, and applies it to
products of uniform random separable permutations, see  \cite[Theorem 1.3]{gorski2026compositionpermutons}.}

\section{Universal random permutons: skew Brownian, Brownian separable, and meandric}
\label{sect:univ-perm-skew}

The skew Brownian permutons mentioned in \Cref{sect:pattern-av-perm} form a two-parameter family of fractal random permutons that have been shown, or are conjectured, to describe the permuton limits of many families of random constrained permutations, especially pattern-avoiding and vincular-pattern-avoiding permutations. They are part of a larger family of permutons constructed from two continuum curves in the plane, called \emph{Schramm--Loewner evolution} curves, and a random area measure called \emph{Liouville quantum gravity}. We start with three examples of discrete families of permutations that will help build intuition before we define the skew Brownian permutons and the more general SLE/LQG permutons.

\begin{figure}[t!]
	\begin{center}
		\includegraphics[width=0.8\textwidth]{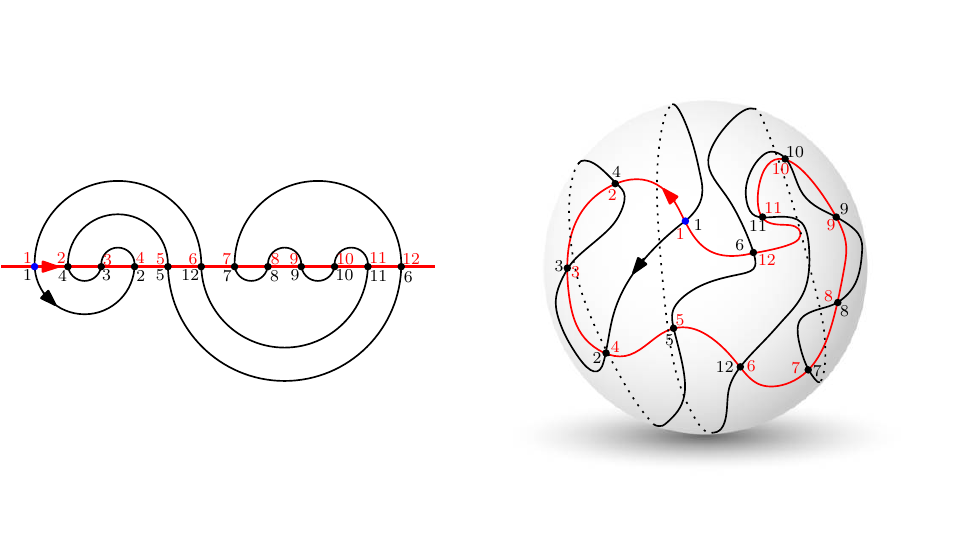}  
		\caption{\label{fig-meander}\textbf{Left:} A meander of size $6$, with the intersection points of the loop and the line labeled by the order in which they are hit by the loop (in black) and the line (in red) starting at the blue vertex. The associated meandric permutation is $\sigma_{\ell}=1,4,3,2,5,12,7,8,9,10,11,6$ in one-line notation; see \Cref{def-meandric-perm}. \textbf{Right:} The same meander of size $6$ viewed as a pair of Jordan curves in the sphere.
		}
	\end{center}
\end{figure}

\medskip

\noindent\textbf{Example 1: Meandric permutations.}
Versions of the following question can be traced back at least to the work of Henri Poincar\'e~\cite{poincare-meander}:
\begin{equation}\label{eq:quest}
	\text{\emph{In how many ways can a simple loop in the plane cross a line a specified number of times?}}
\end{equation}
Despite many efforts, this question remains open after more than a century. Let us start by formalizing it.

\begin{defn} \label{def-meander}
	A \textbf{meander} of size $n\in\mathbb N$ is a simple loop $\ell$, i.e.\ a homeomorphic image of the circle, in $\mathbb R^2$ which crosses the real line $\mathbb R$ exactly $2n$ times and does not touch the line without crossing it. Meanders are viewed modulo orientation-preserving homeomorphisms from $\mathbb R^2$ to $\mathbb R^2$ which take $\mathbb R$ to $\mathbb R$ and whose restrictions to $\mathbb R$ are increasing.
\end{defn} 

\begin{figure}[t!]
	\begin{center}
		\includegraphics[width=.25\textwidth]{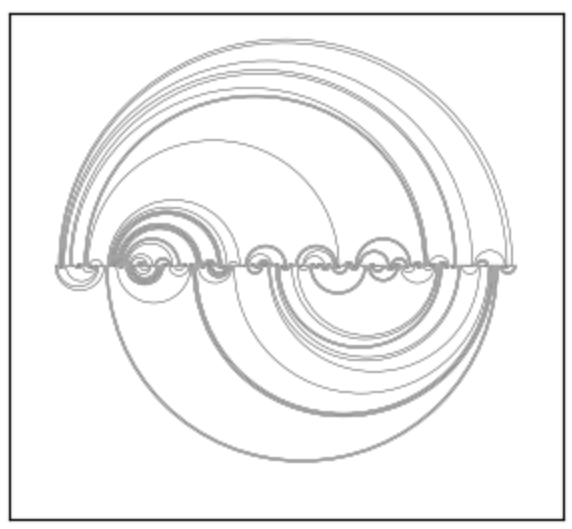}
		\includegraphics[width=.25\textwidth]{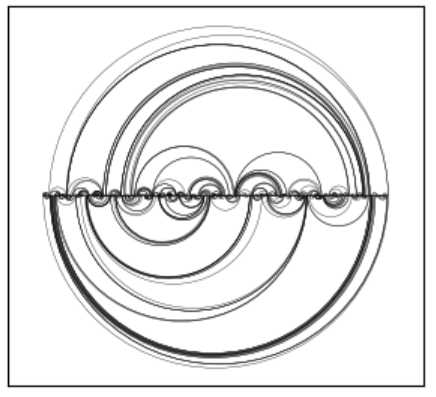}
		\includegraphics[width=.23\textwidth]{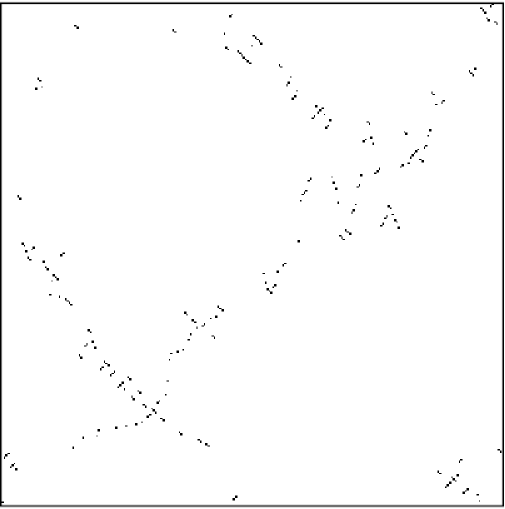}
		\includegraphics[width=.23\textwidth]{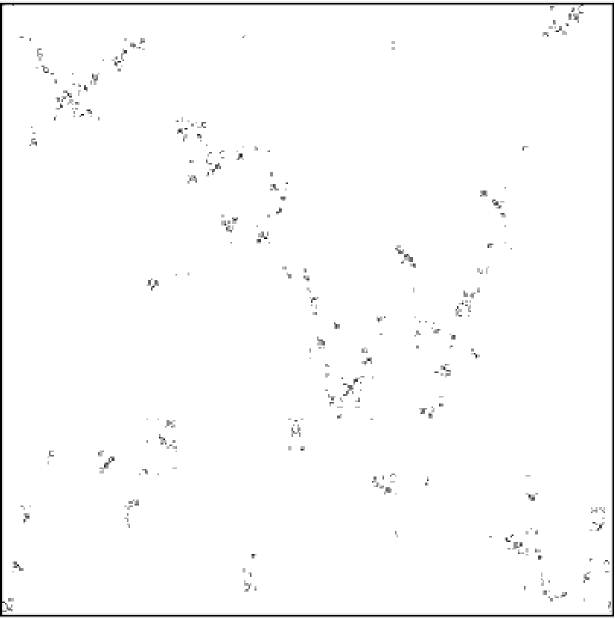} 
		\caption{\label{sim-menders-permuton}\textbf{Left:} Two large uniform meanders $\ell$ of size 256 and 2048. \textbf{Right:} The plots of the two corresponding permutons associated with the meandric permutations $\sigma_{\ell}$. These simulations are obtained using the Markov chain Monte Carlo algorithm from \cite{heitsch2011meander}.}
	\end{center}
\end{figure}

See the left-hand side of~\Cref{fig-meander,sim-menders-permuton} for illustrations of meanders. A meander of size $n$ can equivalently be viewed as an ordered pair of oriented Jordan curves on a sphere, intersecting at exactly $2n$ points and only by crossings; see the right-hand side of~\Cref{fig-meander}. 
It turns out that meanders can be encoded by a subset of permutations, as we shall see now.

\begin{defn} \label{def-meandric-perm}
	For a meander $\ell$ of size $n$, the corresponding \textbf{meandric permutation} $\sigma_{\ell}$ is the permutation of $[2n]$ obtained as follows. Order the $2n$ intersection points of $\ell$ with $\mathbb R$ from left to right. Then order the same intersection points according to the order in which they are hit by $\ell$, starting from the leftmost intersection point and traversing the loop counterclockwise. We set $\sigma_\ell(i)=j$ if the intersection point which is the $i$-th point in the first order is the $j$-th point in the second order.
\end{defn} 

See \Cref{fig-meander} for an example and the right-hand side of \Cref{sim-menders-permuton} for the plots of two large uniform meandric permutations. With the convention above, every meandric permutation fixes $1$. This constraint can be removed by allowing the loop and the line to have different distinguished starting points. We denote by $\operatorname{Mea}_1$ the set of meandric permutations fixing $1$, and by $\operatorname{Mea}$ the corresponding set when this constraint is removed. It is straightforward to check that the map $\ell\mapsto\sigma_{\ell}$ is a bijection from meanders of size $n$ to meandric permutations of size $2n$ fixing $1$. The study of meandric permutations dates back at least to the paper~\cite{rosentiehl-meander}, where they are called ``planar permutations''.

The term ``meander'' was coined by Arnold~\cite{arnold-meander}. Meanders are of interest in physics and computational biology as models of polymer folding~\cite{dgg-meander-folding}. There is a vast mathematical literature devoted to the enumeration of various types of meanders, with connections to combinatorics, theoretical physics, and more recently to the geometry of moduli spaces~\cite{Delecroix-enum-meand}. See~\cite{zvonkin-meander-survey} for a brief recent survey and~\cite{lacroix-meander-survey} for a more detailed account. 

The problem of enumerating \emph{all} meanders of size $n$, that is, the question in~\eqref{eq:quest}, is open and appears quite difficult. Di Francesco, Golinelli, and Guitter~\cite{dgg-meander-asymptotics} conjectured that the number of meanders of size $n$ behaves asymptotically like $C \cdot A^n n^{-\alpha}$, where
\[
    \alpha = \frac{29 + \sqrt{145}}{12} \approx 3.42.
\]
This conjecture was later numerically tested by Jensen and Guttmann~\cite{Jensen-num-meanders}, who proposed the numerical estimate $A \approx 12.26$. The best known rigorous bounds for the constant $A$ are $11.380\leq A \leq 12.806$, as proved in \cite{Albert-bounds-meanders,CDDV26}. 

\medskip

\noindent\textbf{Example 2: Baxter permutations.}
A simpler model of meanders is provided by monotone meanders.

\begin{defn} \label{def-monotone-meander}
	A \textbf{monotone meander} of size $n\in\mathbb N$ is a pair of simple curves $\ell_1$ and $\ell_2$ in $[0,1]^2$ which cross each other exactly $2n-1$ times and such that:
	\begin{itemize}
		\item $\ell_1$ starts on the left-hand side of $[0,1]^2$, ends on the right-hand side of $[0,1]^2$, and always moves to the right. Equivalently, it is the graph of a continuous function from $[0,1]$ to $[0,1]$.
		\item $\ell_2$ starts on the bottom side of $[0,1]^2$, ends on the top side of $[0,1]^2$, and always moves upwards. Equivalently, it is the graph of a continuous function from $[0,1]$ to $[0,1]$ rotated by 90 degrees.
	\end{itemize} 
	As in the case of meanders, one can associate to each monotone meander a \textbf{monotone meandric permutation}. Two monotone meanders are equivalent if they induce the same monotone meandric permutation.
\end{defn} 

\begin{figure}[t!]
	\begin{center}
		\includegraphics[width=0.46\textwidth]{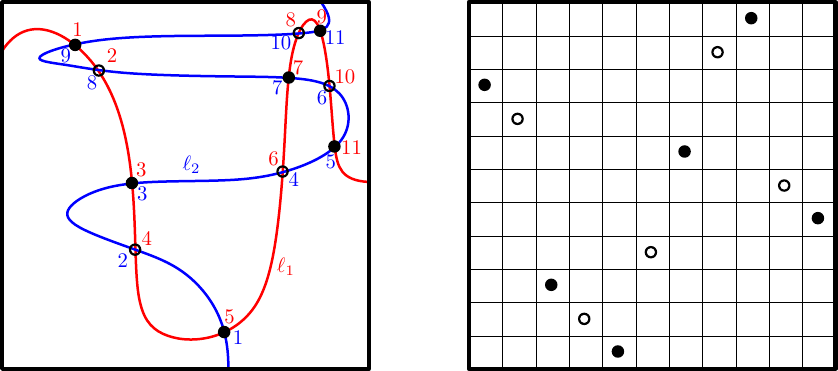} 
		\hspace{0.5 cm} 
		\includegraphics[scale=0.37]{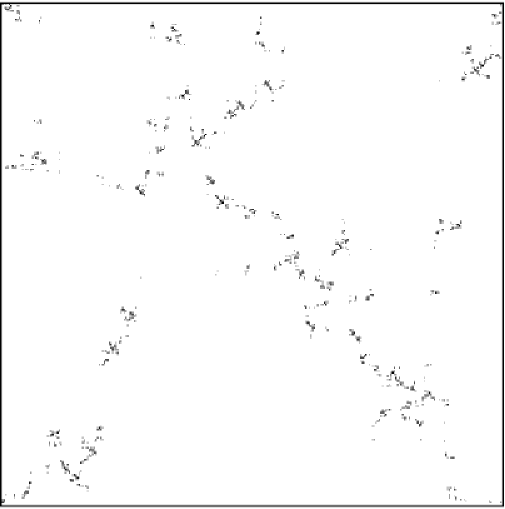}
		\hspace{0.5 cm}
		\includegraphics[scale=0.37]{Baxter_perm_4520}\\
		\caption{\label{fig-monotone-meander}\textbf{Left:} A monotone meander of size $6$, with the intersection points labeled by the order in which they are hit by the curve $\ell_1$ (in red) and the curve $\ell_2$ (in blue). \textbf{Middle-left:} The associated monotone meandric permutation: the Baxter permutation corresponds to the black dots, while the reduced Baxter permutation corresponds to the white dots.
		\textbf{Right:} The diagrams of two uniform Baxter permutations of sizes 3253 and 4520. 	
		}
	\end{center}
\end{figure}

See \Cref{fig-monotone-meander} for an example. A monotone meandric permutation maps even numbers to even numbers and odd numbers to odd numbers. After relabeling, its restriction to the odd numbers is called a \textbf{Baxter permutation}, while its restriction to the even numbers is called a \textbf{reduced Baxter permutation}. It can be shown \cite{baxter-permutation,MR0250516,boyce1981} that the Baxter permutation corresponding to a monotone meandric permutation uniquely determines the monotone meandric permutation itself. Hence, there is a natural bijection between monotone meandric permutations and Baxter permutations. 
{Baxter permutations are more widely known as permutations} avoiding the vincular patterns $2\underline{41}3$ and $3\underline{14}2$, where the underlined entries are required to be adjacent. We denote by $\operatorname{Bax}$ the set of Baxter permutations.

Baxter permutations were introduced by Glen Baxter in 1964 \cite{baxter-permutation} while studying fixed points of commuting functions. They are classical examples of pattern-avoiding permutations and have been intensively studied in the probabilistic and combinatorial literature; see, e.g., \cite{MR0250516,MR491652,MR555815,MR2028288,MR2679559,MR2763051,MR3882946}. They are connected with various other combinatorial structures, such as bipolar orientations~\cite{bbf-baxter-bipolar}, walks in cones~\cite{kmsw-bipolar}, certain pairs of binary trees, triples of non-intersecting lattice paths~\cite{MR2763051}, and domino tilings of Aztec diamonds~\cite{MR2679559}. The diagrams of two large uniform Baxter permutations are plotted on the right-hand side of \Cref{fig-monotone-meander}.

\medskip

\noindent\textbf{Example 3: Separable permutations.}
Separable permutations form a fundamental example of a permutation class. They are defined recursively as the smallest class of permutations containing the singleton permutation $1$ and closed under the direct sum $\oplus$ and skew sum $\ominus$ operations. Equivalently, they are precisely the permutations avoiding the two classical patterns $2413$ and $3142$, that is
\[
\operatorname{Sep}=\operatorname{Av}(2413,3142).
\]
This recursive structure makes separable permutations especially tractable: they admit a natural decomposition by signed trees, are counted by the large Schr\"oder numbers, and often serve as a useful testing ground for questions in pattern avoidance, permutation classes, and combinatorial decompositions~\cite{shapiro1991bootstrap,ghys2017singular}.
\medskip

\begin{figure}[t!]
	\begin{center}
		\includegraphics[width=0.8\textwidth]{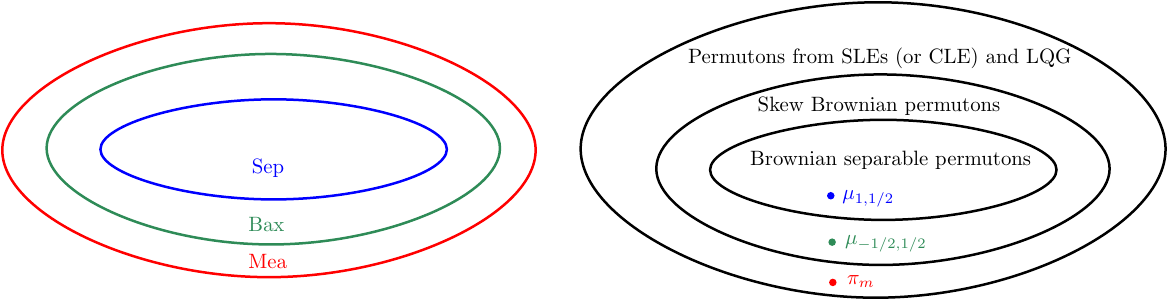} 
		\caption{\label{fig-Venn_permutons}\textbf{Left:} A Venn diagram showing the three families of permutations considered in this section: meandric, Baxter, and separable permutations.
		\textbf{Right:} A Venn diagram showing the corresponding families of limiting random permutons: permutons constructed from SLEs (or CLEs) and LQG, skew Brownian permutons, and Brownian separable permutons. The permuton limit of uniform meandric permutations is conjectured to be the meandric permuton $\bm{\pi}_m$ (\Cref{conj:mean_conv}); the permuton limit of uniform Baxter permutations is the skew Brownian permuton $\bm{\mu}_{-1/2,1/2}$ (\Cref{thm:baxt_conv}); the permuton limit of uniform separable permutations is the Brownian separable permuton $\bm{\mu}_{1,1/2}$ (\Cref{thm:sep_conv}).}
	\end{center}
\end{figure}

Separable permutations are closely related to Baxter permutations. {Comparing the characterization via pattern avoidance, we find that}
\[
\operatorname{Sep}\subseteq\operatorname{Bax}.
\]
Moreover, Baxter permutations are naturally in bijection with monotone
meanders, which can be completed to ordinary meanders by adding one crossing. 
We therefore use the schematic notation
\[
\operatorname{Sep}\subseteq\operatorname{Bax}
\subseteq\operatorname{Mea},
\]
where the second inclusion is understood through {this transformation which} changes the size of the associated permutation. The inclusions are strict; see the left-hand side of \Cref{fig-Venn_permutons}.

\subsection{Constructing permutons in terms of space-filling curves}
\label{sect:space-filling-constr}

We have seen that, in the discrete setting, two intersecting curves naturally give rise to a permutation. The skew Brownian permutons can be defined in a similar spirit directly in the continuum, using a \textbf{Liouville quantum gravity (LQG)} area measure and a pair of coupled \textbf{Schramm--Loewner evolution (SLE)} curves. There is a vast literature concerning SLE and LQG, but no prior knowledge of this literature is needed for the discussion here: we will give only an intuitive description and refer to the literature for precise constructions and definitions.

Fix parameters $\gamma \in (0,2)$ and $\kappa_1,\kappa_2 > 4$. Let $\mu_h$ be the area measure corresponding to a unit area $\gamma$-Liouville quantum sphere. Although its construction is involved, here we only need to know that $\mu_h$ is a random, non-atomic, fractal Borel probability measure on $\mathbb C$ which assigns positive mass to every open subset of $\mathbb C$. See, e.g., \cite[Section 2]{bgs-meander} for a more precise account.

Independently of $\mu_h$, let $(\eta_1,\eta_2)$ be a random pair consisting of a whole-plane space-filling SLE$_{\kappa_1}$ curve and a whole-plane space-filling SLE$_{\kappa_2}$ curve, each going from $\infty$ to $\infty$. For the purposes of this survey, one can think of $(\eta_1,\eta_2)$ as a pair of random space-filling curves in $\mathbb C$, each of which visits almost every point of $\mathbb C$ exactly once. We parametrize each of $\eta_1$ and $\eta_2$ by $\mu_h$-mass, i.e.\ in such a way that
\[
    \mu_h(\eta_1([0,t])) = \mu_h(\eta_2([0,t])) = t, \qquad 0\le t\le 1.
\]
See again~\cite[Section 2]{bgs-meander} for details.

We emphasize that the coupling of the two SLE curves, viewed modulo time parametrization, is arbitrary at this stage. {Even extremal cases where each one of the two curves determines the other or where the two curves are independent can lead to nontrivial objects}.

{Since $\eta_1$ and $\eta_2$ are both space-filling, we claim that there exists a Lebesgue-measurable function $\psi : [0,1] \to [0,1]$ }such that 
\begin{equation} \label{eqn-psi-property}
    \eta_1(t) = \eta_2(\psi(t)),\qquad 0\le t\le 1 .
\end{equation}
The \textbf{permuton associated with $(\mu_h,\eta_1,\eta_2)$} is the random probability measure $\bm{\pi}$ on $[0,1]^2$ defined by
\begin{equation} \label{eqn-permuton-def}
\bm{\pi}(A) = \bm{\pi}_{\mu_h,\eta_1,\eta_2}(A)
= \Leb\left\{ t\in [0,1] : (t,\psi(t)) \in A \right\},
\qquad \text{for every Borel set }A\subset [0,1]^2.
\end{equation}
Although $\psi$ is not necessarily unique, the definition of $\bm{\pi}$ is independent of this choice. Equivalently, for rectangles one has
\[
    \bm{\pi}([a,b]\times[c,d])
    =
    \mu_h\big(\eta_1([a,b])\cap \eta_2([c,d])\big);
\]
see \cite[Lemma 2.8]{bgs-meander}. This formula is often the most intuitive way to remember the construction: the mass of a rectangle in the permuton records the LQG mass of the portion of the surface visited by $\eta_1$ during the time interval $[a,b]$ and by $\eta_2$ during the time interval $[c,d]$.

There are three special cases of the construction in \eqref{eqn-permuton-def} that are particularly important here.

\begin{enumerate}
	
	\item\label{item:mean_perm} The \textbf{meandric permuton} $\bm{\pi}_{m}$ is the permuton $\bm{\pi}$ in~\eqref{eqn-permuton-def} in the case when
    \[
        \gamma = \sqrt{\frac13 \left( 17 - \sqrt{145} \right)} \approx 1.2856,\qquad
        \kappa_1 = \kappa_2 = 8,
    \]
    and $\eta_1$ and $\eta_2$ are independent, viewed as curves modulo time parametrization. This is the conjectural permuton limit of uniform meandric permutations \cite[Conjecture 1.8]{bgs-meander}.
		
	\item\label{item:skew_perm} The \textbf{skew Brownian permutons} $\bm{\mu}_{\rho,q}$, indexed by $(\rho,q) \in (-1,1] \times (0,1)$, are known to be universal limits of several families of constrained permutations~\cite{bassino-separable-permuton,Bassino-universal-perm,bm-baxter-permutation,borga2021permuton}. Their original construction does not involve SLEs or LQG: it is given in terms of a $q$-dependent system of stochastic differential equations driven by a two-dimensional Brownian excursion in the first quadrant with correlation $\rho$; see \cite{borga-skew-permuton}. However, \cite[Theorem 1.17]{borga-skew-permuton} shows that for each $(\rho,q)\in(-1,1)\times(0,1)$ there exists a $q$-dependent coupling of two whole-plane space-filling SLE$_{\kappa}$ curves, with $\kappa=16/\gamma^2$ and
    \[
        \rho = -\cos(\pi \gamma^2/4),
    \]
    such that the permuton \eqref{eqn-permuton-def} coincides with the skew Brownian permuton. The coupling is the \emph{imaginary geometry} coupling of Miller and Sheffield~\cite{ig4}: $\eta_1$ and $\eta_2$ are space-filling counterflow lines of angles $0$ and $\theta(q)-\pi/2$, with $\theta(q)\in(-\pi/2,\pi/2)$, of the same whole-plane Gaussian free field. The relation between $\theta$ and $q$ is now known explicitly, thanks to \cite{ang2025p}:
    \begin{equation}\label{eq:q-theta-rel}
        \frac{q}{1-q}
        =
        \frac{\sin((1/2-\gamma^2/8)(\pi-2\theta))}
        {\sin((1/2-\gamma^2/8)(\pi+2\theta))}.
    \end{equation}
	
	\item\label{item:sep_perm}
    At the boundary $\rho=1$, the skew Brownian permutons $\bm{\mu}_{1,q}$ coincide with the biased Brownian separable permutons $\bm{\mu}_{1-q}$ of \cite{maazoun-separable-permuton}. The connection between Brownian separable permutons and SLE/LQG is still conjectural. It should be a degenerate version of the construction in~\eqref{eqn-permuton-def}, where the SLEs are replaced by two explorations of the same conformal loop ensemble CLE$_4$. The original construction in~\cite{maazoun-separable-permuton} is in terms of a Brownian continuum random tree (CRT) whose internal branch points are decorated with i.i.d.\ $\oplus$ and $\ominus$ signs; in the notation used here, $\mathbb P(\ominus)=q=1-p=1-\mathbb P(\oplus)$.
\end{enumerate}

We point out that when $\kappa = 16/\gamma^2$, a $\gamma$-Liouville quantum surface decorated by an independent space-filling SLE$_\kappa$ enjoys a number of special properties arising from the mating-of-trees theory of~\cite{wedges}; see also \cite{ghs-mating-survey} for a survey. Consequently, this case is more tractable than the general case when $\kappa \neq 16/\gamma^2$, sometimes referred to as the ``mismatched" case.
Skew Brownian permutons fall within this framework, whereas the meandric permuton does not.

\subsection{Permuton limit results}
\label{sect:scaling-limit-results}

We have introduced three families of permutations: meandric, Baxter, and separable permutations. We have also introduced three canonical random fractal permutons: the meandric permuton, the skew Brownian permutons, and the Brownian separable permutons. The next three statements explain how these objects are (or are thought to be) connected.

\begin{conj}[{\cite[Conjecture 1.8]{bgs-meander}}]\label{conj:mean_conv}
	Let $\sigma_n$ be a uniform meandric permutation of size $2n$ fixing 1, equivalently the permutation associated with a uniform meander of size $n$. Then, with respect to the permuton topology,
	\[
	\bm{\mu}_{\sigma_n}\xrightarrow{d} \bm{\pi}_{m},
	\]
	where $\bm{\pi}_{m}$ is the meandric permuton introduced in \Cref{item:mean_perm}.
\end{conj}

\begin{thm}[{\cite[Theorem 1.9]{bm-baxter-permutation}}]\label{thm:baxt_conv}
	Let $\sigma_n$ be a uniform Baxter permutation of size $n$. Then, with respect to the permuton topology,
	\[
	\mu_{\sigma_n}\xrightarrow{d} \bm{\mu}_{-1/2,1/2},
	\]
	where $\bm{\mu}_{-1/2,1/2}$ is the skew Brownian permuton with parameters $\rho=-1/2$ and $q=1/2$ introduced in \Cref{item:skew_perm}. This special skew Brownian permuton is also called the \textbf{Baxter permuton}.
\end{thm}

\begin{thm}[{\cite[Theorem 1.6]{bassino-separable-permuton}}]\label{thm:sep_conv}
	Let $\sigma_n$ be a uniform separable permutation of size $n$. Then, with respect to the permuton topology,
	\[
	\mu_{\sigma_n}\xrightarrow{d} \bm{\mu}_{1,1/2},
	\]
	where $\bm{\mu}_{1,1/2}$ is the Brownian separable permuton with bias parameter $q=1-p=1/2$ introduced in \Cref{item:sep_perm}.
\end{thm}

The Brownian separable permutons $\bm{\mu}_{1,q}$ have been shown to be the permuton limits of many other families of pattern-avoiding permutations~\cite{Bassino-universal-perm,Borga-decorated-perm}. Similarly, the skew Brownian permutons $\bm{\mu}_{\rho,q}$ have been shown to be the permuton limits of several families of vincular-pattern-avoiding permutations~\cite{borga2021permuton, Davit-26-joint-conv}. The appearance of the same continuum objects in different discrete models is the main reason for calling these permutons universal.

\subsection{Properties and further results}

We now collect a few structural results about these limiting permutons. They show that the skew Brownian and Brownian separable permutons are highly singular random measures, while still having computable averaged statistics and rich pattern density behavior.

\subsubsection{Intensity measures}

The skew Brownian permutons $\bm{\mu}_{\rho,q}$ are random probability measures on the unit square with uniform marginals. Recall from Section~\ref{subsec:intensity-measures} that, with such a random permuton, we can associate its intensity measure $\mathbb E[\bm{\mu}_{\rho,q}]$.
The next two results give explicit density formulas for the intensity measure in the Brownian separable and Baxter cases.

\begin{thm}[{\cite[Theorem 1.7]{maazoun-separable-permuton}}]
    For all $q\in(0,1)$, the intensity measure $\mathbb E[\bm{\mu}_{1,q}]$ of the Brownian separable permuton is absolutely continuous with respect to Lebesgue measure on $[0,1]^2$. Moreover, $\mathbb E[\bm{\mu}_{1,q}]$ has density
\begin{equation*}
	p^{1-q}_S(x, y) =
    \int_{\max\{0,x+y -1\}}^{\min\{x,y\}}
    \frac{3q^2(1-q)^2 \,da}
	{2\pi(a(x-a)(1-x-y+a)(y-a))^{3/2}
    {\left(\frac{(1-q)^2}{a}+\frac{q^2}{x-a}+\frac{(1-q)^2}{1-x-y+a}+\frac{q^2}{y-a}\right)^{5/2}}}.
\end{equation*}
\end{thm}

\begin{thm}[{\cite[Theorem 1.3]{bhsy-baxter-permuton}}]\label{thm-baxter-density}
	Write $\bm{\mu}_B:=\bm{\mu}_{-1/2,1/2}$ for the Baxter permuton. Define
	\begin{equation}\label{eqn-rho}
		\rho(t, x, r)
        :=
        \frac{1}{t^2}
        \left(
        \left(\frac{3rx}{2t}-1\right)e^{-\frac{r^2+x^2-rx}{2t}}
        +
        e^{-\frac{(x+r)^2}{2t}}
        \right).
	\end{equation}
	Then the intensity measure $\mathbb E[\bm{\mu}_B]$ is absolutely continuous with respect to Lebesgue measure on $[0,1]^2$. Moreover, it has density
	\begin{equation}\label{eqn-baxter-density}
		p_B(x, y)
        =
        c\int_{\max\{0, x+y-1\}}^{\min\{x, y\}}
        \int_{\mathbb R_{+}^4}
        \rho(y-z, \ell_1, \ell_2)
        \rho(z, \ell_2, \ell_3)
        \rho(x-z, \ell_3, \ell_4)
        \rho(1+z-x-y, \ell_4, \ell_1)
        \, d\ell_1 d\ell_2 d\ell_3 d\ell_4\, dz,
	\end{equation}
	where $c$ is a normalizing constant.
\end{thm}

The integral in \eqref{eqn-baxter-density} is difficult to simplify further since it involves integrating a four-dimensional Gaussian expression over the positive quadrant. One can also obtain formulas for $\mathbb E[\bm{\mu}_{\rho,q}]$ for general $\rho$ and $q$, but the resulting expressions are substantially less explicit; see~\cite[Theorem 3.8]{bhsy-baxter-permuton}.

The proof of \cite[Theorem 1.7]{maazoun-separable-permuton} relies on the construction of the Brownian separable permutons from the Brownian CRT with i.i.d.\ signs. To the best of our knowledge, this proof does not easily extend to the Baxter permuton. The proof of \Cref{thm-baxter-density} instead uses two main ingredients:
\begin{itemize}
	\item the connection between the Baxter permuton, and more generally the skew Brownian permutons, and SLE curves coupled with LQG measures;
	\item exact integrability results for SLE curves and LQG measures, in particular results from the theory of conformal welding of quantum surfaces~\cite{AHS20}.
\end{itemize}

\begin{figure}[t]
	\centering
	\includegraphics[width=0.23\textwidth]{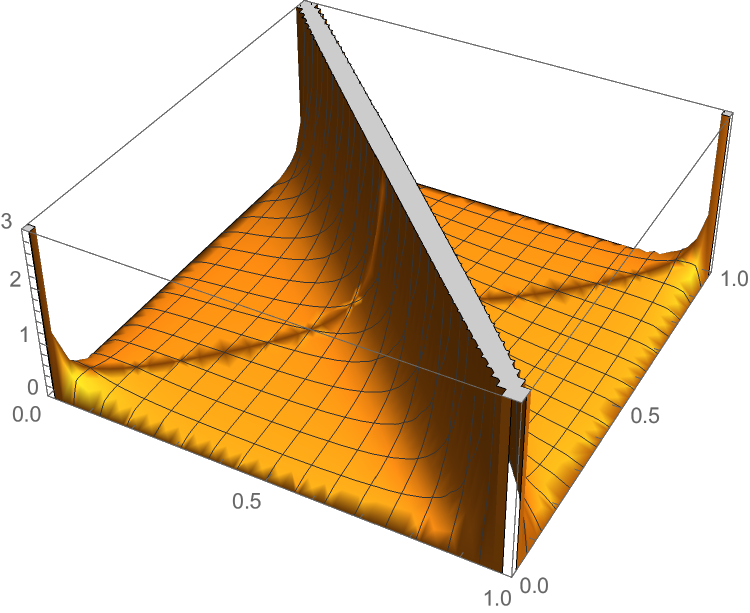}
	\includegraphics[width=0.23\textwidth]{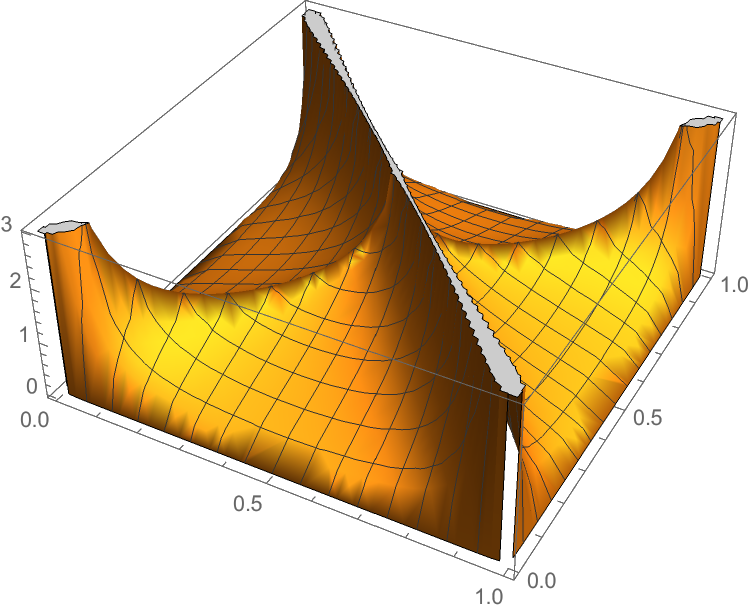}
	\includegraphics[width=0.23\textwidth]{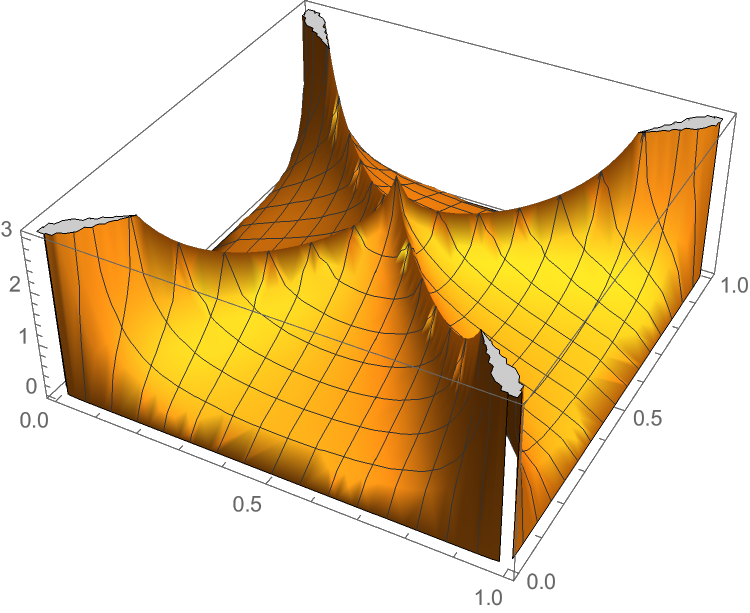}
	\includegraphics[width=0.23\textwidth]{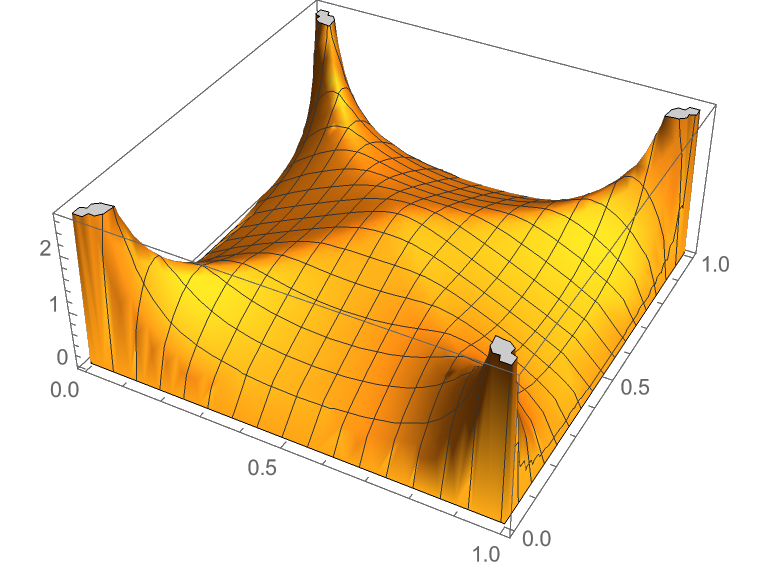}
	\caption{\textbf{From left to right:} Plots of the densities $p_S^{0.1}(x,y)$, $p_S^{0.4}(x,y)$, $p_S^{0.5}(x,y)$, and $p_B(x,y)$.}\label{fig:perm_densities1}
\end{figure}

Explaining these results {is outside the scope of this survey}. We therefore conclude this discussion with numerical plots of $p_B(x,y)$ and $p^q_S(x,y)$; {see Fig.~\ref{fig:perm_densities1}.}

\subsubsection{Pattern occurrences and the proportion of inversions}

The next result shows that there is a sharp difference between the cases $\rho\neq 1$ and $\rho=1$. When $\rho\neq 1$, the skew Brownian permuton $\bm{\mu}_{\rho,q}$ almost surely contains every classical pattern with positive density. By contrast, the Brownian separable permuton $\bm{\mu}_{1,q}$ avoids all non-separable patterns; see \cite[Definition 5.1]{Bassino-universal-perm}.

\begin{thm}[{\cite[Theorem 1.10]{bhsy-baxter-permuton}}]\label{thm:positivity}
	For all $(\rho,q)\in(-1,1)\times(0,1)$ and all patterns $\pi\in\Sym$, it holds that
	\[
		d_\pi(\bm{\mu}_{\rho,q})>0
        \qquad \text{a.s.}
	\]
\end{thm}

As an immediate consequence, if a sequence of random permutations from a proper pattern-avoiding class 
converges to a skew Brownian permuton, then the limit must lie on the Brownian separable boundary.
This is the case of the class of separable permutations, and of many other permutation classes 
studied in~\cite{Bassino-universal-perm,bassino2019scaling}.
{On the other hand, permutation families defined by avoidance of {\em generalized} patterns, such as Baxter, strong-Baxter
or semi-Baxter permutations, may contain all {\em classical} patterns and hence may converge to some
$\bm{\mu}_{\rho,q}$ with $\rho<1$.}

Using SLE and LQG techniques, one can also compute the expected density of inversions in the skew Brownian permutons.

\begin{figure}[t]
	\centering
	\includegraphics[width=0.23\textwidth]{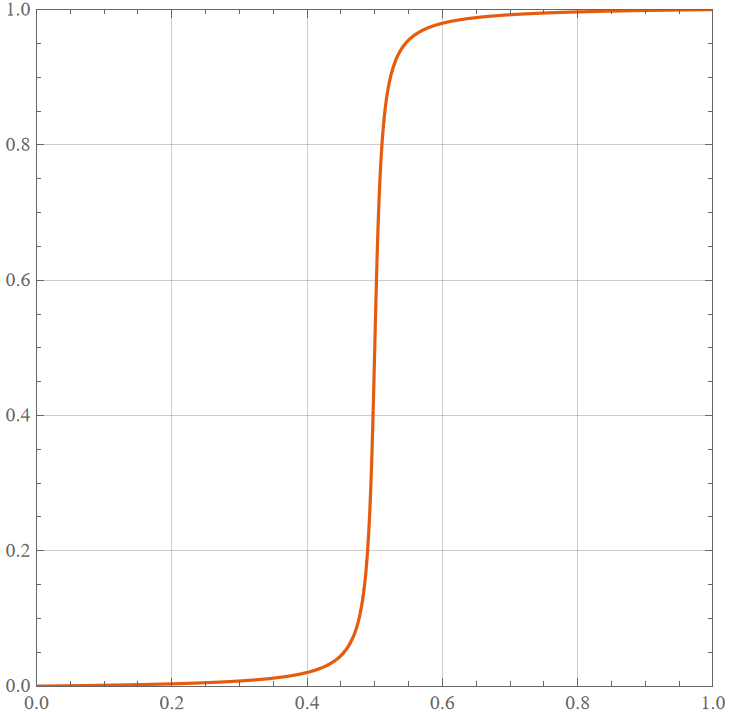}
	\includegraphics[width=0.23\textwidth]{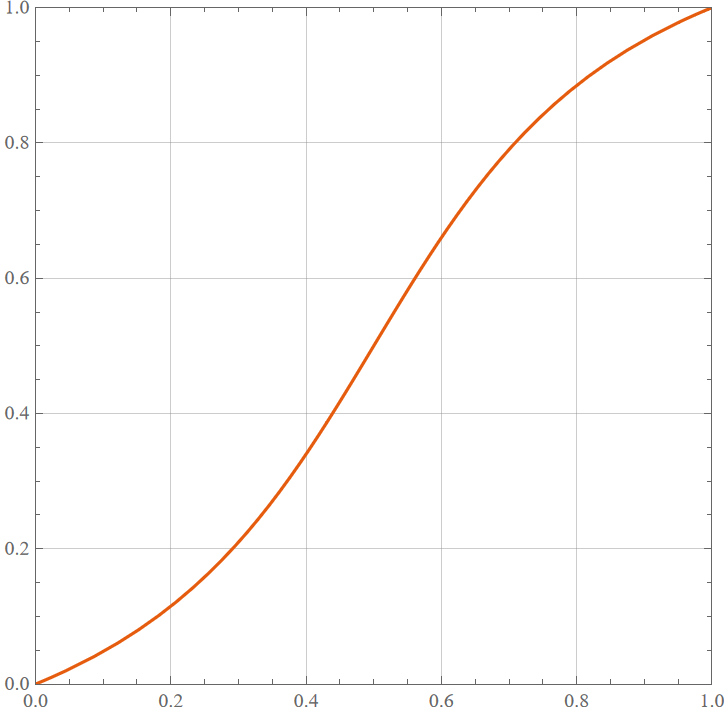}
	\includegraphics[width=0.23\textwidth]{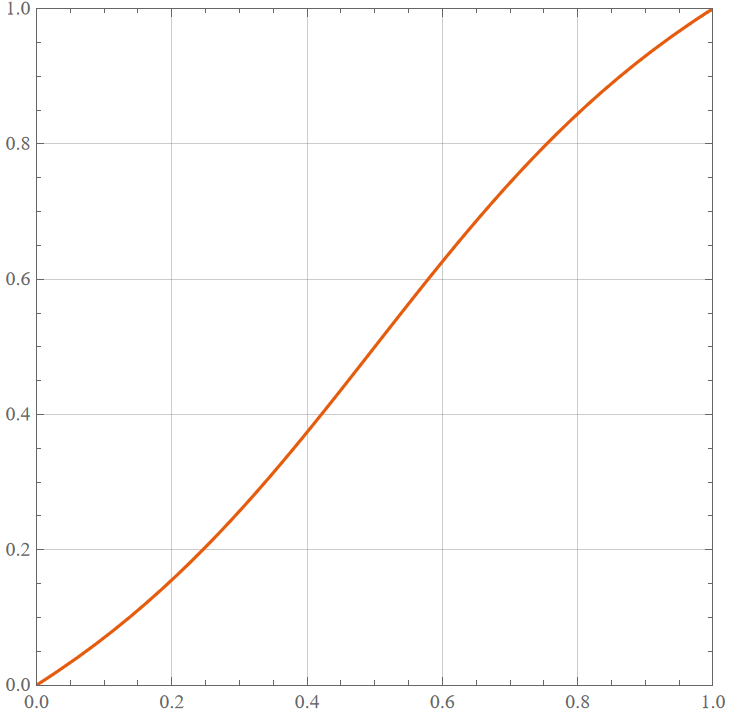}
	\includegraphics[width=0.23\textwidth]{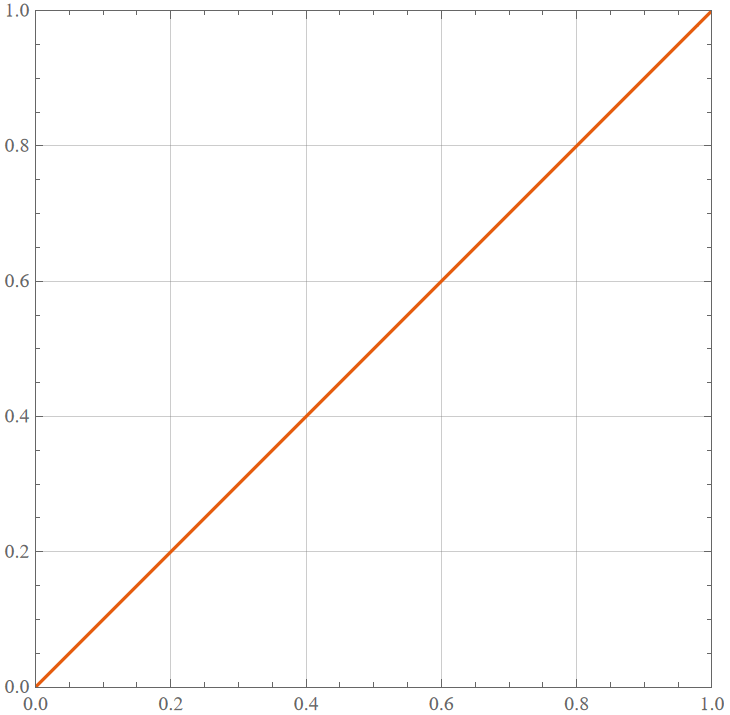}
	\caption{\textbf{From left to right:} The expected proportion of inversions $\mathbb E[d_{21}(\bm{\mu}_{\rho,q})]$ in~\eqref{eq:inversion}, as a function of $q$, for $\rho=-0.999,\,-0.5,\,0,\,1$.}\label{fig:inv_prop}
\end{figure}

\begin{thm}[{\cite[Proposition~1.14]{bhsy-baxter-permuton}}]\label{prop-expected-occ-21}
	Let $(\rho,q)\in(-1,1)\times(0,1)$, and let $\gamma\in(0,2)$ be determined by
    \[
        \rho=-\cos(\pi\gamma^2/4).
    \]
    Let $\theta(q)\in(-\pi/2,\pi/2)$ be related to $q$ by \eqref{eq:q-theta-rel}. Then
	\begin{equation}\label{eq:inversion}
		\mathbb E[d_{21}(\bm{\mu}_{\rho,q})] = \frac{\pi-2\,\theta(q)}{2\pi}.
	\end{equation}
\end{thm}
{In the limit $\rho \to 1$, we get $\mathbb E[d_{21}(\bm{\mu}_{1,q})] =q$, which can also be seen
easily from the original construction of the Brownian separable permuton.}
See \Cref{fig:inv_prop} for plots of the right-hand side of \eqref{eq:inversion} as a function of $q$ for several values of $\rho$.

\subsubsection{Hausdorff dimension of the support}

The support of these random permutons is fractal. A first result in this direction was obtained by Maazoun for the Brownian separable permutons. In the notation of this section, his theorem says that for every $q\in(0,1)$, the support of the biased Brownian separable permuton $\bm{\mu}_{1,q}$ is almost surely totally disconnected and has Hausdorff dimension one; in fact, its one-dimensional Hausdorff measure is bounded above by $\sqrt{2}$ \cite[Theorem 1.5]{maazoun-separable-permuton}. This proof uses the explicit Brownian excursion construction of the Brownian separable permuton: the permuton is realized as the pushforward of Lebesgue measure by the graph of a random measure-preserving function associated with a signed Brownian excursion.
{It was also observed by Riera that, in fact, any permuton limit of a sequence of permutations avoiding a fixed pattern $\tau$
has Hausdorff dimension $1$. This follows from a structural sparsity property of pattern-avoiding permutations developed to prove
the Stanley-Wilf conjecture, see \cite{marcus2004StanleyWilf}.}

For a general permuton $\mu$, the Hausdorff dimension of $\operatorname{supp}(\mu)$ is always at least one. Indeed, since the first marginal of $\mu$ is Lebesgue measure, the projection of $\operatorname{supp}(\mu)$ onto the first coordinate is all of $[0,1]$, and coordinate projection is Lipschitz. Thus, dimension one is the smallest possible value for the support of a permuton. On the other hand, general permutons can have supports of dimension anywhere between $1$ and $2$.

A theorem of Borga, Gwynne, and Sun shows that the same minimal-dimension phenomenon holds throughout the SLE/LQG family of permutons introduced in \Cref{sect:space-filling-constr}. 

\begin{thm}[{\cite[Theorem 1.20]{bgs-meander}}]\label{thm:sle-lqg-support-dim}
	Assume that we are in the setting of \Cref{sect:space-filling-constr}, with arbitrary $\gamma\in(0,2)$, $\kappa_1,\kappa_2>4$, and arbitrary coupling of $\eta_1$ and $\eta_2$, viewed modulo time parametrization. Let $\bm{\pi}$ be the associated SLE/LQG permuton. Then, almost surely,
    \[
        \dim_{\mathrm H}\operatorname{supp}(\bm{\pi})=1.
    \]
\end{thm}

In particular, the closed support of every skew Brownian permuton $\bm{\mu}_{\rho,q}$ with $(\rho,q)\in(-1,1)\times(0,1)$ has Hausdorff dimension one, and the same is true for the meandric permuton. Together with Maazoun's theorem for the Brownian separable case $\rho=1$, this gives a unified picture: all the principal random permutons in \Cref{fig-Venn_permutons} have supports of the smallest possible Hausdorff dimension.

The proof of \Cref{thm:sle-lqg-support-dim} is quite different from Maazoun's proof. Let
\[
    T(\eta_1,\eta_2)
    :=
    \big\{(t,s)\in[0,1]^2:\eta_1(t)=\eta_2(s)\big\}.
\]
This set records all pairs of times at which the two space-filling curves visit the same point.
Since $\operatorname{supp}(\bm{\pi})\subseteq T(\eta_1,\eta_2)$, it is enough to prove an upper bound on the Hausdorff dimension of $T(\eta_1,\eta_2)$. The key estimate controls the probability that the images of two small time intervals under $\eta_1$ and $\eta_2$ intersect. Such an estimate is obtained using metric estimates for Liouville quantum gravity. This is why the result applies not only to the skew Brownian permutons, but also to the meandric permuton and to more general SLE/LQG permutons.

\subsubsection{Permutons determine LQG and SLEs}

The construction in \Cref{sect:space-filling-constr} starts from a rich geometric object: a Liouville quantum gravity surface decorated by two space-filling SLE curves. It then produces a permuton on the unit square by recording, roughly speaking, the correspondence between the time order of the first curve and the time order of the second curve. A priori, one might expect this operation to discard a large amount of information. The permuton does not explicitly remember the locations of the points in $\mathbb C$, the random field $h$, or the two curves as subsets of the plane.

Surprisingly, in the main examples, this loss of information does not occur. Borga and Gwynne prove that, in two natural settings, the permuton determines the entire SLE-decorated LQG surface, up to some simple symmetries~\cite[Theorem 1.1]{borga2024reconstructing}. 

This reconstruction theorem gives another reason why the permuton viewpoint is powerful. In models such as Baxter permutations or, conjecturally, meandric permutations, convergence of the permutation diagrams to a random permuton can be interpreted not merely as a weak limit for a measure on the square, but as a form of convergence toward an SLE-decorated LQG surface.

\subsubsection{Longest increasing subsequences}\label{sect:univ-LIS} 
{We now discuss, in the spirit of Section~\ref{sect:lis-permutons}, the asymptotic behavior of the LIS of permutations $\sigma_n$ converging to some of the fractal permutons discussed here.

The first results of this kind} were
proved by Bassino, Bouvel, Drmota, F\'eray, Gerin, Maazoun, and Pierrot~\cite{bassino2022linear} for the Brownian separable permutons $\bm{\mu}_{1,q}$ and by Borga, Gwynne, and Sun~\cite{bgs-meander} for {general permutons fitting in the framework of ~\eqref{eqn-permuton-def}}.

\begin{thm}[{\cite[Theorem 1.10]{bassino2022linear},\cite[Theorem 1.13]{bgs-meander}}]
\label{thm:bgs-sublinear-lis}
Let $\bm{\pi}$ be either\footnote{The original statement in \cite[Theorem 1.13]{bgs-meander} covers some more general cases.} the meandric permuton or the skew Brownian permuton $\bm{\mu}_{\rho,q}$ with $(\rho,q)\in(-1,1]\times(0,1)$. Let $\sigma_n\in\Sym_n$ be a sequence of random permutations {converging to $\bm{\pi}$.}
Then
\[
    \frac{\LIS{\sigma_n}}{n}
    \xrightarrow[n\to\infty]{\prob}0.
\]
\end{thm}
{Note that the above theorem includes the case $\rho=1$, i.e.~of a sequence of permutations converging to a Brownian separable permuton.
The proof is different in the latter case than for $\rho<1$ and the meandric permuton,
but both use} relatively soft arguments based on the self-similar nature of the permutons. The theorem applies, for example, to uniform separable, Baxter, strong-Baxter and semi-Baxter permutations; all these models have sublinear LIS. Moreover, it gives strong evidence for the sublinear behavior of the LIS of uniform meandric permutations; this would follow from~\Cref{conj:mean_conv}. 

{The weakness of \Cref{thm:bgs-sublinear-lis} is that} it does not identify the correct sublinear exponent.
{As discussed in Section~\ref{sect:lis-permutons}, for a general sequence  $\sigma_n$ converging to a skew Brownian permuton $\mu_{\rho,q}$, we cannot expect a better result
than \Cref{thm:bgs-sublinear-lis}, but for specific sequences, in particular for the permuton samples $\sigma^{\mu_{\rho,q}}_n$, we can ask about the precise behavior of 
$\LIS{\sigma^{\mu_{\rho,q}}_n}$.}

The case $\rho=1$ {is the best understood}.  We write $\bm{\mu}_p$ for the Brownian separable permuton with parameter $p\in(0,1)$, where $p$ is the probability of a $\oplus$ sign in the CRT construction. {This corresponds to the skew Brownian permuton $\mu_{1,1-p}$.}
A first substantial refinement of \Cref{thm:bgs-sublinear-lis} {for $\sigma^{\mu_{p}}_n$} was obtained by Borga, Da Silva and Gwynne~\cite{borga2024power}.
Using Brownian excursion computation and fragmentation theory, 
they proved polynomial upper and lower bounds 
for $LIS(\sigma^{\mu_{p}}_n)$.
Later, Adhikari, Borga, Budzinski, Da Silva, and
S\'enizergues identified the exact exponent and proved a scaling limit result.
\begin{thm}[{\cite[Theorem 1.1]{adhikari2025longest}}]
\label{thm:brownian-separable-exact-lis}
Let $p\in(0,1)$, and let $\bm{\mu}_p$ be the Brownian separable permuton
with parameter $p$.  Conditionally on $\bm{\mu}_p$, let $\sigma_n$ be a
permutation of size $n$ sampled from $\bm{\mu}_p$.  
Let also  $\alpha(p)$ be the unique solution in $(1/2,1)$ of
\[
    -\frac{1}{4^{1/(2\alpha)}\sqrt{\pi}}\,
    \frac{\Gamma\!\left(\frac12-\frac{1}{2\alpha}\right)}
         {\Gamma\!\left(1-\frac{1}{2\alpha}\right)}
    =
    \frac{p}{1-p}.
\]
Then there exists a
non-degenerate random variable $X_p$, which is almost surely positive and
finite and is measurable with respect to $\bm{\mu}_p$, such that
\[
    \frac{\LIS{\sigma_n}}{n^{\alpha(p)}}
    \xrightarrow[n\to\infty]{\mathrm{a.s.}}
    X_p.
\]
\end{thm}
The function $\alpha(p)$ is continuous and increasing, with
$
    \alpha(0^+)=\frac12$,
$
    \alpha(1^-)=1$,
    and
    $\alpha(1/2)\approx0.815$.
The value $p=1/2$ (the standard symmetric case) corresponds to the limit of uniform separable permutations. 
Thus, {the LIS of a sample of
the Brownian separable permuton in this case is approximately $n^{0.815}$,} 
strictly between the classical $\sqrt n$ exponent and the linear scale. 
It is expected that analogous
asymptotics hold for many natural discrete models converging to
$\bm{\mu}_p$, including uniform separable permutations. However, transferring 
such precise sublinear asymptotics from the limiting permuton sample back to a
given discrete model requires a substantial amount of technical work.
\medskip

For general skew Brownian permutons $\bm{\mu}_{\rho,q}$,
only an upper bound is known.
\begin{thm}[\cite{bdg-26-LIS-sbp}]
\label{thm:skew-exact-lis1}
Let $\bm{\mu}_{\rho,q}$ be the skew Brownian permuton
with parameters $(\rho,q)\in(-1,1)\times(0,1)$.  
Then there exists $\alpha^*(\rho,q)<1$ such that
\[
    \frac{\LIS{\sigma^{\bm{\mu}_{\rho,q}}_n}}{n^{\alpha^*(\rho,q)}}
    \xrightarrow[n\to\infty]{\mathbb{P}}
    0.
\]
\end{thm}
Interestingly, there is a curve in the parameter space $(\rho,q)$ on which the LIS is more tractable.
\begin{thm}[\cite{bgsp-26-LIS-sbp}]
\label{thm:skew-exact-lis2}
Let $(\rho,q)\in(-1,1)\times(0,1)$ verify 
\begin{equation}\label{eq:good_parameters}
    q=\frac{2\rho}{2\rho-1} \quad\text{and}\quad  \rho< 0.
\end{equation}
Then with probability tending to one as $n\to\infty$, we have
$\LIS{\sigma^{\bm{\mu}_{\rho,q}}_n}=n^{\alpha+o(1)}$,
where 
\begin{equation}
    \alpha=\alpha(\rho):=\frac{2}{4-\gamma(\rho)^2},\qquad\text{with $\gamma(\rho)\in(0,\sqrt 2)$ defined by\ $\rho = -\cos(\pi \gamma^2/4)$.}
\end{equation}
\end{thm}
Fortunately, the parameters of the limiting permutons of Baxter and strong-Baxter (but not semi-Baxter) permutations satisfy \eqref{eq:good_parameters}.
In particular, the Baxter case corresponds to $(\rho,q)=(-1/2,1/2)$, and the corresponding exponent is $\alpha(-1/2)=3/4$; while the strong-Baxter case corresponds to
$(\rho,q)=(\frac{r-1}{2},\frac{1-r}{2-r})$, where $r\in(0,1)$ is the unique solution of $r^3=(1-r)^2$ and the corresponding exponent is $$\alpha_{sb}=\frac{1}{2-\frac{2}{\pi}\arccos\left(\frac{1-r}{2}\right)}\approx 0.879.$$
The next theorem transfers the results from permuton samples to natural combinatorial models with the same limits, i.e.~uniform Baxter and strong-Baxter permutations.
\begin{thm}
[\cite{bgsp-26-Tightness}]
\label{thm:skew-exact-lis3}
Let $\tau_n$ be a uniform Baxter permutation of size $n$.
With probability tending to one as $n\to\infty$, we have
\[
    \LIS{\tau_n}=n^{3/4+o(1)}.
\]
Let $\tau^{sb}_n$ be a uniform strong-Baxter permutation of size $n$.
With probability tending to one as $n\to\infty$, we have
\[
\LIS{\tau^{sb}_n}=n^{\alpha_{sb}+o(1)}.
\]
\end{thm}

Baxter permutations are known to be in bijection with a family of directed planar maps (i.e. planar maps endowed with an orientation of the edges) called \emph{bipolar orientations}~\cite{bbf-baxter-bipolar}. Under this bijection, an LIS of a Baxter permutation is mapped to a longest directed path from the source to the sink of the corresponding bipolar orientation. This connection is much more general: roughly speaking, it is conjectured that the LIS of the skew Brownian permutons is tightly connected to a new notion of directed distances for Liouville quantum gravity surfaces~\cite{bgsp-26-Tightness},
and the above Theorems~\ref{thm:skew-exact-lis1}--\ref{thm:skew-exact-lis3} are part of a general program aimed at understanding directed distance in planar maps initiated in~\cite{bg-2025-directed}.

\section{Companion theories and open problems}\label{sect:open-prob}

Let us briefly consider several companion limit theories that share the general philosophy of permutons.

\subsection{Companion theories: graphons}\label{sect:open-prob-1}
{
As mentioned in the introduction, permutons were first introduced by analogy with graphons; we refer the reader to~\cite{Lov} for a comprehensive account of graph limit theory,
and to~\cite{Meliot} for a direct comparison of the two theories.
	
A graphon is {a symmetric measurable function} from $\Omega^2$ to $[0,1]$, {where $\Omega$ is a probability space}, {considered up to weak isomorphism}
{(typically $\Omega=[0,1]$, but sometimes other choices are more natural)}. In {graphon theory} the role of pattern densities is played by {\em homomorphism densities} {$t(F,G)$}, which can be defined as the probability that a uniformly random {map from $V(F)$ to $V(G)$ sends every edge of $F$ to an edge of $G$}. {One can alternatively define induced-subgraph densities, by asking in addition}
{that every non-edge of $F$ is sent to a non-edge of $G$}.
 Both families of densities {are related by a triangular change of basis formula,
 and one can use either of them to} characterize dense graph convergence.
	
	A natural (and useful) way to convert a permuton to a graphon extends the following map from
	permutations to graphs, {known as (the complement of) its inversion graph}.  Given $\pi \in S_n$, regarded as a total (re)-ordering of $[n] :=
	\{1,2, \dots, n\}$, intersect that ordering with the natural order on $[n]$ to get a partial
	ordering $P$ of $[n]$.  The resulting {\em comparability graph} $G_P$ on $[n]$, in which
	$i \sim j$ iff {either $i < j$ and $\pi(i) < \pi(j)$, or $j < i$ and $\pi(j) < \pi(i)$}, is our target.  Notice that
	the 12-density of $\pi$ is the edge-density of $G_P$, and more generally the density of the
	pattern $1 2 \dots k$ in $\pi$ is the $k$-clique density {(defined by sampling $k$ distinct vertices)}.
	{For a permuton $\mu$, the corresponding graphon can be represented by}
	{\[
		W_\mu((x,y),(x',y')) = \mathbf{1}_{\{(x-x')(y-y')>0\}},
		\]}
		{where we use $\Omega=[0,1]^2$ instead of the usual $\Omega=[0,1]$ choice.}
		
		{Some pattern densities of $\mu$ can be related to subgraph densities
		of $W_\mu$.}
	{In particular, $d_{12}(\mu)=t(K_2,W_\mu)$, $d_{123}(\mu)=t(K_3,W_\mu)$, and $d_{321}(\mu)=t(K_3,1-W_\mu)$.}
	It follows that the permuton feasible regions for certain sets of patterns are subsets of the
	graphon feasible regions for the corresponding graphs.  Of course most graphons are not the image
	of some permuton in the above-described mapping, but it happens that the feasible regions
	for patterns 12 vs 123 and for patterns 123 vs 321, shown in Fig.~\ref{fig:feas-ref}, {coincide with}
	{the corresponding graphon feasible regions for edge versus triangle density and triangle versus independent-triple density, respectively. The first region was characterized by} Razborov \cite{Raz}
	{and the second coincidence was established by} Elizalde and Noy \cite{ElizNoy}
	{(see also~\cite{Kenyon-Permutations-densities})}.
	
	Graphons have their own variational principle, due to Chatterjee and Varadhan \cite{ChatVar},
	and {Kenyon, Radin, Ren, and Sadun}~\cite{KenRadRenSad} (among other works cited therein) use this {framework to study entropy-maximizing graphons for prescribed edge and $k$-star densities}.
	However, even when feasible regions coincide, entropy-maximization problems
	on graphons and permutons are different.
	In several constrained models, graphon entropy maximizers 
	are constant on every $A_i\times A_j$ for a finite measurable partition $\Omega=\bigsqcup_i A_i$ (they are called {\em multipodal}), while the maximum-entropy permutons (at interior points) {of the $(12,123)$ feasible region} seem to have smooth, analytic density functions.  Whether this contrast
	persists in other feasibility regions, or indeed even in more general circumstances, remains to be seen.}

\subsection{Companion theories: beyond graphons}

{Beyond graphs, {limiting theories have been developed for a variety of combinatorial objects in the past twenty years:  directed graphs and tournaments~\cite{diaconis2008graphons,thornblad2018tournament_limits}, posets~\cite{janson2011poset}, hypergraphs~\cite{elek2012limit_hypergraphs},
 order  types~\cite{goaoc2015limit_order_types}, trees~\cite{dendrons}, Latin squares~\cite{garbe2020limits}, bounded operators~\cite{backhausz2022action_convergence},
non attacking queen configurations~\cite{simkin2023queens}, higher dimensional permutations~\cite{borga2024high_dimensional}, to name a few. These theories describe large discrete objects by analytic limits (functions, measures, metric spaces, \dots). In many of the dense combinatorial settings, convergence is characterized by densities of finite induced substructures; the tree and operator theories use other sampling or convergence notions. Coregliano and Razborov~\cite{CR20} developed a general framework for dense combinatorial limits, under the name of {\em theons}. We discuss in this section theories that are directly linked with permutations, as well as another natural notion of convergence for permutations that captures information at a different scale.}}

While permuton convergence describes the macroscopic geometry of a permutation diagram, {\em local convergence} {focuses on} the neighborhood of a uniformly chosen entry. 
This topology is characterized by the convergence of \emph{consecutive} pattern densities, and its random limiting objects satisfy a natural shift-invariance property \cite{borga2020local}. 
The corresponding feasible regions for consecutive pattern densities are particularly explicit: for patterns of any fixed size, the feasible region of consecutive pattern densities is the cycle polytope of an associated overlap graph~\cite{borga2020feasible}.
In particular, its dimension, vertices, faces, and defining equations can be determined.
Finally, we refer to~\cite{feray26-LDP-consecutive} and references therein for large deviation questions related to local convergence and consecutive patterns in permutations.

In a different direction, the permuton theory has been extended to higher-dimensional permutations. A $d$-dimensional permuton is a probability measure on $[0,1]^d$ whose one-dimensional marginals are uniform, and the basic theory of permuton convergence, sampling, and pattern densities extends to this setting~\cite{borga2024high_dimensional}. This theory has already been applied to Schnyder-wood permutations and to $d$-dimensional separable permutations, producing random limits related to the Brownian and skew Brownian permutons discussed in Section~\ref{sect:univ-perm-skew}; {and more recently (\cite{EPP26}) to uniform random antichains in $[0,1]^d$.}

Coming back to dimension $2$, but leaving the world of permutations by allowing repeated entries and positions leads to the theory of \emph{tablons}:
 the latter has been developed to describe limits of inversion sequences, words over possibly infinite alphabets, and, more broadly, two-dimensional tables of integers~\cite{tablons}.
 Here again, convergence can be characterized through the densities of finite substructures, and interesting quasi-randomness and continuity problems appear.

{Another application of this philosophy is Simkin's limit theory for the $n$-queens problem \cite{simkin2023queens}. An $n$-queens configuration is a permutation matrix with at most one nonzero entry on each diagonal. Its limit objects, called \emph{queenons}, are permutons with additional restrictions in the two diagonal directions. More precisely, after translating the underlying square to $[0,1]^2$, Simkin's space of queenons is the weak closure of step permutons (permutons whose two-dimensional density is constant on each cell of some finite uniform square grid) whose pushforwards under $(x,y)\mapsto x+y$ and $(x,y)\mapsto y-x$ are dominated by Lebesgue measure. An entropy variational principle yields $\log Q(n)=n\log n-\alpha n+o(n)$, where $Q(n)$ is the number of $n$-queens configurations and $\alpha$ is determined by the variational problem. The associated large deviation principle also identifies a deterministic limiting queenon for uniformly random configurations.}

 {Our last example} is the limit theory of Latin squares developed in~\cite{garbe2020limits}. 
The corresponding limit objects, called \emph{Latinons}, {are defined using permutons: this reflects the fact that each row/column of a Latin square is a permutation. 
As for permutons,} one can define densities of finite patterns {which characterize convergence}, 
the resulting limit space is compact, and every Latinon can be approximated by finite Latin squares. 
{Moreover, a quasirandomness criterion similar to those of Section~\ref{sect:quasirandomness} for permutations has been established in~\cite{cooper2022quasirandomLatin}.}

These developments show that permutons constitute a dynamic research field and a source of inspiration, 
even outside the {strict} setting of permutation combinatorics
and random permutations.

\subsection{Open problems}\label{ss:open_problems}
We collect here a list of open problems. For each problem, we mention the section where related results have been discussed.
\medskip

\begin{enumerate}
\setlength{\itemsep}{\medskipamount}
\item \label{quest:feasible} \textbf{Finite-dimensional feasible regions; Section~\ref{sect:feas-reg-1}.}
Given a finite family of patterns $T \subset \mathfrak{S}$, describe
\[
\mathcal{F}_T:=\left\{\bigl(d_\tau(\mu)\bigr)_{\tau \in T}: \mu \in \mathcal{P}\right\}.
\]
In particular, can one describe its boundary, its extremal points, and the permutons realizing them?
{See \cite{Kenyon-Permutations-densities,garbe2023dimensionfeasibleregionpattern} and the next open problem
for the few known results on this question.}

\item \textbf{Exact packing densities; Section~\ref{sect:packing-problem}.}
For a pattern $\tau$, determine{
\[p(\tau) = \lim_{n \to +\infty} \max_{\sigma_n \in S_n} d_\tau(\sigma_n),\]
or equivalently}
\[
p(\tau):=\max_{\mu \in \mathcal{P}} d_\tau(\mu).
\]
{This corresponds to the case $|T|=1$ of Problem~(\ref{quest:feasible}), since in this case $\mathcal{F}_{\{\tau\}}=[0,p(\tau)]$
(except if $\tau=1$ is the one-element permutation).}
{This permuton formulation of the packing problem~\cite{presutti2010packing}, and flag algebra methods~\cite{sliacan2017packing_flag},
 have permitted significant progress on the question, but exact values are still known for very few patterns only.
A related question is the structure of the maximizing permutons. In all cases that have been studied so far,
 the best known lower bounds come from permutons
whose support has only one point in almost all vertical or horizontal lines \cite[Table 2]{sliacan2017packing_flag}.
Is it the case of maximizing permutons in general,
or is there at least a finite bound on the number of points in the support in vertical/horizontal lines,
as for pattern-avoiding permutons (see~\cite{garbe2024avoiding})?}

\item \textbf{Minimal quasirandom-forcing sets; Section~\ref{sect:quasirandomness}.}
Determine the smallest size of a family $T$ of patterns such that
\[
d_\tau(\mu) = \frac{1}{|\tau|!},\qquad \text{for every } \tau \in T,
\]
implies that $\mu$ is the uniform permuton. 
A variant of this problem asks for a linear combination $\sum_{\tau \in T} \alpha_\tau d_\tau(\mu)$ with the smallest possible support size $|T|$
which characterizes quasirandomness.
{Partial results in this direction can be found in the papers \cite{kral2020quasirandom_pattern_sum,kurecka2022lower_bound_quasirandom,kral2024forcing}.}

\item \textbf{Finitely forcible permutons; Section~\ref{sect:quasirandomness}.}
Characterize the permutons that are uniquely determined by finitely many pattern densities.
In particular, which natural packing optimizers are finitely forcible?
{We refer to \cite{glebov15-forcible,glebov17-permuton-testing} for examples of finitely forcible and non finitely forcible permutons.}

\item \textbf{Exact supports of permutons; Section~\ref{sect:perm-support}.}
Coscia, Tassy, and Winkler~\cite{CTW26} show that a regular, closed set $S$ in the unit square is exactly the support of a permuton if and only if it is compliant and expansive.  But their characterization in the general case for closed $S$ (in terms of a ``shavability" condition) leaves a lot to be desired.  Is there a more geometric and/or more tractable one?

\item \textbf{Large deviations for permutation classes; Section~\ref{sect:LDP-sect}.}\label{pb:LDP_avoiding}
{The uniform measure on all permutations satisfies a permuton large deviation principle (LDP) at speed $n$ \cite{trashorras2008large_deviations}, and a corresponding LDP has been established for permutations sampled from a fixed permuton \cite{borga2022large}.
 More recently, LDPs at speed $n$ were proved for permutations avoiding a pattern of size three \cite{budzinski2026constrained_mallows}. These examples suggest that permutation classes {may} admit a broader large deviation theory, but the correct level of generality is not yet clear.}

Let $\mathcal{C}$ be a permutation class with exponentially many elements of each size, and let $\sigma_n$ be uniform in $\mathcal{C}\cap\mathfrak{S}_n$. Does the associated random permuton $\mu_{\sigma_n}$ satisfy a large deviation principle? Is speed $n$ universal for exponentially growing classes, or can other speeds occur? When an LDP exists, compute its rate function.

\item { \textbf{Towards a second-order permuton theory; Section~\ref{ssec:fluctuations}.}
\label{pb:fluctuations}
Establish results relating fluctuations of the cumulative distribution function of a sequence $\mu_n$ of random permutons around that of a limiting permuton $\mu$,
and fluctuations of the corresponding pattern densities $d_{\tau}(\mu_n)$ around $d_{\tau}(\mu)$. Links with the fluctuations of simply indexed permutation statistics,
which form another natural family of continuous functionals on the space of permutons,
 would be interesting as well.}

\item \textbf{Catropy and growth rates; end of Section~\ref{sect:pattern-av-perm}.}\label{pb:catropy}
{Recall that the catropy of a permuton is defined as $C(\mu) = \lim_{n \to \infty} \frac1n h(\sigma^\mu_n)$, if the limit exists, where $h(\sigma^\mu_n)$ is the Shannon entropy of the corresponding permuton sample.}
{For a geometric grid class $\mathcal C$ defined by a matrix $M$, compute and optimize $C(\mu)$ over the permutons obtained by varying the cell dimensions and masses in the corresponding geometric representation, with positive mass uniformly distributed along each prescribed segment. For which matrices is
\[
\sup_{\mu} C(\mu)
=\log_2\operatorname{gr}(\mathcal C),
\]
and when is the supremum attained?
More generally, find natural hypotheses ensuring that a singular permuton has finite catropy, and obtain formulas for its value, extending the results of~\cite{maga2025entropy}. What is the relationship between catropy and the exponential enumeration of permutations in small permuton neighborhoods? In which settings can catropy be used to express the large-deviation rate functions from  Problem~\eqref{pb:LDP_avoiding}?}

\item \textbf{Constrained Mallows models for general avoidance classes; \Cref{expem:cond-mallows}.}
Letting $B$ be a finite collection of forbidden patterns,
what is the probability that a random Mallows permutation (or other exponentially biased measures)
avoids all patterns in $B$? How does a large random Mallows permutation avoiding $B$ behave;
in particular, what is its limiting permuton?
Such problems are a mix between constraint models as in Section~\ref{sect:pattern-av-perm}, and biased models as in Section~\ref{sect:exponentially-biased-permutations}.
See \cite{Pinsky2021Mallows_avoiding,Chelikavada2026fixed_point_biased,Park2025fixed_point_biased,budzinski2026constrained_mallows} for recent results in this line of research.
Large deviation principles for permutation classes, as in the previous open Problem~\ref{pb:LDP_avoiding},
might be a useful tool to answer these questions in the regime $q_n=\exp(\beta a_n/n^2)$,
where $a_n$ is the speed of the corresponding LDP (in all known cases where such an LDP is known, $a_n=n$).

\item \textbf{Finite entropy, uniqueness, and analyticity in the relative
interior; \Cref{sect:Entropy-maximizing-permutons}.}
Fix a finite family $T$ of patterns and let
\[
\mathcal{F}_T:=\left\{\bigl(d_\tau(\mu)\bigr)_{\tau\in T}:\mu\in\mathcal{P}\right\}.
\]
For $\alpha$ in the (relative) interior of $\mathcal{F}_T$, is the constrained entropy
\[
s_T(\alpha):=\max\left\{H(\mu):\mu\in\mathcal{P},\quad d_\tau(\mu)=\alpha_\tau\ \text{for every }\tau\in T\right\}
\]
{always} finite? Under what additional hypotheses is the maximizer unique? (General uniqueness fails in the $123/321$ model of \cite[Section 10]{Kenyon-Permutations-densities}; see also the discussion below Theorem~\ref{thm-entr-ken}) When it is unique, must it have a positive smooth or real-analytic density? The evidence in the existing examples points in these directions, but the general statement is open.

\item \textbf{Enumeration on the boundary of feasible regions; \Cref{sect:Entropy-maximizing-permutons}.}
{Let us consider the problem of enumerating permutations with approximately prescribed pattern densities,
i.e.~the set $S_n^{ \bm \alpha,\eps}$ of Theorem~\ref{thm:fixed_pattern_densities_asymptotic_enumeration}.
The entropy formula of Kenyon, Král', Radin, and Winkler~\cite{Kenyon-Permutations-densities}
 gives the exponential correction to the $n!$ enumeration of all permutations.}
 On parts of the boundary of the feasible regions, where all realizing permutons are singular, however, the entropy is $-\infty$, and this formula no longer identifies the relevant lower-order scale.

Let $\alpha$ lie on the (relative) boundary of $\mathcal{F}_T$, especially in a case where every realizing permuton is singular and $s_T(\alpha)=-\infty$. Determine the correct asymptotic scale of
\[
\# S_n^{\alpha,\varepsilon}.
\]
Can one formulate a refined large deviation principle on an appropriate lower-dimensional 
or singular space of permutons? Pattern-avoiding permutations provide the basic test cases 
(recall Problem~\ref{pb:LDP_avoiding}).

\item \textbf{Permutation classes with trivial limiting permutons; \Cref{sect:pattern-av-perm}.}
While \Cref{sect:pattern-av-perm} focuses on nontrivial limiting permutons,
for many patterns $\tau$, a sequence $\sigma_n$
of uniform $\tau$-avoiding permutations converges to the diagonal or anti-diagonal permuton. It is for example the case for length 3 patterns~\cite{miner-pak-permuton}, all monotone patterns~\cite{hoffman2019scaling}, and
conjecturally for the hard pattern\footnote{The pattern 1324 is the only pattern $\tau$ (up to symmetries) of length 4 for which the number of $\tau$-avoiding permutations of size $n$ is unknown. Computing this number, or finding its asymptotic behavior is a celebrated open problem in the field.} 1324~\cite{madras10-ran-patt-av}, \dots
Although these diagonal permutons are arguably not the most fascinating limiting objects,
having a convergence result which applies to a large family of permutation classes
would be of great interest.

A conjecture, due to Troyka and generalizing previous conjectures of Bevan,
is the following~\cite[Conjecture 3 p.~282]{OberwolfachReportPP2024}. A pattern of size $k$ is called $\oplus$-decomposable if it maps a nontrivial initial interval $\{1,\dots,\ell\}$ with $0<\ell<k$ to itself. For example $2143$ is $\oplus$-decomposable since it maps $\{1,2\}$ to itself. All patterns starting with $1$, or ending with their maximal entry are also $\oplus$-decomposable. Now, let $\tau_1,\dots, \tau_m$ be a finite list of $\oplus$-decomposable patterns and $\sigma_n$ be a uniform random permutation conditioned to avoid all of $\tau_1,\dots, \tau_m$.
Then it is conjectured that $\sigma_n$ converges to the anti-diagonal permuton,
i.e.~to the law of $(U,1-U)$, where $U$ is uniform in $[0,1]$.

\item\label{quest:inv-RS-tableaux} \textbf{A continuous RS-correspondence; \Cref{sect:lis-permutons}.}
As explained in \Cref{sect:lis-permutons}, given a density $\rho$ on $[0,1]^2$,
Sjöstrand~\cite{sjostrand2023monotone} proved the existence of a limiting RS shape for the permutations $\sigma_n^\rho$ sampled from $\rho$ (after rescaling by $\sqrt n$ in both directions).
Applying Sjöstrand's results to appropriate restrictions of $\rho$,
one can improve this result and prove the existence of limiting surfaces for the RS {\em tableaux} associated with $\sigma_n^\rho$. This allows to associate with $\rho$
a pair $(f^\rho_P,f^\rho_Q)$ of $2$-dimensional surfaces representing the limiting tableaux.

This raises many fascinating questions. Sjöstrand's proof uses subadditivity results and does not allow computation of limits. Can we find some continuous analogue, e.g.~of the insertion procedure of RS correspondence, allowing more efficient computation of $(f^\rho_P,f^\rho_Q)$? Can this procedure be inverted? In other words, does this continuous RS correspondence define a bijection from permuton densities on $[0,1]^2$ to suitable pairs of surfaces?
Finally let us mention a beautiful conjecture from the paper of Sjöstrand \cite[Conjecture 3.5]{sjostrand2023monotone}: the RS shape associated to the uniform density on the diamond $\{|x-1/2|+|y-1/2| \le 1\}$ should simply be an isosceles triangle.

\item \textbf{Binary search trees of permuton samples; \Cref{ss:other-stats}.}
Let $\sigma^\mu_n$ denote the random permutation of size $n$ sampled from a permuton $\mu$, and $f$ be some permutation statistic. For some $f$, it is reasonable to believe that $\mathbb E[f(\sigma^\mu_n)]$ is minimized by the uniform permuton $\mu$. An example of such statistics is the height of the binary search tree associated with $\sigma^\mu_n$: we refer to Devroye~\cite{devroye1986note} for its behavior in the uniform case, and to Corsini--Dubach--Féray~\cite{corsini25-BST-samples} for results on permuton samples, but the problem of whether the uniform permuton asymptotically minimizes $\mathbb E[f(\sigma^\mu_n)]$ is open.

\item\label{quest:meandric-permuton} \textbf{The meandric permuton; \Cref{sect:scaling-limit-results}.}
Let $\sigma_n$ be the meandric permutation associated with a uniform meander of size $n$. Prove that
\[
\mu_{\sigma_n} \xrightarrow[n\to\infty]{d} \bm{\pi}_{\mathrm{m}},
\]
where $\bm{\pi}_{\mathrm{m}}$ is the meandric permuton. See~\Cref{conj:mean_conv}, {or \cite[Conjecture~1.8]{bgs-meander}}.

\item\label{quest:meandric-intrinsic} \textbf{An intrinsic description of the meandric permuton; \Cref{sect:scaling-limit-results}.}
The current definition of the meandric permuton $\bm{\pi}_{\mathrm{m}}$ 
uses Liouville quantum gravity and Schramm--Loewner evolution. 
Borga, Gwynne, and Sun ask in \cite[Problem~7.1]{bgs-meander} {for an
alternative construction or characterization.
For example, can $\bm{\pi}_{\mathrm{m}}$ be characterized directly through Brownian objects, a recursive distributional identity, or its rerooting symmetries and pattern densities? 
Finding such a characterization could be a key step in solving Problem~\ref{quest:meandric-permuton}.
See also~\cite{bw-26-meanders} for further discussions.}

\item \textbf{LIS exponents of skew Brownian permutons; Section~\ref{sect:univ-LIS}.}
Let $\bm{\mu}_{\rho,q}$ be a skew Brownian permuton with $(\rho,q)\in(-1,1]\times(0,1)$. Prove that there exists a deterministic exponent $\alpha(\rho,q)\in[0,1]$ such that
\[
\mathbb{P}\big(\LIS{\sigma_n^{\bm{\mu}_{\rho,q}}}=n^{\alpha(\rho,q)+o(1)}\big)\xrightarrow[n\to\infty]{} 1,
\]
and compute $\alpha(\rho,q)$. More generally, does there exist a nondegenerate scaling limit after normalization by $n^{\alpha(\rho,q)}$? 
This problem was first formulated in the paper \cite{borga2024power} (see Conjecture~1.7 there).
Results in the Brownian separable case ($\rho=1$) and for $(\rho,q)$ on the curve $\rho=\frac{q}{2(q-1)}$ with $\rho\leq 0$ have been obtained since then; 
see \cite{adhikari2025longest,bgsp-26-LIS-sbp}
or  Theorems~\ref{thm:brownian-separable-exact-lis} and~\ref{thm:skew-exact-lis2} here.

\item \textbf{LIS exponents and Hausdorff dimension of increasing subsets; \Cref{sect:univ-LIS}.} Let $\mu$ be a permuton and consider the LIS of  its permuton sample $\LIS{\sigma_n^\mu}$. In \cite{adhikari2025longest,bgsp-26-LIS-sbp} (see Theorems~\ref{thm:brownian-separable-exact-lis} and~\ref{thm:skew-exact-lis2} here), for some fractal random permuton $\mu$, an exponent $\alpha=\alpha(\mu)$  has been identified such that $\LIS{\sigma_n^\mu}=n^{\alpha+o(1)}$.  In these cases, this exponent $\alpha$ is expected to be the a.s.~dimension of an increasing subset in the support of $\mu$ (the support itself has Hausdorff dimension 1 in this case).

Prove that it is indeed the case and find good conditions so that $\LIS{\sigma^\mu_n}=n^{\frac{\alpha}{d}+o(1)}$, where $\alpha$ is the Hausdorff dimension of the largest increasing subset of the support of $\mu$ and $d$ is the Hausdorff dimension of the support of $\mu$. 

\item\label{quest:lis-exp} \textbf{Transferring Brownian separable LIS asymptotics to discrete models; \Cref{sect:univ-LIS}.}
Adhikari, Borga, Budzinski, Da Silva, and Sénizergues identify the exact exponent and the almost-sure scaling limit for permuton samples of the Brownian separable permuton (recall Theorem~\ref{thm:brownian-separable-exact-lis}). They also formulated a transfer conjecture for the discrete model \cite[Conjecture~1.4]{adhikari2025longest}.

Let $\sigma_n$ be a uniform separable permutation, or more generally a uniform permutation in a class whose permuton limit is the Brownian separable permuton $\bm{\mu}_p$. Prove that
\[
\frac{\log\LIS{\sigma_n}}{\log n} \xrightarrow[n\to\infty]{\mathbb{P}} \alpha(p),
\]
where $\alpha(p)$ is the exponent for samples from $\bm{\mu}_p$ appearing in \Cref{thm:brownian-separable-exact-lis}. Can the almost-sure scaling limit for samples from $\bm{\mu}_p$ be transferred to the corresponding discrete models?
The same question can be asked for skew Brownian permutons, Theorem~\ref{thm:skew-exact-lis3} being a first result in this direction.
This kind of question is in the spirit of Section~\ref{sect:lis-permutons},
where we study LIS of random permutation models by comparison with samples of their limiting permutons.
\end{enumerate}

\section*{Acknowledgements}

We warmly thank all those who have contributed to the study of permutons, whose collective work forms the foundation of this survey. {A particular thanks to our collaborators for many years of stimulating work on the topic. We would like to thank also Duncan Dauvergne, Colin Defant, Persi Diaconis, Victor Dubach, Alejandro Mojales,  Sumit Mukherjee, and Leo Petrov for sending feedback on a first version of this survey; with a special nod to Persi for a wonderful 3-page long email with many kind words and even more (legitimate) complaints.}

{The simulation for Grothendieck permutations in the first panel in Figure~\ref{fig:fig-intro-2} was obtained using the code available at \url{https://lpetrov.cc/simulations/2025-01-26-grothendieck-shenanigans/}.}

{J.B.\ is partially supported by NSF grant DMS-2441646 and an Alfred P.\ Sloan Research Fellowship.} {V.F. is partially supported by ANR grant LOUCCOUM, number ANR-24-CE40-7809.}

\bigskip

\noindent\textbf{AI usage statement.} An AI tool was used for literature search, for developing the code for simulations, to improve the wording of a few narrative paragraphs and during the revision of this manuscript to check for possible errors, typographical mistakes, and language issues. All suggested changes were reviewed by the authors, who take full responsibility for the final content of the paper. The exposition, synthesis of the literature, and overall perspective on the field of permutons are the authors’ own; AI was not used to develop these or any original mathematical arguments or proofs presented in this survey.

\newpage

\appendix

\section{Generating permuton pictures}\label{sect-pictures}

\noindent{\emph{Note.} The code discussed in this section is available (in Mathematica and Python) at the following link: \url{https://www.jacopoborga.com/2026/10/12/generating-permuton-pictures/}.}

\bigskip

\begin{figure}[t]
	\centering
	\includegraphics{PermutonPlot-Luce}
	\hspace{1cm}
	\includegraphics[width=0.37\textwidth]{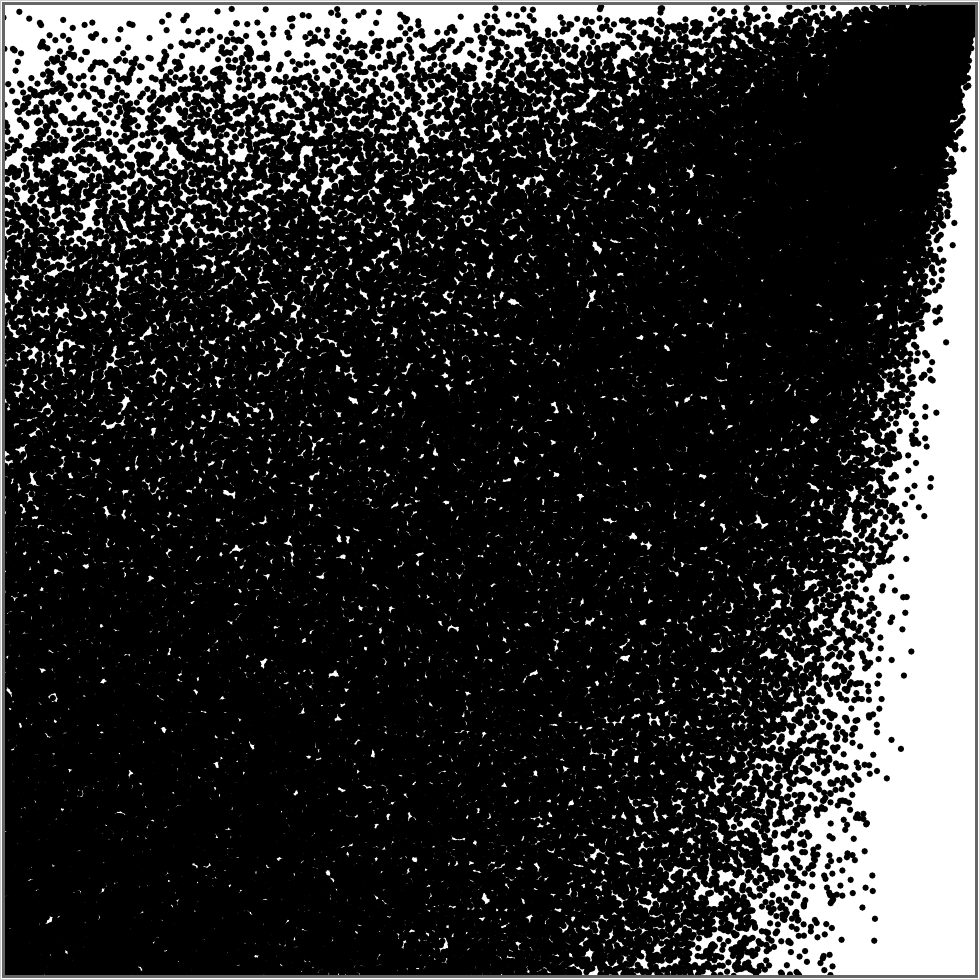}
	\caption{{Two visualizations of the same realization of a Luce-distributed permutation $\pi$ of size $n=100,000$ (recall \Cref{fig:mallows_luce}), both with target image size $300$. \textbf{Left:} The output of \texttt{\detokenize{PermutonPlot[perm, 300, 2]}}. \textbf{Right:} The standard point diagram produced by \texttt{\detokenize{ListPlot[perm, AspectRatio -> 1,  PlotRange -> {{1, Length[perm]}, {1, Length[perm]}}, Axes -> False,  Frame -> True, FrameTicks -> None, PlotStyle -> Black,  ImageSize  -> 300]}}. At this fixed resolution, the markers in the \texttt{ListPlot} representation have a non-negligible pixel footprint and overlap heavily, causing large regions to become saturated in black and obscuring variations in local density. In contrast, \texttt{PermutonPlot} aggregates points at the output resolution and encodes local permuton mass by shades of gray, thereby retaining these density variations and displaying the macroscopic permuton shape more clearly.} \label{fig:comp-plots}}
\end{figure}

We discuss here how the permuton pictures in this paper have been generated. For a permutation $\pi$ of size $n$, the associated permuton $\mu_{\pi}$ assigns Lebesgue measure of mass $1/n$ to $n$ boxes of size $1/n\times1/n$. The difficulty is that, for a very large permutation, the individual squares representing the permutation entries may be much smaller than a single pixel in the final picture.  It is therefore impossible to display the finite permuton literally at screen resolution, and a naive plot might create an accumulation of too large black dots, saturating the picture; see, e.g.,~\Cref{fig:comp-plots}.

Instead, we first replace the $n\times n$ permutation matrix by a coarser $b\times b$ grid, where
\[
	b=\min(n,R),
\]
and $R$ denotes the resolution of the final image. In each grid square we record how many permutation points fall inside it. When $n\le R$, no information is lost and this simply reproduces the usual permutation diagram. When $n>R$, several permutation points may contribute to the same visible region, so these counts provide a natural approximation to the local mass of the permuton.

To make the large-$n$ picture easier to read, we then average these counts slightly over neighboring grid squares using a Gaussian weight. {The appropriate amount of averaging depends not only on $n$ and the output resolution, but also on the expected dimension of the support of the limiting permuton. For a permuton with a two-dimensional density, local averaging is useful for suppressing sampling noise and revealing variations in the density. For a permuton concentrated on a lower-dimensional set, however, too much averaging can blur the geometric structure that one wants to display.}

{The plotting function therefore takes an optional parameter $d\in[1,2]$, interpreted as the expected Hausdorff dimension of the support of the limiting candidate permuton. The Gaussian smoothing scale is chosen heuristically according to}
\[
	\sigma_d = \left( \frac{0.6\,b^d}{(2\pi)^{d/2}n} \right)^{1/d},
\]
{together with an additional factor which makes the smoothing vanish when the permutation is already well resolved at the chosen output resolution, and with an upper cutoff to avoid excessive blurring.} Near the boundary of the square, this averaging is corrected so that points close to an edge are not artificially given less weight.

{The scaling above is motivated by the following heuristic. If the limiting measure is concentrated on a set of dimension $d$, then a neighborhood of pixel radius $\sigma$ should contain roughly a proportion $(\sigma/b)^d$ of the total mass, and hence about}
\[
	n\left(\frac{\sigma}{b}\right)^d
\]
{permutation points. Keeping this quantity of constant order leads to the scaling $\sigma\asymp b n^{-1/d}$. The precise numerical constants used in the code are only visualization choices.}

Finally, regions carrying more permutation mass are represented by darker shades of gray. The conversion from local mass to darkness is chosen so that a few exceptionally dense regions do not force the rest of the picture to appear almost white. The final image is therefore a resolution-dependent approximation of the permuton in which the darkness of each region reflects the amount of permutation mass concentrated there.

{For each fixed output resolution, the normalized coarse-grid mass array is compatible with permuton convergence: if $\mu_{\pi_m}$ converges weakly to $\mu$, then the mass recorded in each grid cell converges to the mass assigned to that cell by $\mu$. The subsequent Gaussian averaging and grayscale transformation are visualization devices applied to this coarse-grained mass array.}

{The code (see the note above for the link) implements this idea both in Mathematica and Python.} For example, use \texttt{PermutonPlot[perm]} or \texttt{PermutonPlot[perm, 1000]} for a higher-resolution output {(the first optional argument is the desired image size in pixels, while a second optional argument specifies the expected Hausdorff dimension $d$ of the support of the limiting candidate permuton). The default value is $d=2$. Thus, for example, one may use \texttt{PermutonPlot[perm, 300, 1]} for a candidate whose limiting support is expected to have Hausdorff dimension $1$.}

\bigskip

{\noindent\textbf{Adaptive smoothing for permutons with unknown or mixed-dimensional limits.}
The preceding choice of a single dimension parameter $d$ is less suitable when the dimension of the limiting support of the permuton is unknown, or when the limiting permuton has both an absolutely continuous and a singular part, as in the Grothendieck, record-biased, or runsort examples (see Figure~\ref{fig:fig-intro-2}). In these situations, the function \texttt{PermutonPlotAdaptive} chooses the smoothing scale from local point spacings rather than from a prescribed dimension. 

For each permutation point, it measures the distance to its $k$-th nearest other point (with $k=32$ by default) and uses one third of this distance as the Gaussian smoothing scale, subject to a maximum of two pixels. As before, smoothing is reduced when $n$ is close to the output resolution and is omitted when $n\le R$. This tends to produce narrow kernels on densely sampled curves and wider kernels in diffuse regions.}

\newpage

\bibliographystyle{bibli_perso}
\bibliography{Bib_jb,Bib_additional}

\end{document}